\documentclass[10pt]{article}
\usepackage[T1]{fontenc}
\usepackage[utf8]{inputenc}
\usepackage{float}
\usepackage{amsmath,amsthm,graphicx,mathtools}
\usepackage{subcaption}
\usepackage{amssymb}
\usepackage{tikz}
\usetikzlibrary{patterns,patterns.meta,arrows.meta,calc,3d,perspective,fadings,shadows,cd}
\usepackage[a4paper, width=150mm, top=30mm, bottom=30mm]{geometry}
\usepackage{hyperref}

\theoremstyle{plain}
	\newtheorem{thm}{Theorem}[section]
	\newtheorem{cor}[thm]{Corollary}
	\newtheorem{lem}[thm]{Lemma}
	\newtheorem{prop}[thm]{Proposition}
	\newtheorem*{thm:gen_Haas}{Theorem~\ref{thm:gen_Haas}}
	\newtheorem*{prop:simple_harnack}{Proposition~\ref{prop:simple_harnack}}
	\newtheorem*{prop:sharper_upper_bound}{Proposition~\ref{prop:sharper_upper_bound}}
	\newtheorem*{thm*}{Theorem}
	\newtheorem*{thm:ItViro}{Theorem~\ref{thm:ItViro}}
	\newtheorem*{prop:partial_eps_as_wave}{Proposition~\ref{prop:partial_eps_as_wave}}
	\newtheorem*{thm:growth_rate}{Theorem~\ref{thm:growth_rate}}

\theoremstyle{definition}
    
	\newtheorem{dfn}[thm]{Definition}

	\newtheorem{conv}[thm]{Conventions}
	\newtheorem{rem}[thm]{Remark}
	
	\newtheorem{ex}[thm]{Example}

\newcommand{\C}{\mathbb{C}}
\newcommand{\R}{\mathbb{R}}
\newcommand{\Q}{\mathbb{Q}}
\newcommand{\Z}{\mathbb{Z}}
\newcommand{\T}{\mathbb{T}}
\newcommand{\F}{\mathbb{F}}
\renewcommand{\P}{\textnormal{P}}
\newcommand{\Proj}{\mathbb{P}}
\newcommand{\Sed}{\mathcal{S}\hspace{-.4ex}\textit{ed}}

\newcommand{\relint}{\textnormal{rel.int.}}
\newcommand{\coker}{\textnormal{coker}}
\newcommand{\vol}{\textnormal{vol}}

\newcommand{\card}{\textnormal{card}}
\newcommand{\Hom}{\textnormal{Hom}}

\renewcommand{\O}{\mathcal{O}}
\newcommand{\K}{\mathcal{K}}
\newcommand{\df}{\mathrm{d}}
\newcommand{\D}{\nabla}

\newcommand{\X}{{\R X_\varepsilon}}

\newcommand{\x}{\textnormal{x}}

\renewcommand{\mod}{\textnormal{mod}\,}

\newcommand{\Arg}{\mathcal{E}}
\newcommand{\bv}{\textnormal{bv}}
\newcommand{\Aff}{\textnormal{Aff}}

\newcommand{\rk}{\textnormal{rk}\,}
\newcommand{\Rad}{\textnormal{Rad}}
\newcommand{\rad}{\textnormal{rad}}

\newcommand{\oo}{o}

\newcommand{\dv}{\mathrm{v}}

\newcommand{\cW}{\mathcal{W}}

\newcommand{\cZ}{\mathcal{Z}}
\newcommand{\cF}{\mathcal{F}}
\newcommand{\la}{\langle}
\newcommand{\ra}{\rangle}

\newcommand{\M}[1]{M_{#1}}
\newcommand{\N}[1]{N_{#1}}
\newcommand{\capp}{{\scriptstyle{\,\cap\,}}}

\newcommand{\cupp}{{\scriptstyle{\,\cup\,}}}

\title{\bf On the First Page of the Renaudineau--Shaw Spectral Sequence}
\author{Jules Chenal\footnote{Universitetet i Oslo, Norge, julesche[at]math.uio.no}}
\date{September 2026}

\begin{document}

\maketitle

\begin{abstract}
	In this article, we study the topology of \emph{T-hypersurfaces}. These cellular hypersurfaces are generalisations of the cellular complexes used to describe the real loci of varieties obtained via \emph{Viro's primitive patchworking method}. As such, they have deep ties to real geometry. Renaudineau and Shaw introduced a spectral sequence computing their homology from which they derived upper bounds on their Betti numbers using tropical geometry. Here, we study the first page this spectral sequence and express the action of its boundary operators on a distinguished subspace as cap-products with \emph{Mikhalkin--Zharkov waves}. Then, we characterise those T-hypersurfaces with a maximal number of connected components with respect to the Renaudineau--Shaw inequality with a combinatorial condition on the building triangulation and a system of combinatorial differential equations on the sign distribution. This extends a theorem of Haas for curves. In addition, we study the growth rate of the expected number of connected components and provide examples of T-hypersurfaces of the projective spaces satisfying these conditions in every degree and dimension. 
\end{abstract}

\tableofcontents

\paragraph{Acknowledgements}
	This research was funded, in part, by l’Agence Nationale de la Recherche (ANR), project ANR-22-CE40-0014. The author was also supported by the Trond Mohn Foundation project “\,Algebraic and Topological Cycles in Complex and Tropical Geometries\,”. For the purpose of open access, the author has applied a CC-BY public copyright licence to any Author Accepted Manuscript (AAM) version arising from this submission.

\section{Introduction}

One of the main themes of real algebraic geometry is the study of the topology of the real loci of real algebraic varieties. A typical question of the field, tracing back to Hilbert's 16\textsuperscript{th} problem, is to classify the possible embedded topologies of algebraic hypersurfaces of the real projective space $\Proj^{n+1}$ with a fixed degree $d\geq 0$. Such a hypersurface $\mathcal{X}$ is defined by the vanishing of a real homogeneous polynomial $f\in \R[x_0,...,x_{n+1}]$ of degree $d$. This is a very ambitious problem and one would be happy to simply derive algebraic constraints on the topology of the real solution set $\mathcal{X}(\R)$: its \emph{real locus}. In contrast, the possible complex solution sets $\mathcal{X}(\C)$ are all homeomorphic provided $\mathcal{X}$ is non-singular. This follows Thom Isotopy Lemma. Perhaps the most general topological constraint on the real loci of algebraic varieties is the Smith--Thom Inequality, cf. \cite[Theorem 2.3.1]{degtyarev_topological_2000},
\begin{equation}\label{eq:SM}
	\sum_{q\geq 0} \dim H^q(\mathcal{X}(\R);\F_2) \leq \sum_{q\geq 0} \dim H^q(\mathcal{X}(\C);\F_2).
\end{equation} 
It applies to every real variety, not only to hypersurfaces of the projective spaces. Those varieties for which (\ref{eq:SM}) is an equality are called \emph{maximal varieties} and are of great interest to real geometers as the available examples are rather scarce. Let us assume forward that $\mathcal{X}$ is non-singular. Therefore, the right hand side of (\ref{eq:SM}) only depends on the degree $d$. In this framework, the maximality problem asks whether there exist maximal hypersurfaces of every degree in every projective space. 

This illustrates the dichotomic approach often used in the study of the topology of the real loci of real algebraic varieties. First establishing constraints, then testing their optimality by constructing examples. To this day, the most powerful construction technic remains \emph{Viro's patchworking method}. In its most general form, Viro's theorem provides a topological description of the real locus of certain hypersurfaces of a toric variety defined by a section of the line bundle spanned by a moment polytope of the toric variety. We refer to \cite{viro_patchworking_2006, risler_construction_1993} and \cite[Chapter 11, \S5]{gelfand_discriminants_1994} for in depth expositions of the theorem. 

\begin{figure}[H]
	\centering
	\begin{subfigure}[t]{.5\textwidth}
		\centering
		\input{tropical_subdivision.tex}
		\caption{The subdivision of $\mathbb{P}^2(\T)$ induced by the tropical hypersurface defined by the tropical Laurent polynomial $f=\max(0,x+2,y+2,2x+2,x+y+3,2y+2,3x,2x+y+2,x+2y+2,3y)$}
		\label{fig:trop_cub}
	\end{subfigure}
	\hfill
	\begin{subfigure}[t]{.45\textwidth}
		\centering
		\input{harnack_dist.tex}
		\caption{The sign distribution $\varepsilon$.}
	\end{subfigure}
	\hfill
	\begin{subfigure}[t]{.45\textwidth}
		\centering
		\input{tropical_patchwork.tex}
		\caption{The curve $V_\varepsilon\subset \mathbb{P}^2(\R)$ associated to the tropical vanishing locus of $f$ and $\varepsilon$ (n.b. the boundary edges have to be glued in an antipodal way).}
		\label{fig:patch_V}
	\end{subfigure}
	\hfill
	\begin{subfigure}[t]{.5\textwidth}
		\centering
		\includegraphics[width=.75\textwidth]{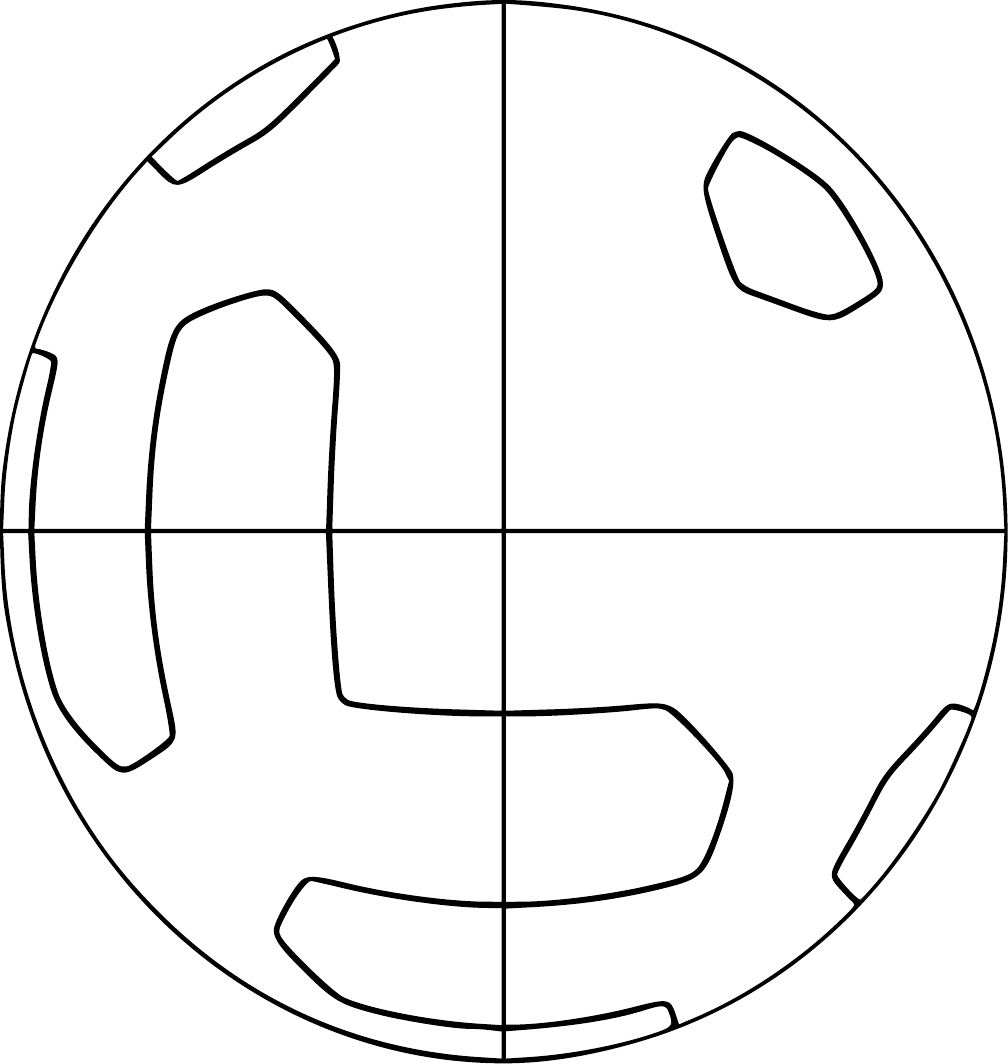}
		\caption{The real locus of the plane cubic of equation $1+x^3+y^3+t^2(x+y+x^2+y^2+x^2y+xy^2)-t^3xy=0$ with $t=e^{29/2}$ provided by Viro's Theorem in the parametrisation $\{x^2+y^2\leq 1\}\rightarrow \mathbb{P}^2(\R)$ given by $(x,y)\mapsto [x^{29}:y^{29}:(1-x^2-y^2)^{29/2}] $.}
	\end{subfigure}
	\caption{The patchwork of a tropical cubic.}
	\label{fig:patch_cub}
\end{figure}

In the present article, we will only be concerned with a sub-method called \emph{primitive patchworking}, especially its more recent formulation within the scope of \emph{tropical geometry}. This theory is an algebraic geometry over the \emph{tropical semi-field} $\T=\{-\infty\}\coprod \R$ where the tropical addition is given by the maximum and the tropical multiplication by the classical addition. A tropical polynomial is then a convex piecewise affine function with integral slope. Every projective toric variety $\mathcal{Y}$ naturally possesses a \emph{tropical locus} $\mathcal{Y}(\T)$ which is, in particular, a cellular compactification of the tropical torus $(\T^\times)^{n+1}=\R^{n+1}$. With the help of the moment map, one shows that $\mathcal{Y}(\T)$ is isomorphic to $P$ as a CW-complex. As in classical geometry, a tropical hypersurface $\mathcal{X}(\T)$ of $\mathcal{Y}(\T)$ is defined by a single equation “$f=0\,$” where $f$ is a tropical Laurent polynomial. The closures of the connected components of $\mathcal{Y}(\T)\setminus\mathcal{X}(\T)$ naturally endow $\mathcal{Y}(\T)$ with a cellular structure of which $\mathcal{X}(\T)$ can be seen as the support of a subcomplex; cf. Figure~\ref{fig:trop_cub}. We refer to \cite{maclagan_introduction_2015, brugalle_brief_2015, mikhalkin_tropical_nodate} for references on tropical geometry.

In this context, primitive patchworking starts with a lattice polytope $P\subset \R^{n+1}$ yielding a smooth toric variety $\mathcal{Y}$, and a non-singular tropical hypersurface $\mathcal{X}(\T)\subset Y(\T)$ whose defining equation has Newton polytope $P$. Then, to every \emph{sign distribution} $\varepsilon : P\cap \Z^{n+1}\rightarrow \F_2$, one can associate a PL-smooth hypersurface $V_\varepsilon\subset \mathcal{Y}(\R)$ endowed with a natural cellular surjection $V_\varepsilon\rightarrow \mathcal{X}(\T)$. To do so, we first notice that $\mathcal{Y}(\R)$ is homeomorphic to the gluing of $2^{n+1}$ copies of $\mathcal{Y}(\T)$ indexed by the element of $\N{2}$, cf. the application $\phi$ of (\ref{eq:phi}), and that the $(n+1)$-cells of $\mathcal{Y}(\T)\setminus \mathcal{X}(\T)$ are in natural bijection with the lattice points of $P$. Then, in the copy of $\mathcal{Y}(\T)$ indexed by $u\in \N{2}$, $V_\varepsilon$ is the union of the $n$-cells of $\mathcal{X}(\T)$ bounding two $(n+1)$-cells whose corresponding vertices $v_1,v_2\in P\cap M$ satisfy $\varepsilon(v_1)+\la v_1,u\ra \neq \varepsilon(v_2)+\la v_2,u\ra$; see Figure~\ref{fig:patch_V}. Viro's theorem guaranties that there exists a real Laurent polynomial with Newton polytope $P$ that defines an algebraic hypersurface of $\mathcal{Y}$ whose real locus is ambiantly isotopic to $V_\varepsilon$. See Figure~\ref{fig:patch_cub}.  

	Moreover, the equation of $\mathcal{X}(\T)$ also yields a \emph{primitive triangulation}\footnote{A lattice triangulation whose maximal simplices have minimal lattice normalised volume $1/(n+1)!$. In particular, the vertices of $K$ are exactly the lattice points of $P$.} $K$ of $P$, called its \emph{dual subdivision}. A primitive triangulation $K$ of $P$ is dual to a tropical hypersurface if and only if it is \emph{convex}; see Definition~\ref{dfn:prim_conv}. For any primitive triangulation $K$ of $P$, we can use the duality property to build a cellular subdivision $Y$ of $P$ and a subcomplex $X\leq Y$ which are isomorphic to $(\mathcal{Y}(\T),\mathcal{X}(\T))$ in the convex case; see Definitions~\ref{dfn:XY} and \ref{dfn:X}. This allows us to generalise the construction of $V_\varepsilon$ to non-convex triangulations. Following Definition~\ref{dfn:t-hyp}, to every sign distribution $\varepsilon : P\cap \Z^{n+1}\rightarrow \F_2$, we associate a PL-smooth hypersurface $\R X_\varepsilon \subset \mathcal{Y}(\R)$ endowed with a cellular surjection $\X\rightarrow X$. Such a hypersurface $\X$ is called a \emph{T-hypersurface}. We also note that the cellular complex $Y$ induces a cellular structure $\R Y$ on $\mathcal{Y}(\R)$ of which $\X$ is a subcomplex. Even though we loose ties to both tropical and real geometry when $K$ is not convex, T-hypersurfaces enjoy most of the properties of the hypersurfaces $V_\varepsilon$ such as the Smith--Thom inequality (\ref{eq:SM}). Moreover, the theory of tropical homology has a natural extension to the complexes $X$ which satisfies all main structure theorems. 

In \cite{renaudineau_bounding_2023}, Renaudineau and Shaw used the tropical framework to construct a spectral sequence $(E^r_{p,q}(\X))_{p,q,r\geq 0}$ computing the singular homology of $\X$ with coefficients in $\F_2$. Their construction was extended to the non-convex setting in \cite{brugalle_combinatorial_2024,chenal_poincare_2026}. They showed that its first page is isomorphic to the \emph{tropical homology} of $X$. The tropical homology groups $(H_{p,q}(X;\F_2))_{p,q\geq 0}$, introduced in \cite{itenberg_tropical_2019} and extended beyond convexity in \cite{brugalle_combinatorial_2024,chenal_poincare-lefschetz_2025}, mimic the Hodge theory of complex varieties. They enjoy many of the properties of the classical Hodge groups and it can be showed that, in our case, $\dim H_{p,q}(X;\F_2)$ equals the corresponding Hodge number of a generic complex hypersurface of $\mathcal{Y}$ with Newton polytope $P$. We provide an account of the relevant properties of tropical homology in Section~\ref{sec:trop_homo}. The Renaudineau--Shaw spectral sequence then implies the \emph{Renaudineau--Shaw Inequalities}
\begin{equation}\label{eq:RSI}
	\dim H_q(\X;\F_2) \leq \sum_{p\geq 0} \dim H_{p,q}(X;\F_2) \quad\textnormal{for all $q\geq 0$},
\end{equation} 
which refine (\ref{eq:SM}) for \emph{real hypersurfaces near the tropical limit}, i.e. obtained via primitive patchwork; cf. \cite[Theorem 1.4]{renaudineau_bounding_2023}. The construction of the Renaudineau--Shaw spectral sequence is detailed in Section~\ref{sec:RSSS}.

Here, we investigate deeper the first page of the spectral sequence. As Renaudineau and Shaw pointed out in \cite{renaudineau_bounding_2023}, there is only one chain complex with \emph{a priori} non-trivial boundary operators in $(E^1_{p,q}(\X))_{p,q\geq0}$. It is located on the line $p+q=n$. Using the Borel--Viro isomorphism, \cite[Definition 4.9]{renaudineau_bounding_2023}, we obtain boundary operators
\begin{equation*}
	\partial_\varepsilon: H_{p,q}(X;\F_2)\longrightarrow H_{p+1,q-1}(X;\F_2) \quad\textnormal{for all $p+q=n$}.
\end{equation*}
Our first result is to express the restriction of those operators to the subspaces spanned by the so-called \emph{face cycles} as cap-products with \emph{wave classes}.

Every interior simplex $\sigma^p$ of the dual triangulation $K$ yields a face cycle $f_{\sigma^p}$ and a class $[f_{\sigma^p}]\in H_{p,n-p}(X;\F_2)$; cf. Definition~\ref{dfn:face_cycle}. Moreover, the collection of these classes, for $p$ fixed, spans the kernel of the morphism $H_{p,n-p}(X;\F_2)\rightarrow H_{p,n-p}(Y;\F_2)$ induced by the inclusion $X\leq Y$; see Proposition~\ref{prop:face_cycles_kernel}. These morphisms are all surjective, this follows the Tropical Lefschetz Hyperplane Section Theorem; cf. \cite[Theorem 1.2]{arnal_lefschetz_2021}, \cite[Proposition 3.2]{brugalle_combinatorial_2024}, or \cite[Theorem 2]{chenal_poincare-lefschetz_2025}. In addition, the group $H_{p,n-p}(Y;\F_2)$ vanishes unless $n=2p$.

The wave cohomology, introduced in \cite{mikhalkin_tropical_2014}, associates a bigraded ring $(H^q(X;\cW_p))_{p,q\geq 0}$ to the complex $X$. The cap-product of \cite[Proposition 9]{mikhalkin_tropical_2014} turns $(H_{p,q}(X;\F_2))_{p,q\geq 0}$ into a module over $(H^q(X;\cW_p))_{p,q\geq 0}$. For instance, the cap-product with a wave class $w\in H^1(X;\cW_1)$ acts similarly to one of the boundary operators $\partial_\varepsilon$
\begin{equation*}
	[c]\in  H_{p,q}(X;\F_2)\mapsto w\capp[c]\in H_{p+1,q-1}(X;\F_2) \quad\textnormal{for all $p,q$}.
\end{equation*}
These constructions are recalled in Section~\ref{sec:waves}. We prove that if $e^p$ is a \emph{bounded cell} of $X$, see Definition~\ref{dfn:XY}, then $H^p(e^p;\cW_p)$ has dimension $1$; cf. Lemma~\ref{lem:cohomo_wave_bounded_cell}. As every $1$-dimensional $\F_2$ vector space is canonically isomorphic to $\F_2$, we define, for every $w\in H^p(X;\cW_p)$, a number $w|_{e^p}^{\F_2}\in \F_2$ that vanishes if and only if $w|_{e^p}=0$; cf. Definition~\ref{dfn:wave_numbers}. 

To every sign distribution $\varepsilon : P\cap \Z^{n+1}\rightarrow \F_2$, we associate a \emph{Hessian wave} $\D^2\varepsilon\in H^1(X;\cW_1)$; cf. Definition~\ref{dfn:hessian_wave}. It measures where $\varepsilon$ fails to be an affine function of the coordinates of the lattice points of $P$. If we denote by $\chi_{\sigma^p}:P\cap \Z^{n+1}\rightarrow \F_2$ the function with value $1$ exactly on the vertices of $\sigma^p$, we show that:
\begin{prop:partial_eps_as_wave}
	For all interior simplices $\sigma^p$, we have $\partial_\varepsilon [f_{\sigma^p}]=\D^2(\varepsilon+\chi_{\sigma^p})\capp[f_{\sigma^p}]$.
\end{prop:partial_eps_as_wave}
This allows us to characterise the couples $(K,\varepsilon)$ for which the inequality (\ref{eq:RSI}) with $q=0$ is an equality. In this case, we say that $\X$ has a maximal number of connected components with respect to the Renaudineau--Shaw inequality. We let $X_0\leq X$ denote the subcomplex spanned by the support of the face cycles $f_v$ for all interior vertices $v$ of $K$. The complex $X_0$ contains another distinguished subcomplex $X_1$ whose precise definition can be found in Definition~\ref{dfn:X0}. Moreover, we say that a triangulation $K$ is \emph{$\rho$-uniform}, see Definition~\ref{dfn:rho_unif}, if there is a, necessarily unique, wave $\rho\in H^1(X_0;\cW_1)$ whose restriction to the support of $f_v$ equals $\D^2\chi_v|_{\mathrm{Supp}(f_v)}$ for every interior vertex $v$. We provide a geometrical characterisation of this property in Proposition~\ref{prop:crit_rho_unif}. With these definitions, we can state our main result:

\begin{thm:gen_Haas}
	Let $n\geq 2$. The hypersurface $\X$ has a maximal number of connected components with respect to the Renaudineau--Shaw inequality if and only if \textnormal{(1)} $K$ is $\rho$-uniform; \textnormal{(2)} the morphism $H_1(\X;\F_2)\rightarrow H_1(\R Y;\F_2)$ induced by the inclusion $\X\subset \R Y$ is onto; and \textnormal{(3)} the sign distribution $\varepsilon$ satisfies the equations
	\begin{equation*}\tag{\ref{eq:equation_max_numb_comp}}
		\D^2\varepsilon|_{X_1}=\rho|_{X_1}\quad \mathrm{and}\quad \D^2\varepsilon|_{X_0}\cupp\rho=\rho^2.
	\end{equation*}
\end{thm:gen_Haas}

It extends Haas' Theorem for curves to higher dimensions; cf. \cite[Theorem 7.3.0.10]{haas_real_1997} or \cite[Theorem 4.6]{bertrand_haas_2017}. However, some translation is required to clearly see the parallel with Haas' Theorem. 

First, if $n=1$, $K$ is automatically $\rho$-uniform. Moreover, the second condition is not required for curves (and in fact usually not satisfied, even for maximal curves). It accounts for higher dimensional phenomena: when $n=2$ it counterbalances the fact that $1$-dimensional face cycles form a proper subspace of $H_{1,1}(X;\F_2)$, while in higher dimensions, it ensures the vanishing of the adequate higher boundary operators of the spectral sequence. The latter vanishing criterion was actually the topic of our article \cite{chenal_poincare_2026} and is conjecturally always satisfied.

Furthermore, Haas' Theorem was stated in the vocabulary of twisted edges, a notion originated from the shape of the boundary of the \emph{amoeba} of the algebraic curve constructed by Viro. Haas' original formulation stated that:

\begin{thm*} \textnormal{\cite[Theorem 7.3.0.10]{haas_real_1997}, \cite[Theorem 4.6]{bertrand_haas_2017}, or Theorem~\ref{thm:haas_classic}}.
	A curve $\X$ is maximal (i.e. has a maximal number of connected components) if and only if \textnormal{(1)} none of the edges of $X_1$ are twisted; and \textnormal{(2)} every cycle of $X$ contains an even number of twisted edges.
\end{thm*}

We discuss the extension of the concept of twisted edges to $n\geq 2$ in Section~\ref{sec:twists}. Note that for $n\geq 2$, it is, in full generality, more subtle than for curves. Nonetheless, given a curve $\X$, a bounded edge $e^1$ of $X$ is twisted if and only if $(\D^2\varepsilon+\rho)|^{\F_2}_{e^1}=1$. Therefore, the first equation of Theorem~\ref{thm:gen_Haas} extends the first condition of Haas' Theorem.

The second equation, however, does not \emph{a priori} make sense for curves as it happens in $H^2(X_0;\cW_2)$, which vanishes for dimensional reasons when $n=1$. However, if $n\geq 2$ and $e^2\subset X_0$, the number $(\D^2\varepsilon|_{X_0}+\rho)\cupp\rho|^{\F_2}_{e^2}$ is given by a weighted sum of the numbers $(\D^2\varepsilon+\rho)|^{\F_2}_{e^1}$ over the edges $e^1$ of boundary of $e^2$. Therefore, $\D^2\varepsilon|_{X_0}\cupp\rho=\rho^2$ can be interpreted as asking that some weighted numbers of twists along the boundaries of $2$-cells of $X_0$ are even, thus corresponding to Haas' second condition for in dimension $n=1$ every cycle of $X$ is contained in $X_0$. We can provide computations of the weights, however, they will be non-canonical and depend on some choices (but the weighted sums themselves will not). This is why we prefer to keep this cohomological formulation of the equation; see Remark~\ref{rem:haas_twists}. 

Solving the equations of Theorem~\ref{thm:gen_Haas} is not necessarily an easy task. However, we can consider the stronger equation
\begin{equation*}\tag{\ref{eq:simple_harnack}}
	\D^2\varepsilon|_{X_0}=\rho|_{X_0},
\end{equation*}
which, of course, only makes sense for $\rho$-uniform triangulations. It obviously implies both equations of Theorem~\ref{thm:gen_Haas} but also implies the surjectivity condition. In addition, (\ref{eq:simple_harnack}) allows us to derive further topological information about the T-hypersurfaces constructed by its solutions.
\begin{prop:simple_harnack}
	If $\varepsilon$ is a solution Equation~(\ref{eq:simple_harnack}), then $\X$ is made of \textnormal{(1)} $h^{n,0}(X)$ spheres contained in the torus of $\mathcal{Y}(\R)$ where they bound disjoint balls; and \textnormal{(2)} an additional connected component that intersects every toric divisor of $\mathcal{Y}(\R)$. When (\ref{eq:simple_harnack}) has a solution, we say that $K$ is \emph{simply integrable}. 
\end{prop:simple_harnack}
The advantage of (\ref{eq:simple_harnack}) over (\ref{eq:equation_max_numb_comp}) is our ability to provide a cohomological criterion for its resolution; see Proposition~\ref{prop:solvability_simple_harnack}. In the case of curves, the \emph{Harnack distribution} given in \cite[\S Patchworking Harnack Curves]{itenberg_patchworking_1996} is a universal solution of (\ref{eq:simple_harnack}). It is actually a solution of an even stronger equation; see Theorem~\ref{thm:quadratic_distrib}. In addition, we quantify the defect of $\rho$-uniformity of a triangulation $K$ with a natural number $\kappa(K)$. The latter integer can be used to sharpen the Renaudineau--Shaw inequality: 
\begin{prop:sharper_upper_bound}
	For all $\varepsilon$, we have  $\dim H_0(\X;\F_2)\leq 1+\dim H_{0,n}(X;\F_2) - \kappa(K)$.
\end{prop:sharper_upper_bound}
We also use our computations to determine the growth rate of the expected number of connected components of random T-hypersurfaces $(\X(d))_{d\geq 1}$ of $\mathcal{Y}(\R)$ with respective Newton polytopes $dP$. The explicit description of the stochastic process can be found in Section~\ref{sec:growth}.
\begin{thm:growth_rate}
	In Landau's notation, $\mathbb{E}[b_0(\X(d))]=\Theta(d^{n+1})$, i.e. it is eventually bounded above and below by positive multiples of $d^{n+1}$.
\end{thm:growth_rate}

This result is in strong contrast with the behaviour of random algebraic hypersurfaces following Kostlan's distributions. For instance, when one chooses a random real hypersurface of $\mathcal{Y}$ in the linear system $|\mathcal{O}(dP)|$ with Kostlan's measure, Gayet and Welschinger showed that the growth rate of every Betti number of the real locus is $\Theta(d^{(n+1)/2})$; cf. \cite{gayet_lower_2014,gayet_betti_2016}. From this point of view, typical T-hypersurfaces differ greatly from typical real hypersurfaces. This observation is strengthened by Ancona's result \cite[Theorem 1.2]{ancona_exponential_2022}. The latter states that the Kostlan's measure of the set of hypersurfaces in $|\mathcal{O}(dP)|$ whose total Betti number is of order of $d^{n+1}$ tends to $0$ faster than any negative power of $d$.

Finally, in Section~\ref{sec:examples}, we study $\rho$-uniformity and solvability of (\ref{eq:simple_harnack}) for three families of triangulations of the projective spaces. In particular, we show that the \emph{Itenberg--Viro triangulations}, inspired by \cite[\S 5]{itenberg_topology_1997}, are all $\rho$-uniform, cf. Proposition~\ref{prop:rho_unif_IV}, and that  Equation~(\ref{eq:simple_harnack}) has an explicit solution:

\begin{thm:ItViro}
	For all $n\geq 2$, and all $d\geq 1$, the Itenberg-Viro triangulation $IV^n_d$ is simply integrable. In particular, the distribution
	\begin{equation*}
		h^n(x_1,...,x_n)\coloneqq \sum_{i<j}x_ix_j,
	\end{equation*}
	is a solution of (\ref{eq:simple_harnack}).
\end{thm:ItViro}

To conclude this introduction, we hope that this work, especially the interpretation of the boundary operators of the first page in terms of Mikhalkin--Zharkov waves, will help in the more than twenty years old endeavour of constructing maximal hypersurfaces of the projective spaces of every dimension and every degree using primitive patchwork.    

\section{Preliminaries}

\begin{conv}
	In the whole article, we will consider two $(n+1)$-dimensional lattices $M$ and $N$ endowed with a duality $\langle-,-\rangle :M\otimes N \rightarrow \Z$. We will denote $N\otimes\R$ by $\N{\R}$ and $N\otimes\F_2$ by $\N{2}$. The same notations apply to $M$. Moreover, if $V$ is a finite-dimensional vector space over $\F_2$,
	\begin{enumerate}
		\item $V^*$ denotes its dual $\Hom(V;\F_2)$ and their duality pairing is denoted by $\la-,-\ra$,
		\item the \emph{volume element} of $V$ is the unique non-zero vector of $\bigwedge^dV$ where $d=\dim V$.
	\end{enumerate}
	Finally, if $P$ is a convex polyhedron of $\M{\R}$, $TP$ denotes its tangent space.   
\end{conv}

\subsection{Cellular Complexes and Cellular Homology}

Every finite simplicial complex defines a regular CW-complex called is \emph{geometric realisation}. In this article, we will assimilate these two objects. Let $K$ be a simplicial complex. The \emph{star} of a simplex $\sigma\in K$ is the union of the interiors of the simplices of $K$ containing $\sigma$.
 
If $\sigma\leq \tau$ are two simplices, $\tau-\sigma$ denotes the simplex of $\tau$ spanned by the vertices of $\tau$ not contained in $\sigma$. If $\sigma,\tau$ are two disjoint simplices the \emph{join} $\sigma*\tau$ is the abstract simplex on the vertices of $\sigma$ and $\tau$. We use the notation $\sigma*\tau\in K$ to signify that (1) $\sigma$ and $\tau$ are disjoint; and (2) there is a simplex whose vertex set is the union of those of $\sigma$ and $\tau$. We denote such simplex by $\sigma*\tau$. When we write $\sigma^k\in K$, it means $\sigma^k$ is a simplex of $K$ of dimension $k$ or also a $k$-simplex of $K$. We use a similar convention for the cells of any CW-complex. 

\begin{dfn}[Barycentric Subdivision]
	Let $K$ be a finite simplicial complex. Recall that its \emph{barycentric subdivision} is the simplicial complex whose vertices are the simplices of $K$ and whose simplices are given by chains of simplices $\{\sigma_0<\cdots<\sigma_k\}$. Figure~\ref{fig:bar_sub_tri} depicts the barycentric subdivision of a triangle.
\end{dfn}

\begin{dfn}[Cubical Subdivision]	 
	Let $K$ be a finite simplicial complex. The \emph{cubical subdivision} of $K$ is a regular CW-complex subdivision of K whose open $k$-cells are indexed by pairs of simplices $[\sigma,\tau]$ where $\sigma$ is a $k$-codimensional face of $\tau$, and are defined by
	\begin{equation*}
		[\sigma,\tau]\coloneqq \underset{\hspace{-16ex}\sigma<\sigma_1<\cdots<\sigma_{k-1}<\tau}{\bigcup \hspace{.5ex}\{\sigma<\sigma_1<\cdots<\sigma_{k-1}<\tau\}}.
	\end{equation*} The cubical subdivision is, by definition, coarser than the barycentric subdivision. The subdivision of its closed $k$-cells have the combinatorics of $k$-dimensional cubes triangulated into $k!$ simplices of dimension $k$. The cubical subdivision of a triangle is represented in Figure~\ref{fig:cub_sub_tri}.
\end{dfn}

\begin{figure}[ht]
	\centering
	\begin{subfigure}[t]{0.45\textwidth}
		\centering
		\begin{tikzpicture}[scale=2]
			\draw[thick] (0,0) -- (0,1) -- (1,0) -- cycle ;
			\draw[thick] (.33,.33) -- (0,0) ;
			\draw[thick] (.33,.33) -- (.5,0) ;
			\draw[thick] (.33,.33) -- (1,0) ;
			\draw[thick] (.33,.33) -- (.5,.5) ;
			\draw[thick] (.33,.33) -- (0,1) ;
			\draw[thick] (.33,.33) -- (0,.5) ;
			\fill (0,0) circle (.03);
			\fill (.5,0) circle (.03);
			\fill (1,0) circle (.03);
			\fill (.5,.5) circle (.03);
			\fill (0,1) circle (.03);
			\fill (0,.5) circle (.03);
			\fill (.33,.33) circle (.03);
		\end{tikzpicture}
		\caption{Barycentric subdivision of a triangle.}
		\label{fig:bar_sub_tri}
	\end{subfigure}
	\begin{subfigure}[t]{0.45\textwidth}
		\centering
		\begin{tikzpicture}[scale=2]
			\draw[thick] (0,0) -- (0,1) -- (1,0) -- cycle ;
			\draw[thick] (.33,.33) -- (.5,0) ;
			\draw[thick] (.33,.33) -- (.5,.5) ;
			\draw[thick] (.33,.33) -- (0,.5) ;
			\fill (0,0) circle (.03);
			\fill (.5,0) circle (.03);
			\fill (1,0) circle (.03);
			\fill (.5,.5) circle (.03);
			\fill (0,1) circle (.03);
			\fill (0,.5) circle (.03);
			\fill (.33,.33) circle (.03);
		\end{tikzpicture}
		\caption{Cubical subdivision of a triangle.}
		\label{fig:cub_sub_tri}
	\end{subfigure}
	\caption{Subdivisions of a triangle.}
\end{figure}

\begin{dfn}[Cellular Cosheaves and Cellular Sheaves]
	Let $E$ be a regular CW-complex, $\mathrm{Cell}\,E$ denote the small category of its underlying poset of cells, and $R$ be a ring. A \emph{cellular cosheaf} (resp. \emph{sheaf}) of $R$-modules on $E$ is a contravariant (resp. covariant) functor:
	\begin{equation*}
		\cF:\mathrm{Cell}\,E \rightarrow \left(\begin{array}{c}\textrm{Category of}\\ R\textrm{-Modules}\end{array}\right),
	\end{equation*}
	A morphism $f$ between two cosheaves (or two sheaves) on $E$ is a natural transformation. The kernel, image, and cokernel of $f$ are formed “cell-wise”.
\end{dfn}

A cosheaf (resp. sheaf) on the cubical subdivision of a simplicial complex will be called a \emph{cubical cosheaf (resp. sheaf) on $K$}. Such an object is an assignment of an $R$-module to every pair of adjacent simplices $\sigma\leq \tau$ that is covariant in the first coordinate and contravariant in the second. 

\begin{dfn}\label{dfn:homo_cosheaf}
	Let $E$ be a regular CW-complex, $\cF$ be a cellular cosheaf (resp. $\mathcal{G}$ be a cellular sheaf), and $k$ be a non-negative integer. The space of cellular $k$-chains with values in $\cF$ (resp. of $k$-cochains with values in $\mathcal{G}$) is the vector space
	\begin{equation*}
		C_k(E;\cF)\coloneqq\bigoplus_{e^k\in E} \cF(e^k) \quad\textrm{resp.}\quad C^k(E;\mathcal{G})\coloneqq\prod_{e^k\in E} \mathcal{G}(e^k).
	\end{equation*} 
	A $k$-chain with values in $\cF$ is a formal sum $c=\sum_{e^k\in E} c_{e^k}\otimes e^k$ with $c_{e^k}\in\cF(e^k)$ while a $k$-cochain with values in $\mathcal{G}$ is seen as a map $\alpha:e^k\mapsto \alpha(e^k)\in \mathcal{G}(e^k)$. The collections of the vector spaces $(C_k(E;\cF))_{k\geq 0}$ and $(C^k(E;\mathcal{G}))_{k\geq 0}$ form a chain and a cochain complex respectively. Their respective boundary $\partial$ and differential $\df$ are defined, for all $e^{k\pm1}\in E$, by the following formulæ
	\begin{equation*}
		(\partial c)_{e^{k-1}}=\sum_{e^k\geq e^{k-1}}\cF(e^{k}\leftarrow e^{k-1}) (c_{e^k}) \quad\textrm{resp.}\quad \df\alpha(e^{k+1})=\sum_{e^{k}\leq e^{k+1}}\mathcal{G}(e^k\rightarrow e^{k+1}) \big(\alpha(e^k)\big).
	\end{equation*}  
	The homology and cohomology of these complexes are respectively denoted by $(H_k(E;\cF))_{k\geq 0}$ and $(H^k(E;\mathcal{G}))_{k\geq 0}$. If $R=\F_2$, their \emph{Betti numbers} are defined as $b_k(E;\cF)\coloneqq\dim H_k(E;\cF)$ and $b^k(E;\mathcal{G})\coloneqq \dim H^k(E;\mathcal{G})$ for all $k\geq 0$.
	We note that morphisms of (co)sheaves induce morphisms of (co)chain complexes and thus (co)homological morphisms.
	Whenever $c$ is a cycle, we denote its homology class by $[c]$. We use a similar notation for cocycles.
	
\end{dfn}

Let $\cF$ be a cellular cosheaf (resp. sheaf) on $E$ and $S$ be a subdivision of $E$. One can subdivide $\cF$ into a cosheaf (resp. sheaf) on $S$ via the following procedure: for every open cell $s$ of $S$ denote $e(s)$ the unique open cell of $E$ that contains $s$, then define $\cF(s)\coloneqq\cF(e(s))$ and $\cF(s_1\rightarrow s_2)\coloneqq \cF(e(s_1)\rightarrow e(s_2))$. It yields quasi-isomorphisms of complexes
\begin{equation*}
	\begin{array}{rcl}
		C_k(E;\cF) & \rightarrow & C_k(S;\cF) \\
		v\otimes e^k  & \mapsto & \sum_{s^k\subset e^k} v\otimes s^k
	\end{array} \quad\mathrm{resp.}\quad \begin{array}{rcl}
		C^k(S;\cF) & \rightarrow & C^k(E;\cF) \\
		\alpha  & \mapsto & (e^k\mapsto \sum_{s^k\subset e^k}\alpha(s^k)).
	\end{array}
\end{equation*}
If $K$ is a simplicial complex, $S$ is its cubical subdivision, and $\cF$ is a cubical cosheaf (resp. sheaf) on $K$, we denote by $(\Omega_k(K;\cF))_{k\geq 0}$ and $(\Omega^k(K;\mathcal{G}))_{k\geq 0}$ the chain complexes $(C_k(S;\cF))_{k\geq 0}$ and $(C^k(S;\mathcal{G}))_{k\geq 0}$ respectively. 

\begin{rem}\label{rem:cellular_sheaves_alexandrov}
		Let $K$ be a finite regular CW-complex with underlying topological space $|K|$. We denote by $A(K)$ the \emph{Alexandrov space} of $K$, i.e. the finite topological space whose points are the cells of $K$ and for which a basis of open sets is given by the sets $\{e'\geq e\}$ for all $e\in K$. The map $f:K\rightarrow A(K)$, that sends a point $x\in K$ to the unique open cell that contains it, is continuous and open. In this framework, a cellular sheaf $\cF$ on $K$ is equivalent to a sheaf on $A(K)$. The application $H^k(|K|;f^*\cF)\rightarrow H^k(A(K);\cF)$ is an isomorphism for all $k$. This follows the identity $f_*f^*\cF=\cF$ and the use of the Leray spectral sequence. One can also note that the cellular complex $(C^k(K;\cF))_{k\geq 0}$ is the \v{C}ech complex associated to the Leray cover $\{\{e\geq v\}:v \textnormal{ vertices}\}$. For a more thorough exposition of the theory of cellular sheaves and cosheaves one can look at \cite[\S Constructible Sheaves]{kashiwara_sheaves_2010} (under the name \emph{constructible sheaves}), \cite{shepard_cellular_1985}, or even \cite{chenal_poincare-lefschetz_2025} where we discuss \emph{dihomologic cosheaves} which is a generalisation of cubical cosheaves to the setting of regular CW-complexes.
\end{rem}

The particular combinatorics of the cubical subdivision of a simplicial complex allows for the definition of a cellular cup-product without making additional choices. The same is true for the barycentric subdivision but it introduces a far greater number of cells. 

\begin{dfn}\label{dfn:cup}
	Let $\cF$, $\mathcal{G}$, $\mathcal{H}$ be three cubical sheaves of $R$-modules on a simplicial complex $K$, and $m:\cF\otimes_R \mathcal{G}\rightarrow \mathcal{H}$ be a morphism of sheaves. The associated \emph{cup-product} is the $R$-bilinear morphism $\cupp_m:\Omega^p(K;\cF)\otimes_R \Omega^q(K;\mathcal{G}) \rightarrow \Omega^{p+q}(K;\mathcal{H})$ defined by
	\begin{equation*}
		\alpha\cupp_m\beta:[\sigma^a;\sigma^b]\mapsto \sum_{\sigma^b\geq \sigma^c\geq \sigma^a} m\big(f(\alpha(\sigma^a;\sigma^c))\otimes g(\beta(\sigma^c;\sigma^b))\big)\in \mathcal{H}(\sigma^a;\sigma^b)
	\end{equation*}
	where $b-a=p+q$, $c-a=p$, $f=\cF([\sigma^a;\sigma^c]\rightarrow [\sigma^a;\sigma^b])$, and $g=\mathcal{G}([\sigma^c;\sigma^b]\rightarrow [\sigma^a;\sigma^b])$. It satisfies the Leibniz rule $\df(\alpha\cupp_m\beta)=(\df\alpha)\cupp_m\beta+(-1)^p \alpha\cupp_m(\df\beta)$ and induces a product in cohomology. Moreover, if $s:\cF\otimes_R\mathcal{G}\rightarrow \mathcal{G}\otimes_R\cF$ is the swap, we have $[\beta]\cupp_{ms}[\alpha]=(-1)^{pq}[\alpha]\cupp_m[\beta]$. If $m$ is the law of a sheaf of associative algebras, the induced cohomology algebra is associative as well. In the sequel, the product $m$ will be clear from the context and we will omit it in the notation. 
\end{dfn}

\begin{dfn}\label{dfn:contra}
	Let $A,B,C$ be three finite-dimensional vector spaces over $\F_2$. To every product $\cupp:A\otimes B\rightarrow C$ one can associate two \emph{contractions}: the two morphisms $\capp:C^*\otimes A\rightarrow B^*$ and $\capp:B\otimes C^*\rightarrow A^*$ such that for all $\alpha\in A$, all $\beta\in B$, and all $c\in C^*$
	\begin{equation*}
		\la \alpha,\beta\capp c\ra = \la \alpha\cupp \beta,c\ra = \la \beta,c\capp \alpha\ra.
	\end{equation*}
\end{dfn}

Elementary computations show that if $\cupp$ is the law of an associative algebra $A$, then we have $(c\capp\alpha)\capp\beta=c\capp(\alpha\cupp\beta)$ and $\alpha\capp(\beta\capp c)=(\alpha\cupp\beta)\capp c$. If, in addition, $A$ is graded over an Abelian group\footnote{Here $L$ will be either $\Z$ or $\Z^2$.} $L$ and $\alpha$ is homogeneous of degree $p\in L$, $c\mapsto \alpha\capp c$ and $c\mapsto c\capp \alpha$ are graded operations of degree $-p$. Furthermore, if $A$ is commutative, $c\capp\alpha=\alpha\capp c$. 

Combining Definitions~\ref{dfn:cup}~and~\ref{dfn:contra}, we get \emph{cap-products} associated to $m$, both on the (co)chains level and in (co)homology. The (co)chain cap-products also satisfy Leibniz formulæ easily derived from the cup-product Leibniz formula and the defining property of the contractions. Naturally, the induced products in (co)homology are the contractions of the cohomological cup-product. 

In the following, we will only almost exclusively consider (co)sheaves of $\F_2$-vector space, hence if we do not mention any ring $R$, we assume $R=\F_2$.

\subsection{Primitive Patchworks and T-Hypersurfaces}\label{sec:prim_patch}

\paragraph{Tropical Hypersurfaces}
Let $P$ be a lattice polytope in $\M{\R}$ yielding a smooth toric variety $\mathcal{Y}$, we say that $P$ is \emph{smooth}. This means that (1) $P$ is simple, i.e. every codimension $q$ face of $P$ is the intersection of exactly $q$ faces of codimension $1$; and (2) for every vertex $V$ of $P$, $M$ is spanned by the lattice points contained in the tangent lines to the edges of $P$ that are adjacent to $V$; see \cite[\S 2.1]{fulton_introduction_1993} or \cite[Theorem 2.4.3]{cox_david_a_toric_2011}.

The variety $\mathcal{Y}$ is then defined by the normal fan $C\subset \N{\R}$ of $P$. If $\mathbb{K}\in\{\R;\C\}$, let $(\mathbb{K},\times)$ denote its multiplicative monoid. For every cone $c\in C$, $c^+$ denotes $\{u\in \N{R}\;|\;\la v,u\ra\geq 0,\,\forall v\in c\ra$. Then, the $\mathbb{K}$-locus $\mathcal{Y}(\mathbb{K})$ is homeomorphic to the inductive limit
\begin{equation*}
	\mathcal{Y}(\mathbb{K})\approx \underset{c\in C}{\mathrm{ind.lim}}\, \Hom_\mathrm{mon.}(c^+\cap M;(\mathbb{K},\times)).
\end{equation*}
The \emph{tropical locus} $\mathcal{Y}(\T)$ is defined with this formula $\mathcal{Y}(\T)\coloneqq\mathrm{ind.lim}_{c\in C} \Hom_\mathrm{mon.}(c^+\cap M;(\T,\times))$. Note that, with this definition, $\mathcal{Y}(S)$ makes sense for any semi-group such as $(\R_+,\times)$. Since $\log:(\R_+;\times)\rightarrow (\T;\times)$ is an isomorphism of topological monoids, we have a homeomorphism $\log:\mathcal{Y}(\R_+)\rightarrow \mathcal{Y}(\T)$. The tropical loci of the toric subvarieties of $Y$ induce a filtration of $\mathcal{Y}(\T)$. This filtration is the skeletal filtration of a cellular structure for which every closed cell is the tropical locus of a subvariety. The moment map of $P$,
\begin{equation*}
	\mu:y\in \mathcal{Y}(\C)\mapsto \frac{\sum_{v\in P\cap M} |y^v|v}{\sum_{v\in P\cap M}|y^v|}\in P
\end{equation*}
is a homeomorphism when restricted to $\mathcal{Y}(\R_+)$ and the composition $\mu\circ\log^{-1}:\mathcal{Y}(\T)\rightarrow P$ is even an isomorphism of cellular complexes. 

Let $(a_v)_{v\in P\cap M}\in\T^{P\cap M}$ not all tropically zero, i.e. not all equal to $-\infty$. This corresponds to a tropical section of the line bundle $\O(P)$ of $\mathcal{Y}$. The tropical hypersurface $\mathcal{X}(\T)\subset \mathcal{Y}(\T)$ defined by this section is the closure in $\mathcal{Y}(\T)$ of the non-differentiability locus of the function $f:\N{\R}\rightarrow \T$ defined by 
\begin{equation*}
	f(u)=\max_{v\in P\cap M}\big(\la v,u\ra + a_v\big).
\end{equation*} The closure of this non-differentiability locus is also called the \emph{tropical vanishing locus} of $f$. Since it is the closure of a polyhedral complex, it is endowed with a cellular structure. Moreover, when the convex hull of the lattice points $v\in P\cap M$ for which $a_v\neq -\infty$ is $P$ itself, then the closures of the connected components of $\mathcal{Y}(\T)\setminus \mathcal{X}(\T)$ endow $\mathcal{Y}(\T)$ with a cellular structure of which $\mathcal{X}(\T)$ is a subcomplex. Using the moment map, one can also see this as a cellular subdivision of $P$; see Figure~\ref{fig:trop_cub2}. In addition, $f$ also defines a polyhedral subdivision $K$ of $P$, called it \emph{dual subdivision}\footnote{In the framework of Poincaré duality, the subdivision of $\mathcal{Y}(\T)$ induced by $\mathcal{X}(\T)$ is the actual dual subdivision of $K$; cf. \cite[\S 64]{munkres_elements_1984}.}. A cell of $K$ is the projection along $\M{\R}\times\R\rightarrow \M{\R}$ of a compact face of the polyhedron 
\begin{equation*}
	\mathrm{Conv.Hull} \bigcup_{\substack{v\in P\cap M\\a_v\neq -\infty}} \{v\}\times [-a_v;+\infty).
\end{equation*}
See \cite[Definition~2.3.5]{mikhalkin_tropical_nodate} for instance. Figure~\ref{fig:dual_sub_cub} depicts an example.  The hypersurface $\mathcal{X}(\T)$ is said to be \emph{smooth} when $K$ is a \emph{primitive triangulation}.  

\begin{figure}[ht]
	\centering
	\begin{subfigure}[t]{.45\textwidth}
		\centering
		\input{tropical_subdivision.tex}
		\caption{The subdivision of $\mathcal{Y}(\T)$ induced by $\mathcal{X}(\T)$.}
		\label{fig:trop_cub2}
	\end{subfigure}
	\hfill
	\begin{subfigure}[t]{.45\textwidth}
		\centering
		\input{dual_subdivision.tex}
		\caption{The dual subdivision $K$.}
		\label{fig:dual_sub_cub}
	\end{subfigure}
	\caption{The tropical curve of tropical projective plane induced by the tropical polynomial $\max(0,x+2,y+2,2x+2,x+y+3,2y+2,3x,2x+y+2,x+2y+2,3y)$.}
\end{figure}

\begin{dfn}\label{dfn:prim_conv}
	 A \emph{primitive triangulation} $K$ of $P$ is a polyhedral subdivision of $P$ into integral simplices of minimal integral volume\footnote{The Lebesgue mesure normalised so that a parallelogram on a basis of the lattice $M$ has unit volume.}, i.e. $\frac{1}{(n+1)!}$.
	 
	 A triangulation $K$ of $P$ is \emph{convex} if there is a convex piecewise affine function $g:P\rightarrow \R$ whose domains of affinity are exactly the simplices of $K$.
\end{dfn}

A polyhedral subdivision of $P$ arrises as the dual subdivision of a tropical hypersurface if and only if it is convex and its vertices are lattice points. Using the duality property, one can still build replacements $Y$ and $X$ of $\mathcal{Y}(\T)$ and $\mathcal{X}(\T)$ respectively, when $K$ is not convex. The complex $Y$ is a cellular subdivision of $P$ and $X$ is a subcomplex of $Y$. In the convex case, $Y$ and $X$ are cellularly isomorphic to $\mathcal{Y}(\T)$ and $\mathcal{X}(\T)$. Moreover, these complexes $X\leq Y$, together with the triangulation $K$, retain enough information to be able to compute the tropical homology of $\mathcal{Y}(\T)$ and $\mathcal{X}(\T)$. Therefore, one can extend tropical homology beyond the convex setup as in \cite{brugalle_combinatorial_2024}.  

\paragraph{Replacement Cellular Complexes}
 
Let $K$ be a primitive triangulation of a polytope $P$.

\begin{dfn}\label{dfn:XY}
	Let $\sigma^q$ be a simplex of $K$ and $Q$ be a $p$-face of $P$ with $\sigma^q\subset Q$. The \emph{dual cell} of $\sigma^q$ in $Q$, denoted by $[\sigma^q,Q]$, is the support of the cubical complex 
	\begin{equation*}
		\bigcup_{\sigma^q\leq\sigma^p\subset Q}[\sigma^q,\sigma^p].
	\end{equation*}
	This is a closed ball of dimension $p-q$. We denote $[\sigma^q,P]$ by $(\sigma^q)^*$ and refer to it as the \emph{dual cell} of $\sigma^q$. The collection $\{[\sigma^q,Q]:\sigma\subset Q\}$ is the set of closed cells of a regular subdivision $Y$ of $P$. This cellular structure is called the \emph{coarse structure} of $Y$ while the cubical subdivision of $K$, which naturally subdivides $Y$, is its \emph{fine structure} or \emph{fine subdivision}. A coarse cell of $Y$ is \emph{bounded} if it is dual to an interior simplex of $K$. See Figure~\ref{fig:coarse_fine_YX}.
\end{dfn}

\begin{dfn}
	The cosheaf $\Sed$ on $P$ associates to a face $Q$ of $P$ the group
		\begin{equation*}
			\Sed(Q)\coloneqq\big\{v\in N \;|\; \la\alpha,v\ra=0,\, \forall \alpha\in TQ\cap M\big\}\otimes\F_2,
		\end{equation*}
	If $Q_1\leq Q_2$ are faces of $P$, the transition morphism $\Sed(Q_2\leftarrow Q_1)$ is the inclusion (n.b. by the subdivision construction, $\Sed$ also yields cosheaves on $K$ and both structures of $Y$). 
\end{dfn}

\begin{dfn}\label{dfn:X}
	The hypersurface $X\leq Y$ is the subcomplex spanned by the dual cells of the edges of $K$. As for $Y$, $X$ has a coarse and a fine subdivision. See Figure~\ref{fig:coarse_fine_YX}. 
\end{dfn}

\begin{figure}[H]
	\centering
	\begin{subfigure}[t]{.45\textwidth}
		\centering
		\input{coarse_sub.tex}
		\caption{The coarse subdivision $Y$ of $P$. The subcomplex $X$ is in bold.}
		\label{fig:coarse_sub}
	\end{subfigure}
	\hfill
	\begin{subfigure}[t]{.45\textwidth}
		\centering
		\input{fine_subd.tex}
		\caption{The fine subdivision of $Y$, i.e. the cubical subdivision of $K$.}
		\label{fig:fine_subd}
	\end{subfigure}
	\caption{The complexes $X$ and $Y$ associated to the primitive triangulation $K$ of Figure~\ref{fig:dual_sub_cub}.}
	\label{fig:coarse_fine_YX}
\end{figure}

\begin{dfn}
	Let $Z\in\{X;Y\}$. If $\cF$ is a cellular cosheaf (resp. sheaf) on $Z$, we denote by $\Omega_*(Z;\cF)$ (resp. $\Omega^*(Z;\cF)$) the (co)chain complex on the fine subdivision of $Z$.
\end{dfn}

\paragraph{Real Loci of Toric Varieties} Let $\mu^{-1}:P\rightarrow \mathcal{Y}(\R_+)$ denote the inverse of the moment map. The real locus $\mathcal{Y}(\R)$ can be recovered as the gluing of $2^{n+1}$ copies of $P$. Indeed, we have a quotient map
\begin{equation}\label{eq:phi}
	\begin{array}{rcl}
		\phi:\N{2}\times P & \longrightarrow & \mathcal{Y}(\R)\\
		(u;x) & \longmapsto & (-1)^u\cdot \mu^{-1}(x)
	\end{array}
\end{equation}
for which $(u_0,x_0)$ is identified with $(u_1,x_1)$ if and only if $x_0=x_1$ and $u_0-u_1$ belongs to the subspace $(TQ^\perp\cap N)\otimes \F_2$ where $Q$ is the smallest face of $P$ containing $x_0$; cf. \cite[Theorem 5.4]{gelfand_discriminants_1994}, \cite[\S 4.1]{fulton_introduction_1993}, or \cite{davis_convex_1991}. This defines a cellular structure $\R P$ on $\mathcal{Y}(\R)$. Moreover, every subdivision of $P$ induces a subdivision of $\R P$. This way, $\R P$ is endowed with a triangulation $\R K$, the coarse structure $\R Y$ and its fine subdivision. In the future, we will always see $P$ as embedded in $\mathcal{Y}(\R)$ via $\mu^{-1}$. See Figure~\ref{fig:pol_coord}.     

\begin{figure}[ht]
	\centering
	\input{polar_coordinate_real_locus.tex}
	\caption{The morphism $\phi$ gluing four copies of the standard simplex $P$ endowed with the subdivision $Y$ from Figure~\ref{fig:coarse_sub} into the real projective plane with its cellular structure $\R P$ featuring the associated subdivision $\R Y$ (n.b. arrowed edges have to be glued together according to the displayed orientation).}
	\label{fig:pol_coord}
\end{figure}

In addition, $\N{2}$ acts on $\mathcal{Y}(\R)$ as the $2$-torsion of the real locus $\N{\R^\times}$ of the algebraic torus and, in that regard, $\phi$ is equivariant. Therefore, if $E\subset \mathcal{Y}(\R)$ and $u\in\N{2}$, $u\cdot E\subset \mathcal{Y}(\R)$ denotes the image of $E$ by the action of $u$. Hence, if $E\subset P$, $u\cdot E=\phi(\{u\}\times E)$. We can also note that the cellular structures $\R P$, $\R K$, and $\R Y$ are all equivariant. Moreover, for all coarse cells $e^q\in Y$, $\Sed(e^q)\subset \N{2}$ is the isotropy group of the cell $e^q$. 

\paragraph{T-Hypersurfaces} We consider $K$ a primitive triangulation of $P$. Since $K$ is primitive, its vertices are exactly the lattice points contained in $P$.

\begin{dfn}\label{dfn:t-hyp}
	A \emph{sign distribution} $\varepsilon$ is a cochain in $C^0(K;\F_2)$. Traditionally, a sign distribution takes its values in $\{\pm1\}$ but here it will be more convenient to work with $\F_2$ instead. Naturally, we identify $\F_2$ with $\{\pm1\}$ using the unique group isomorphism $\eta\in \F_2\mapsto (-1)^\eta\in\{\pm 1\}$, so that $0$ is $+1$ and $1$ is $-1$.
	
	Let $\varepsilon$ be a sign distribution, $\sigma^1$ be an edge of $K$, $v_0,v_1$ be its vertices, and $e^n\in Y$ be its dual cell. We set
	\begin{equation*}
		\Arg(e^n)\coloneqq \{u\in \N{2}\,|\, \la v_0-v_1,u\ra + \df\varepsilon(\sigma^1)=1\}\subset \N{2}.
	\end{equation*} 
	The \emph{T-hypersurface} associated to $K$ and $\varepsilon$ is the subcomplex $\X\leq \R Y$ defined by
	\begin{equation*}
		\X\coloneqq \bigcup_{e^n\in X} \bigcup_{u\in \Arg(e^n)} u\cdot e^n.
	\end{equation*}
	In other words, $\X\cap P$ is made of the duals of the edges $\sigma^1$ of $K$ carrying opposite “signs” i.e. such that $\df\varepsilon(\sigma^1)=1$. To determine $\X\cap u\cdot P$ for $u\in\N{2}$, one has to apply the same rule after replacing $\varepsilon$ by $v\mapsto \varepsilon(v)+\la v,u\ra$. See Figure~\ref{fig:t_curve} for an example. 
\end{dfn}

\begin{figure}[H]
	\centering
	\input{t_curve.tex}
	\caption{The T-curve $\X$ associated to the sign distribution $\varepsilon$ depicted in the upper right corner. In each of the four copies of $P$, we have represented $\X\cap u\cdot P$ for all $u\in \N{2}$. In addition, we have attached a small copy of the triangulation $K$ with the appropriate modification of $\varepsilon$ as explained in Definition~\ref{dfn:t-hyp}.}
	\label{fig:t_curve}
\end{figure}

\begin{prop} \textnormal{\cite[Proposition 4.11]{brugalle_combinatorial_2024}}
	The complex $\X$ is a PL-manifold.
\end{prop}

\begin{rem}
	Let $\omega_{\X}\in C^1(\R K;\F_2)$ denote the unique cocycle satisfying 
	\begin{equation*}
		\omega_{\X}(u\cdot \sigma^1)= \la v_0-v_1,u \ra + \df\varepsilon (\sigma^1),
	\end{equation*}
	for all $\sigma^1\in K$ with vertices $v_0,v_1$ and all $u\in \N{2}$. The hypersurface $\X$ is the the \emph{dual hypersurface} of $\omega_\X$, i.e. a cell $e^n$ of $\R Y$ belongs to $\X$ if and only if the dual edge $\sigma^1$ of $\R K$ satisfies $\omega_\X(\sigma^1)=1$. As a consequence, the cohomology class $[\omega_\X]\in H^1(\R Y;\F_2)$ is Poincaré dual to the image of the fundamental class of $\X$ in $H_n(\R Y;\F_2)$. 
\end{rem}

\begin{dfn}
	By construction, the restriction of the moment map induces a surjective cellular map $\mu:\X\rightarrow X$. The \emph{real phase structure} $\Arg$ is the cosheaf of sets on $X$ such that 
	\begin{equation*}
		\Arg(e^q)\coloneqq \{u\in \N{2}\;|\;u\cdot e^q\in \X\}/\Sed(e^q), \quad\textnormal{for all coarse cells }e^q\in X.
	\end{equation*}
	This set is in natural bijection with the $\mu$-lifts of $e^q$ that belong to $\X$. 
\end{dfn}

\begin{thm} \textnormal{\cite[Theorem 4.3.A]{viro_patchworking_2006}, \cite[Théorème 4.2]{risler_construction_1993}, or \cite[Theorem 5.6]{gelfand_discriminants_1994}.}
	If $K$ is convex, there exists an algebraic hypersurface $\mathcal{X}$ of $\mathcal{Y}$ with Newton polytope $P$ whose real locus $\mathcal{X}(\R)$ is ambiantly isotopic to the support of $\X$.
	
	More precisely, if $K$ is the dual subdivision of the tropical hypersurface defined by the tropical Laurent polynomial $u\in \N{\R}\mapsto \max_{v\in P\cap M}(\la v,u\ra+a_v)\in \R$, then $\mathcal{X}$ can be taken as the vanishing locus of the Laurent polynomial
	\begin{equation*}
		f_t=\sum_{v\in P\cap M} (-1)^{\varepsilon(v)}t^{a_v}x^v,
	\end{equation*}
	for $t>0$ sufficiently large.  
\end{thm}

\begin{figure}[H]
	\centering
	\includegraphics[width=.35\textwidth]{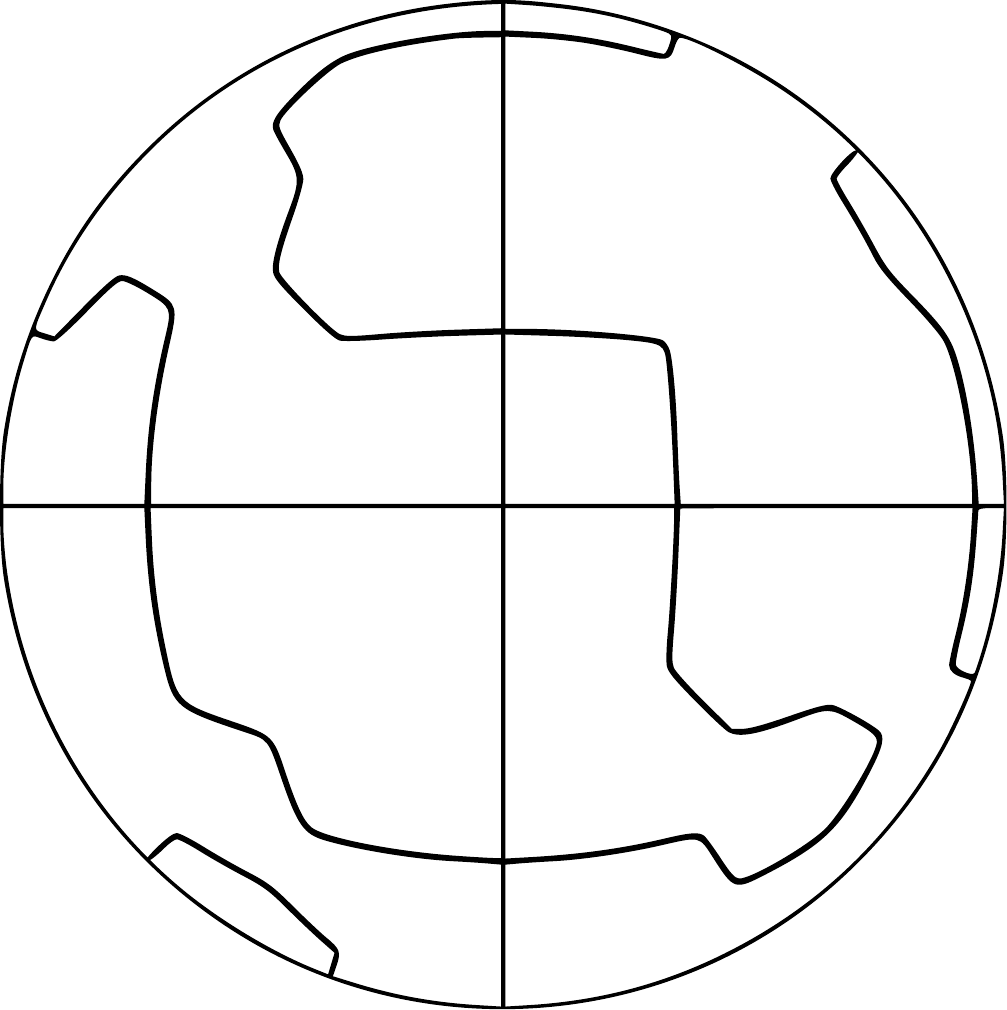}
	\caption{The plane cubic of equation $f_t=0$ of Example~\ref{ex:viro_patch} with $t=e^{29/2}$ in the parametrisation $(x;y)\in\{x^2+y^2\leq 1\}\mapsto [x^{29}:y^{29}:(1-x^2-y^2)^{29}]\in \mathbb{P}^2(\R)$ of the real projective plane.}
	\label{fig:alg_patch}
\end{figure}

\begin{ex}\label{ex:viro_patch}
	To illustrate Viro's Theorem, one can consider the triangulation $K$ of Figure~\ref{fig:dual_sub_cub} and the T-curve of Figure~\ref{fig:t_curve}. Then, an equation of $\mathcal{X}$ is provided by 
	\begin{equation*}
		f_t=1+x^3+y^3+t^2(-x-y-x^2-y^2+x^2y-xy^2)-t^3xy
	\end{equation*}
	for $t>0$ sufficiently large. Figure~\ref{fig:alg_patch} depicts the real locus $\mathcal{X}(\R)$ for a certain value of $t>0$. One can compare the result with Figure~\ref{fig:t_curve}. 
\end{ex}

\subsection{Tropical Homology}\label{sec:trop_homo}

\noindent \textbf{Tropical (Co)Sheaves.} These objects were introduced in \cite{itenberg_tropical_2019} to define an analog of Hodge theory for tropical varieties. They form a family of cellular (co)sheaves of Abelian groups associated to every tropical variety and allow us to define tropical (co)homology. We use their generalisations to the complexes $Y$ and $X$ as in \cite{brugalle_combinatorial_2024} and only provide the definition of their $\F_2$ version.

\begin{dfn}\label{dfn:trop_cosheaves}
	Let $P$ be a smooth lattice polytope of $\M{\R}$, $K$ be a primitive triangulation of $P$ and $X,Y$ the complexes of Definitions~\ref{dfn:XY} and \ref{dfn:X}. Let $p\geq 0$ be an integer,
	\begin{enumerate}	
		\item The cosheaf $\cF_p^Y$ on the coarse structure of $Y$ is the quotient $\bigwedge^p \N{2}/\Sed$, where $\N{2}$ is understood as a constant cosheaf.
		\item The cosheaf $\cF_p^X$ on the coarse structure of $X$ is the sub-cosheaf of $\cF^Y_p|_X$ such that if $[\sigma^{q},Q]$ is a coarse cell of $X$ (i.e. $q\geq 1$), we have 
		\begin{equation*}
			\cF_p^X(\sigma^q,Q)\coloneqq \sum_{\sigma^1\leq\sigma^q}{\textstyle \bigwedge^p}\Big(((T\sigma^1)^\perp\cap\N{})\otimes\F_2/\Sed(\sigma^1,Q)\Big)\subset \cF^Y_p(\sigma^q,Q).
		\end{equation*}
		An element of a group of one of these two cosheaves is referred to as a \emph{framing}. 
		\item The sheaf $\cF^p_Y$ (resp. $\cF_X^p$) denotes $\Hom(\cF^p_Y;\F_2)$ (resp. $\Hom(\cF^p_X;\F_2)$).
	\end{enumerate}
	Let $p,q\geq 0$ be integers and $Z\in\{X;Y\}$. We denote by $H_{p,q}(Z;\F_2)$ and (resp. $H^{p,q}(Z;\F_2)$) the $q$\textsuperscript{th} (co)homology group of $\cF^Z_p$ (resp. $\cF_Z^p$). We denote by $i_{p,q}:H_{p,q}(X;\F)\rightarrow H_{p,q}(Y;\F_2)$ the morphism induced by the inclusion $X\leq Y$. Moreover, $h^{p,q}(Z)$ denotes $\dim H^{p,q}(Z;\F_2)$ (n.b. this is the same as $\dim H_{p,q}(Z;\F_2)$ by the Universal Coefficient Theorem).
\end{dfn}

\paragraph{Products} By definition, $(\cF^p_Y)_{p\geq0}$ is a sheaf of graded algebras, namely the exterior algebra of $\cF_Y^1$. This sheaf of algebras is associative and, since we work over $\F_2$, also commutative.

\begin{dfn}
	Let $\sigma^p$ be a simplex of $K$. We denote by $\omega(\sigma^p)\in \bigwedge^p\M{2}$, the volume element of the vector subspace $(T\sigma\cap M)\otimes \F_2\subset \M{2}$.
\end{dfn}

\begin{prop}
	\textnormal{\cite[Propositions 4.10 and 4.12]{chenal_poincare-lefschetz_2025}}. For all coarse cells $[\sigma^q,Q]$ of $X$, the group $\cF^X_p(\sigma^q,Q)$ is the kernel of the map
	\begin{equation*}
		 u\in \cF^Y_p(\sigma^q,Q) \mapsto \omega(\sigma^q)\cdot u \in \cF^Y_{p-q}(\sigma^q,Q),
	\end{equation*}
	where “$~\cdot~$” denotes the contraction of the algebra structure of $(\cF_Y^p(\sigma^q,Q))_{p\geq 0}$; see Definition~\ref{dfn:contra}. 
\end{prop}

This implies that, for every coarse cell $[\sigma^q,Q]$ of $X$, the kernel of the surjective morphism $\bigoplus_{p\geq 0}\cF^p_Y(\sigma^q,Q)\rightarrow  \bigoplus_{p\geq 0}\cF^p_X(\sigma^q,Q)$ is the ideal generated by $\omega(\sigma^q)$, i.e. $(\cF^p_X)_{p\geq 0}$ is also a sheaf of graded algebras. Following Definition~\ref{dfn:cup}, we have \emph{cup-products} on the tropical cohomology of $X$ and $Y$. If $\alpha\in \Omega^{q_1}(X;\cF^{p_1}_X)$ and $\beta\in \Omega^{q_2}(X;\cF^{p_2}_X)$, we compute the cup-product as follows
\begin{equation}\label{eq:cupp}
		(\alpha\cupp\beta)(\sigma^a,\sigma^{b})\coloneqq \sum_{\sigma^a\leq\sigma^{c}\leq \sigma^{b}} \alpha(\sigma^a,\sigma^{c})\wedge\beta(\sigma^{c},\sigma^{b}),
\end{equation}
with $b-a=q_1+q_2$ and $c-a=q_1$ and being understood that both $\alpha(\sigma^a,\sigma^c)$ and $\beta(\sigma^c,\sigma^{b})$ have to be respectively transferred in $\cF_X^{p_1}(\sigma^a,\sigma^b)$ and $\cF_X^{p_2}(\sigma^a,\sigma^b)$ before being wedged. Since $(\cF^p_X)_{p\geq 0}$ and $(\cF^p_Y)_{p\geq 0}$ are sheaves of associative and commutative graded algebras over $\F_2$, the tropical cohomologies of $X$ and $Y$ with coefficient in $\F_2$ are associative and commutative bigraded $\F_2$-algebra. In addition, the restriction 
\begin{equation*}
	(\Omega^q(Y,\cF^p_Y))_{p,q\geq 0}\rightarrow (\Omega^q(X,\cF^p_X))_{p,q\geq 0},
\end{equation*}
is a morphism of bigraded algebras. Following Definition~\ref{dfn:contra}, they also yield \emph{cap-products} defined by the relations $\la \alpha, \beta\capp c\ra=\la \alpha\cupp\beta,c\ra$ and $\la\beta, c\capp\alpha\ra=\la \alpha\cupp\beta,c\ra$. With these conventions, the associativity rules $\alpha\capp(\beta\capp c)=(\alpha\cupp\beta)\capp c$ and $(c\capp \alpha)\capp\beta=c\capp(\alpha\cupp\beta)$ both hold. A direct computation provides the formulæ for the cap-products of a chain $c\in \Omega_{q_1+q_2}(Z;\cF^Z_{p_1+p_2})$ and cochains $\alpha\in \Omega^{q_1}(Z;\cF_Z^{p_1})$ and $\beta\in \Omega^{q_2}(Z;\cF^{p_2}_Z)$
\begin{equation}\label{eq:cap_prod}
	\left\{\begin{array}{lll}
		 (\beta\capp c)_{\sigma^{k},\sigma^{l}}=\sum_{\sigma^{m}\geq\sigma^{l}}\beta(\sigma^{l},\sigma^{m})\cdot c_{\sigma^k,\sigma^{m}}, & \scriptstyle\mathrm{with}\;l-k=q_1 \;\mathrm{and}\; l-m=q_1+q_2 \\[5pt]
		(c\capp \alpha)_{\sigma^k,\sigma^{l}}=\sum_{\sigma^{j}\leq\sigma^k} c_{\sigma^{j},\sigma^{l}} \cdot \alpha(\sigma^{j},\sigma^k), & \scriptstyle\mathrm{with}\;l-k=q_2 \;\mathrm{and}\; l-j=q_1+q_2
	\end{array}\right.
\end{equation} 
where the “$~\cdot~$” denote the contractions of the wedge used in (\ref{eq:cupp}). Since the cohomological cup-product is commutative, so is the homological cap-product, i.e. $[\alpha]\capp [c]=[c]\capp[\alpha]$. 

Tropical homology satisfies Poincaré duality; cf. \cite[Theorem 5.3]{jell_lefschetz_2018} and \cite[Theorem 3.3]{brugalle_combinatorial_2024}. This means that the cap-product with the fundamental class $[X]$ yields isomorphisms
$H^{p,q}(X,\F_2) \rightarrow H_{n-p,n-q}(X,\F_2)$, for all integers $0\leq p,q\leq n$. If $e^n$ is a coarse cell of $X$, $\cF_n^X(e^n)$ has dimension $1$ and is spanned by the volume element $\dv_{e^n}$ of $\cF_1^X(e^n)$. The fundamental class $[X]$ is represented by the cycle
\begin{equation*}
	\sum_{e^n\in X} \dv_{e^n}\otimes e^n.
\end{equation*}
Poincaré duality allows for the definition of a non-degenerate \emph{intersection pairing} 
	\begin{equation*}
		H_{p,q}(X;\F_2)\otimes H_{n-p,n-q}(X;\F_2) \rightarrow \F_2
	\end{equation*}
denoted by $[a]\cdot[b]$ for all classes $[a]$ and $[b]$. It is defined by $[a]\cdot[b]=\langle [\alpha],[b]\rangle$ where the class $[\alpha]$ satisfies $[\alpha]\cap [X]=[a]$.

\subsection{The Renaudineau--Shaw Spectral Sequence}\label{sec:RSSS}

The original construction of the homological Renaudineau--Shaw spectral sequence can be found in \cite[\S4]{renaudineau_bounding_2023}. Here, we recall its quickest construction. For a presentation of the cohomological construction, we refer to \cite[\S5]{chenal_poincare_2026}. 

\begin{dfn}\label{dfn:sign_cosheaves}
	Let $\K^{\R Y}$ denote the cosheaf of group algebras on $Y$ defined by $\F_2[\cF_1^Y]$ and $e^q$ be a coarse cell of $Y$. Since $\cF_1^Y=\N{2}/\Sed$, the space $\K^{\R Y}(e^q)$ is the vector space freely generated by the lifts of $e^q$ along $\mu:\R Y\rightarrow Y$. If $u\in \cF^Y_1(e^q)$, we denote the associated generator of $\K^{\R Y}(e^q)$ by $\x^u$.
	
	Let $\varepsilon$ be a sign distribution on $K$ and $e^q$ be a coarse cell of $X$. We denote by $\K^{\X}(e^q)$ the subspace of $\K^{\R Y}(e^q)$ spanned by lifts of $e^q$ contained in $\X$. This way $e^q\mapsto \K^\X(e^q)$ forms a cellular cosheaf $\K^\X$ on $X$. Equivalently, $\K^\X$ can be described as the cosheaf of vector spaces spanned by the cosheaf of sets $\mathcal{E}$.  
\end{dfn}

\begin{prop} \textnormal{\cite[Proposition 3.17]{renaudineau_bounding_2023}.} The two chain complexes $(C_q(\R Y;\F_2))_{q\geq 0}$ and $(C_q(Y;\K^{\R Y}))_{q\geq 0}$ are naturally isomorphic. Moreover, this natural isomorphism restricts to an isomorphism between $(C_q(\X;\F_2))_{q\geq 0}$ and $(C_q(X;\K^{\X}))_{q\geq 0}$. The same statement applies to the fine structures. 
\end{prop}

For all finite dimensional $\F_2$-vector spaces $V$, the group algebra $\F_2[V]$ is a local ring whose maximal ideal is the kernel of the augmentation map $\F_2[V]\rightarrow \F_2$ sending $\x^v$ to $1$ for every $v\in V$. Since it is a natural transformation, it allows for the definition of a filtration of the cosheaves of Definition~\ref{dfn:sign_cosheaves}.

\begin{dfn}
	Let $p\geq 0$, the cosheaf $\K^{\R Y}_{p}\subset \K^{\R Y}$ is the $p$\textsuperscript{th} power of the maximal ideal. This filtration of cosheaves induces a spectral sequence $(E^r_{p,q}(\R Y))_{p,q,r\geq 0}$ computing the homology of $\R Y$ with coefficients in $\F_2$. 
	
	Likewise, for $\varepsilon$ a sign distribution on $K$, the intersections $\K^{\X}_{p}\coloneqq \K^\X\cap \K^{\R Y}_p|_X$, for all integers $p\geq 0$, define a filtration of $\K^\X$. The associated spectral sequence is denoted by $(E^r_{p,q}(\X))_{p,q,r\geq 0}$ and is called the \emph{Renaudineau--Shaw spectral sequence}. It computes the homology of $\X$ with coefficients in $\F_2$.
	
	The inclusion $i:\K^\X\rightarrow \K^{\R Y}|_X$ induces a morphism of spectral sequences denoted by $(i^r_{p,q}:E^r_{p,q}(\X)\rightarrow E^r_{p,q}(\R Y))_{p,q,r\geq 0}$.     
\end{dfn}

\begin{rem}
	The Renaudineau--Shaw spectral sequence is graded so that the $r$\textsuperscript{th} boundary operator has bidegree $(r,-1)$, i.e. $\partial^r_{p,q}:E^r_{p,q}(\X)\rightarrow E^r_{p+r,q-1}(\X)$. This also has the following consequence, if $(F_pH_q(\X;\F_2))_{p,q\geq 0}$ denotes the induced filtration of the homology of $\X$, we have the exact sequences
	\begin{equation*}
		0\rightarrow F_{p+1}H_q(\X;\F_2) \rightarrow F_{p}H_q(\X;\F_2) \rightarrow E^\infty_{p,q}(\X)\rightarrow 0,\quad\textnormal{for all }p,q\geq 0,
	\end{equation*} 
	which are natural in this case, but differ from the usual conventions used for the total homology of a bicomplex.
	 
	\emph{A priori}, there should be coarse and fine spectral sequences respectively associated with the coarse and fine subdivisions of $X$ and $Y$. However, since the $(r+1)$\textsuperscript{th} page of a	spectral sequence is the homology of the $r\textsuperscript{th}$ page, the coarse and fine spectral sequences are naturally isomorphic from the 1\textsuperscript{st} page, i.e. the morphism of spectral sequences induced by subdivision provides isomorphisms $E^r_{p,q}(\X)_\mathrm{coarse}\rightarrow E^r_{p,q}(\X)_\mathrm{fine} $ for all $p,q\geq 0$ and all $r\geq 1$.
\end{rem}

\begin{prop}\label{prop:borel_viro} \textnormal{\cite[Lemma 4.8 and Proposition 4.10]{renaudineau_bounding_2023} and \cite[Proposition 4.22]{chenal_poincare_2026}.}
	For all integers $p\geq 0$, there is a unique surjective morphism of cosheaves $\bv_p:\K^{\R Y}_{p}\rightarrow \cF^Y_p$ such that for all coarse cells $e^q\in Y$ and all $u_0,...,u_p\in\cF^Y_1(e^q)$,
	\begin{equation*}
		\bv_p\left( \x^{u_0}\prod_{i=1}^p(1+\x^{u_i}) \right)= \bigwedge_{i=1}^p u_i.
	\end{equation*}
	Furthermore, its kernel is $\K^{\R Y}_{p+1}$ and it sends $\K^\X_{p}$ onto $\cF^X_p$. They are called the \emph{Borel--Viro morphisms}. In addition, they induce the following commutative square
	\begin{equation}\label{eq:borel-viro_isomorphisms}
		\begin{tikzcd}
			E^1_{p,q}(\R Y) \ar[r,"\bv_{p,q}" above, "\cong" below] & H_{p,q}(Y;\F_2) \\
			E^1_{p,q}(\X) \ar[r,"\bv_{p,q}" below, "\cong" above] \ar[u,"i^1_{p,q}" left] & H_{p,q}(X;\F_2) \ar[u,"i_{p,q}" right]
		\end{tikzcd}
	\end{equation}
	for all integers $p,q\geq 0$. We call the applications $\bv_{p,q}$ the \emph{Borel--Viro isomorphisms}.
\end{prop}

\begin{dfn}\label{dfn:partial_varepsilon}
	For all $p,q\geq 0$ and all $\varepsilon\in C^0(K;\F_2)$, $\partial_\varepsilon:H_{p,q}(X;\F_2)\rightarrow H_{p+1,q-1}(X;\F_2)$ denotes the unique morphism for which the following diagram commutes
	\begin{equation*}
		\begin{tikzcd}
			H_{p,q}(X;\F_2) \ar[r,"\partial_\varepsilon"] & H_{p+1,q-1}(X;\F_2) \\
			E^1_{p,q}(\X) \ar[u,"\bv_{p,q}","\cong" swap] \ar[r,"\partial^1_{p,q}" swap] & E^1_{p+1,q-1}(\X) \ar[u,"\cong", "\bv_{p,q}" swap]
		\end{tikzcd}
	\end{equation*}
\end{dfn}

In \cite{chenal_poincare_2026}, we show that the dual spectral sequence $(E^{p,q}_r(\X))_{p,q,r\geq 0}$ is a spectral ring satisfying Poincaré duality. This allows us to endow $(E^1_{p,q}(\X))_{p,q\geq 0}$ with an intersection product. Furthermore, the adjoint Borel--Viro isomorphism $(\bv^*_{p,q}:H^{p,q}(X;\F)\rightarrow E_1^{p,q}(\X))_{p,q\geq 0}$ is a bigraded algebra isomorphism. In particular, this implies, with the Leibniz rule, that we have
\begin{equation}\label{eq:auto_adj}
	[a]\cdot\partial_\varepsilon[b]=\partial_\varepsilon[a]\cdot [b],
\end{equation}
for all $[a]\in H_{p,q}(X;\F_2)$ and $[b]\in H_{n-p-1,n-q+1}(X;\F_2)$. Hence, $\partial_\varepsilon$ is self adjoint with respect to the intersection pairing; cf. \cite[Theorem 5.29]{chenal_poincare_2026}.

\section{Waves and Face Cycles}

\subsection{Face Cycles} Let $P$ be a smooth $(n+1)$-dimensional polytope of $\M{\R}$ endowed with a primitive triangulation $K$. The Tropical Lefschetz Hyperplane Section Theorem states that ${i_{p,q}:H_{p,q}(X;\F_2)\rightarrow H_{p,q}(Y;\F_2)}$ is surjective, for all $p,q\geq0$ with $p+q=n$; see \cite[Theorem 1.2]{arnal_lefschetz_2021}, \cite[Proposition 3.2]{brugalle_combinatorial_2024} or \cite[Theorem 2]{chenal_poincare-lefschetz_2025}. In the following, we introduce a generating set for the kernel of $i_{p,q}$. These generators are called face cycles and were independently used in \cite{toussaint_real_2023}. An elementary computation shows the following lemma.

\begin{lem}\label{lem:omega_closed}
	The $1$-cochain $(\sigma^1\mapsto \omega(\sigma^1))\in C^1(K;\M{2})$ is closed (and in fact even exact).
\end{lem}

\begin{rem}
	If $\dv\in\bigwedge^{n+1}\N{2}$ is the volume element of $\N{2}$ and $\sigma^1$ is an edge of $K$, the volume element of $\cF_1((\sigma^1)^*)$ is $\omega(\sigma^1)\cdot\dv$. Thus, the fundamental class $[X]$ can be rewritten as
	\begin{equation*}
		[X]=\sum_{\sigma^1\in K} (\omega(\sigma^1)\cdot\dv)\otimes (\sigma^1)^*.
	\end{equation*}
	N.b. Lemma~\ref{lem:omega_closed} confirms it is a cycle. 
\end{rem}

\begin{dfn}
	For all simplices $\sigma^p\in K$, we choose\footnote{They are many ways to do so. The most convenient one is to choose the volume element of a complementary subspace of reduction modulo $2$ of $(T\sigma)^\perp\cap N$.} ${\omega^*(\sigma^p)=u_1\wedge\cdots\wedge u_p}$ a totally decomposed $p$-vector of $\N{2}$ satisfying $\la\omega(\sigma^p),\omega^*(\sigma^p)\ra=1$.
	In addition, $\chi_{\sigma^p}\in C^0(K;\F_2)$ denotes the cochain whose value is $1$ on the vertices of $\sigma^p$ and $0$ elsewhere. 
	Finally, we denote by $\mathbf{1}_{\sigma^p}\in C^p(K;\F_2)$ the $p$-cochain whose is $1$ on $\sigma^p$ and $0$ elsewhere (n.b. if $p=0$, $\chi_{\sigma^p}=\mathbf{1}_{\sigma^p}$). 
\end{dfn}

\begin{dfn}\label{dfn:face_cycle}
	Let $p,q\geq0$ be two integers satisfying $p+q=n$. To every interior $p$-simplex $\sigma^p$ of $K$ is associated a \emph{face cycle} $f_{\sigma^p}\coloneqq \omega^*(\sigma^p)\otimes\partial(\sigma^p)^*\in C_q(X;\cF_p^X)$. See Figure~\ref{fig:face_cycle}.
\end{dfn}

Since $\cF_0^X$ is the constant cosheaf $\F_2$, $H_{0,n}(X;\F_2)=H_n(X;\F_2)$. The space $X$ is homotopic to a wedge of $n$-spheres, one for each interior vertex of $K$. The class $[f_v]$ is the class of the sphere corresponding to $v$.

\begin{lem}\label{lem:generating_family}
	For all interior simplices $\sigma^p$, $f_{\sigma^p}$ is a cycle, the class $[f_{\sigma^p}]$ does not depend on the choice of $\omega^*(\sigma^p)$, and $i_{p,n-p}[f_{\sigma^p}]=0$.
\end{lem}

\begin{proof}
	Let $\sigma^p$ be an interior simplex of $K$, and $q=n-p$. By definition, $f_{\sigma^p}\in C_q(X;\cF_p^X)$ seen as an element of $C_q(Y;\cF_p^Y)$ is the boundary of
	\begin{equation}\label{eq:boundary_face_cycle}
		b=\omega^*(\sigma^p)\otimes(\sigma^p)^* \in C_{q+1}(Y;\cF_p^Y).
	\end{equation}
	Hence, it is a cycle and $i_{p,q}[f_{\sigma^p}]$ vanishes. Also, should we choose another $\omega^*$, say $\omega^*_1$, the chain
	\begin{equation*}
		b-b_1=(\omega^*(\sigma^p)-\omega^*_1(\sigma^p))\otimes(\sigma^p)^*
	\end{equation*}
	belongs to $C_{q+1}(X;\cF_p^X)$ since $\la\omega(\sigma^p),\omega^*(\sigma^p)-\omega^*_1(\sigma^p)\ra$ vanishes by hypothesis. The boundary of $b-b_1$ is
	\begin{equation*}
		\partial(b-b_1)=f_{\sigma^p}-\omega^*_1(\sigma^p)\otimes\partial(\sigma^p)^*,
	\end{equation*}
	and therefore the class $[f_{\sigma^p}]$ does not depend on $\omega^*$. 
\end{proof}

\begin{figure}[H]
	\centering
	\begin{subfigure}[t]{.45\textwidth}
		\centering
		\input{face_cycle_triangulation_snippet}
		\caption{The closed star of a vertex $v$ and the portion of $X$ around it.}
	\end{subfigure}
	\begin{subfigure}[t]{.45\textwidth}
		\centering
		\input{face_cycle_vertex}
		\caption{The face cycle associated to $v$.}
	\end{subfigure}
	\begin{subfigure}[t]{.45\textwidth}
		\centering
		\input{face_cycle_edge}
		\caption{The face cycle associated to $\sigma^1$.}
	\end{subfigure}
	\begin{subfigure}[t]{.45\textwidth}
		\centering
		\input{face_cycle_triangle}
		\caption{The face cycle associated to $\sigma^2$.}
	\end{subfigure}
	\caption{The closed star of a vertex in a tridimensional triangulation $K$ and the face cycles associated to some adjacent simplices (n.b. the framings have to be though of as in $\bigwedge^*\N{2}$).}
	\label{fig:face_cycle}
\end{figure}

\begin{prop}\label{prop:face_cycles_kernel}
	 For all $p+q=n$, $\la[f_\sigma]:\sigma^p\textnormal{ interior simplices of }K\ra$ is the kernel of $i_{p,q}$.
\end{prop}

\begin{proof}
	The proposition is a corollary of \cite[Lemma~4.22 and Theorem~4.23]{chenal_poincare-lefschetz_2025}. They assert that both the tropical homologies of $Y$ and $X$ can be computed as the cohomology of cellular sheaves on $K$. Let $Z$ denote either $Y$ or $X$ and define, as in the proof of aforementioned Theorem 4.23, the cellular sheaves of local tropical homology
	\begin{equation*}
		\mathcal{G}^Z_p : \sigma^q \mapsto H_n(Y;Y-\sigma^q;\cF_p^Z) \quad\textnormal{for all }p\geq 0,
	\end{equation*}
	where $Y$ is endowed with its fine structure (i.e. is the cubical subdivision of $K$) and $\cF_p^X$ here means its extension by $0$ over $Y$. 
	The inclusion $i_p:\cF_p^X\rightarrow \cF_p^Y$ of cosheaves on the fine structure of $Y$ gives rise an inclusion $j_p:\mathcal{G}_p^X\rightarrow \mathcal{G}_p^Y$ of sheaves on $K$. By \cite[Theorem~3.3]{chenal_poincare-lefschetz_2025}, we have a commutative diagram
	\begin{equation*}
		\begin{tikzcd}
			H_{q+1}(Y;\coker(i_p)) \ar[r,"\partial"] \ar[d,"\cong"] & H_{p,q}(X;\F_2)\arrow[d,"\cong"] \arrow[r,"i_{p,q}"] & H_{p,q}(Y;\F_2)\arrow[d,"\cong"] \\
			H^{n-q}(K;\coker(j_p)) \ar[r,"\df" swap] & H^{n+1-q}(K;\mathcal{G}^X_p ) \arrow[r,"j_{p}^{n+1-q}" swap] & H^{n+1-q}(K;\mathcal{G}^Y_p)
		\end{tikzcd}
	\end{equation*}
	where the vertical arrows are canonically provided by \cite[Theorem~3.3]{chenal_poincare-lefschetz_2025} (n.b. $\coker(j_p)$ coincides with the sheaf $\sigma^q\mapsto H_{n+1}(Y;Y-\sigma^q;\coker(i_p))$). By \cite[Lemma 4.22]{chenal_poincare-lefschetz_2025}, for all simplices $\sigma^q\in K$, we have the exact sequence
	\begin{equation*}
		0 \rightarrow \mathcal{G}^X_p(\sigma^q) \rightarrow \mathcal{G}^Y_p(\sigma^q)\cong {\textstyle \bigwedge^{p-s}}\M{2}/\Sed(\sigma^q) \xrightarrow{\omega(\sigma^q)\cdot-} {\textstyle \bigwedge^{p-s-q}}\M{2}/\Sed(\sigma^q), 
	\end{equation*}
	where $s$ denotes the rank of $\Sed(\sigma^q)$. Thus, the cokernel sheaf of $j_p$ is supported on the $p$-skeleton of $K$. Furthermore, for all simplices $\sigma^p$, $\coker(j_p)(\sigma^p)$ is either $\F_2$ or $0$, the latter occurring if and only if $\sigma^p$ lies on the boundary of $P$. Therefore, the group of $p$-cocycles with coefficients in $\coker(j_p)$ is isomorphic to $C^p(K;\partial K;\F_2)$ and we have an exact sequence
	\begin{equation*}
		C^p(K;\partial K;\F_2) \xrightarrow{f} H_{p,q}(X;\F_2) \rightarrow H_{p,q}(Y;\F_2) \rightarrow 0,
	\end{equation*} 
	induced from the exact sequence $H^p(K;\coker(j_p)) \xrightarrow{\df~} H^{p+1}(K;\mathcal{G}_p^X) \rightarrow H^{p+1}(K;\mathcal{G}_p^Y) \rightarrow 0$.
	The morphism $f$ maps the characteristic cochain $\mathbf{1}_{\sigma^p}$ to the class $[f_{\sigma^p}]$. Indeed, under \cite[Lemma~4.22]{chenal_poincare-lefschetz_2025}, the reduction of the cochain 
	\begin{equation*}
		b=\sum_{\sigma^{n+1}\geq \sigma^p}\omega^*(\sigma^p)\otimes[\sigma^p;\sigma^{n+1}] \in \Omega_{n+1-p}(Y;\cF_p^Y),
	\end{equation*}
	in $\Omega_{n+1-p}(Y;\coker(i_p))$ is the image of the cocycle $\mathbf{1}_{\sigma^p}$ seen in $Z^p(K;\coker(j_p))$ by the quasi-isomorphism $C^*(K;\coker(j_p))\rightarrow \Omega_{n+1-*}(Y;\coker(i_p))$ of \cite[Theorem~3.3]{chenal_poincare-lefschetz_2025}. Therefore, $f(\mathbf{1}_{\sigma^p})$ is the class of $\partial b$ seen in $\Omega_{n-p}(X;\cF_p)$, i.e. $[f_{\sigma^p}]$; cf. Equation~(\ref{eq:boundary_face_cycle}).
\end{proof} 

\begin{lem}\label{lem:intersec_numb}
	Let $\upsilon$ be an interior vertex and $\sigma^n$ be an interior simplex of $K$. We have
	\begin{equation*}
		[f_\upsilon]\cdot [f_{\sigma^n}]=\sum_{\sigma^{n+1}\geq \upsilon,\sigma^n}\la \omega(\sigma^{n+1}-\upsilon),\omega^*(\sigma^n)\ra.
	\end{equation*}
\end{lem}

\begin{proof}
	One way to compute this intersection product is to find a cocycle $\alpha_v$ whose class is Poincaré dual to $[f_v]$ and then evaluate it on $f_{\sigma^n}$. Let $\dv$ be the volume element of $\N{2}$ and $\sigma^p\leq \sigma^q$ be two simplices of $K$. Following \cite[Lemma 2.3]{brugalle_combinatorial_2024}, we have
	\begin{equation*}
		\dim \cF^X_n(\sigma^p, \sigma^q)=\binom{n+1-s}{n}-\binom{n+1-s-p}{n-p},
	\end{equation*}
	where $s=\dim\Sed(\sigma^q)$. Thus, $\cF^X_{n}(\sigma^p,\sigma^q)$ vanishes whenever $\sigma^q$ lies in $\partial K$, since in this case $s\geq 1$. If $\sigma^q$ is interior, i.e. $s=0$, $\dim \cF^X_n(\sigma^p,\sigma^q)=p$. Therefore, if $\sigma^q$ is interior, we have the exact sequence
	\begin{equation}\label{seq:Fn}
		C_2(\sigma^p;\F_2)\xrightarrow{\partial} C_1(\sigma^p;\F_2) \xrightarrow{\pi_{\sigma^p,\sigma^q}} \cF^X_n(\sigma^p;\sigma^q) \rightarrow 0
	\end{equation}
	where $\pi_{\sigma^p,\sigma^q}$ sends the edge $\sigma^1\leq\sigma^p$ to $\omega(\sigma^1)\cdot\dv$. Indeed, Definition~\ref{dfn:trop_cosheaves} ensures $\pi_{\sigma^p,\sigma^q}$ is surjective and Lemma~\ref{lem:omega_closed} garanties the image of $\partial$ lies in the kernel of $\pi_{\sigma^p,\sigma^q}$. Then, the exactness follows the equality $\dim\coker(\partial) =\dim\cF^X_n(\sigma^p,\sigma^q)$. We note that the definition of $\pi_{\sigma^p,\sigma^q}$ makes the following square commute for all $\sigma^{p_1}\leq\sigma^{p_2}\leq \sigma^{q_2}\leq\sigma^{p_1}$ with $\sigma^{q_2}$ interior,
	\begin{equation}\label{eq:comm_pi}
		\begin{tikzcd}
			C_1(\sigma^{p_1};\F_2)\ar[d,"\subset" swap] \ar[r,"\pi_{\sigma^{p_1},\sigma^{q_1}}"]& \cF^X_n(\sigma^{p_1},\sigma^{q_1}) \ar[d,"\cF^X_n({[\sigma^{p_1},\sigma^{q_1}]\leftarrow [\sigma^{p_2},\sigma^{q_2}]})"] \\
			C_1(\sigma^{p_2};\F_2) \ar[r,"\pi_{\sigma^{p_2},\sigma^{q_2}}" swap]& \cF^X_n(\sigma^{p_2},\sigma^{q_2}) 
		\end{tikzcd}
	\end{equation}
	 Let us consider the cochain $\alpha_v\in \Omega^0(X;\cF_X^n)$ defined by the property
	\begin{equation}\label{eq:def_alpha}
		\pi_{\sigma^p,\sigma^p}^*(\alpha_v(\sigma^p,\sigma^p))=\df\chi_v|_{\sigma^p},\;\forall \sigma^p\in K\setminus \partial K \;\textrm{and}\;\alpha_v(\sigma^p,\sigma^p)=0,\;\forall \sigma^p\in \partial K, 
	\end{equation} 
	where $\pi^*_{\sigma^p,\sigma^p}:\cF_X^n(\sigma^p,\sigma^p)\rightarrow C^1(\sigma^p;\F_2)$ is the adjoint of $\pi_{\sigma^p,\sigma^p}$ and $\chi_v\in C^0(K;\F_2)$ is the characteristic cochain of $v$. Let us show that $\alpha_v$ is closed. If $[\sigma^p,\sigma^{p+1}]$ is a fine edge of $X$, i.e. $p\geq 1$, then
	\begin{equation}\label{eq:dalpha}
		\df\alpha_v(\sigma^p,\sigma^{p+1})=\alpha_v(\sigma^{p+1},\sigma^{p+1})-\alpha_v(\sigma^p,\sigma^p),
	\end{equation}
	where both $\alpha_v(\sigma^{p+1},\sigma^{p+1})$ and $\alpha_v(\sigma^p,\sigma^p)$ are understood as their image in $\cF^n_X(\sigma^p,\sigma^{p+1})$. There are three cases depending on whether $\sigma^p$ and $\sigma^{p+1}$ belong to $\partial K$. If $\sigma^{p+1}$ belongs to $\partial K$, so does $\sigma^p$ and both terms vanish in Equation~(\ref{eq:dalpha}). If $\sigma^p$ belongs to $\partial K$ and $\sigma^{p+1}$ does not, (\ref{eq:comm_pi}) and (\ref{eq:def_alpha}) ensure that $\pi_{\sigma^p,\sigma^{p+1}}^*(\df\alpha_v(\sigma^p,\sigma^{p+1}))=\pi_{\sigma^{p+1},\sigma^{p+1}}^*(\alpha_v(\sigma^{p+1},\sigma^{p+1}))|_{\sigma^{p}}=\df\chi_v|_{\sigma^{p}}$. It vanishes for $v$ is an interior vertex and $\sigma^{p}$ is contained in $\partial K$. The last case, where neither $\sigma^{p+1}$ nor $\sigma^p$ belong to $\partial K$, follows the computation
	\begin{align*}
		\pi_{\sigma^p,\sigma^{p+1}}^*(\df\alpha_v(\sigma^p,\sigma^{p+1}))&=\pi_{\sigma^{p+1},\sigma^{p+1}}^*(\alpha_v(\sigma^{p+1},\sigma^{p+1}))|_{\sigma^p}-\pi_{\sigma^p,\sigma^p}^*(\alpha_v(\sigma^p,\sigma^p))\\ &=\df\chi_v|_{\sigma^p}-\df\chi_v|_{\sigma^p}=0.
	\end{align*} 
	Now, the fundamental class of $X$ is represented by the cycle $\sum_{\sigma^1\in K} (\omega(\sigma^1)\cdot\dv)\otimes (\sigma^1)^*$. Following formula (\ref{eq:cap_prod}), 
we have
	\begin{equation*}
		\left(\sum_{\sigma^{n+1}\geq \sigma^1} (\omega(\sigma^1)\cdot \dv)\otimes[\sigma^1,\sigma^{n+1}]\right)\capp \alpha_v =\sum_{\sigma^{n+1}\geq \sigma^1} (\alpha_v(\sigma^1,\sigma^1)\wedge\omega(\sigma^1))\cdot \dv \otimes [ \sigma^1,\sigma^{n+1}].
	\end{equation*}
	But the definition of $\pi_{\sigma^p,\sigma^q}$ implies $\alpha\wedge\omega(\sigma^1)=\pi_{\sigma^p,\sigma^q}^*(\alpha)(\sigma^1)\dv^*$ for all $\sigma^p\leq\sigma^q$ with $\sigma^p\notin \partial K$, all edges $\sigma^1\leq\sigma^p$, and all $\alpha\in \cF_X^n(\sigma^p,\sigma^q)$. Thus,
	\begin{equation*}
		\left(\sum_{\sigma^{n+1}\geq \sigma^1} (\omega(\sigma^1)\cdot \dv)\otimes[\sigma^1,\sigma^{n+1}]\right)\capp \alpha_v=\sum_{\sigma^{n+1}\geq \sigma^1} \df\chi_v(\sigma^1)\otimes[\sigma^1,\sigma^{n+1}] =f_v,
	\end{equation*}
	and $[\alpha_v]\capp[X]=[X]\capp[\alpha_v]=[f_v]$. Therefore,
	\begin{equation*}
		[f_v]\cdot[f_{\sigma^n}]=\la[\alpha_v],[f_{\sigma^n}]\ra= \la\alpha_v,f_{\sigma^n}\ra = \sum_{\sigma^{n+1}\geq \sigma^n} \la\alpha_v(\sigma^{n+1},\sigma^{n+1}),\omega^*(\sigma^n)\ra.
	\end{equation*}
	Let $\sigma^{n+1}\in K$, since it is primitive, $\omega(\sigma^1)\wedge \sum_{v'\in \sigma^{n+1}}\chi_{v}(v')\omega(\sigma^{n+1}-v')=\df\chi_v(\sigma^1)$ for all edges $\sigma^1\leq\sigma^{n+1}$. Thus,
	\begin{equation*}
		[f_v]\cdot[f_{\sigma^n}]=\sum_{\sigma^{n+1}\geq v,\sigma^n} \la\omega(\sigma^{n+1}-v),\omega^*(\sigma^n)\ra,
	\end{equation*}
	and the proposition follows.
\end{proof}

There are three cases for computing the intersection $[f_v]\cdot[f_{\sigma^n}]$:  (1) $v$ and $\sigma^n$ are “far away”, i.e. no $(n+1)$-simplex contains both, and we have $[f_v]\cdot[f_{\sigma^n}]=0$; (2) the join $v*\sigma^n$ is an $(n+1)$-simplex of $K$ and $[f_v]\cdot[f_{\sigma^n}]=1$; and (3) $v$ is a vertex of $\sigma^{n}$. We relate the intersection number of the last case to structural numbers of the triangulation $K$.

\begin{lem}\label{lem:rho}
	Let $\sigma^n$ be an interior simplex of $K$ and $v_1,v_2$ be the apices of its two adjacent $(n+1)$-simplices. There are uniquely defined numbers $(\rho_{\sigma^n}(v))_{v\in\sigma^n}$ in $\F_2$ such that 
	\begin{equation*}
		\sum_{v\in\sigma^n} \rho_{\sigma^n}(v) = 0 \quad\mathrm{and}\quad
		\sum_{v\in\sigma^n} \rho_{\sigma^n}(v)v=v_1-v_2\;(\mod 2).
	\end{equation*}
\end{lem}

\begin{proof}
	Let $v_0$ be a vertex of ${\sigma^n}$. Since $K$ is primitive, $v_2-v_0$ equals $-(v_1-v_0)+u$ where $u$ is a vector in $T{\sigma^n}\cap M$. Moreover, the collection $(v-v_0)_{v\in {\sigma^n}\setminus\{v_0\}}$ is a basis of $T{\sigma^n}\cap M$. Thus, there are unique integers $(a_v)_{v\in {\sigma^n}\setminus\{v_0\}}$ such that 
	\begin{equation*}
		v_2-v_0=-(v_1-v_0)+\sum_{v\neq v_0}a_v(v-v_0).
	\end{equation*}
	Hence, if $v\neq v_0$, $\rho_{\sigma^n}(v)=a_v\,(\mod 2)$ and $\rho_{\sigma^n}(v_0)=\sum_{v\neq v_0} a_v\, (\mod 2)$.  
\end{proof}

\begin{prop}\label{prop:inter_rho}
	If $v$ is a vertex of ${\sigma^n}$, $[f_v]\cdot[f_{\sigma^n}]=\rho_{\sigma^n}(v)$. 
\end{prop}

\begin{proof}
	Let $(\df x_1,...,\df x_{n+1})$ be a basis of $M$ and $(\partial x_1,...,\partial x_{n+1})$ be the dual basis of $N$. Up to an affine change of coordinates, we can assume that $v=0$, ${\sigma^n}$ is the convex hull of $\{0,\df x_1,...,\df x_n\}$ and the two maximal simplices of $K$ that contain ${\sigma^n}$ are the convex hulls of $\{0,\df x_1,...,\df x_n,\df x_{n+1}\}$ and $\{0,\df x_1,...,-\df x_{n+1}+\sum_{i=1}^na_i\df x_i\}$. Thus, by Lemma~\ref{lem:intersec_numb}, we have
		\begin{equation*}
			[f_v]\cdot[f_{\sigma^n}]= \Big\la (\df x_2-\df x_1)\wedge\cdots\wedge(\df x_2-\df x_1)\wedge\sum_{i=1}^na_i\df x_i, \partial x_1\wedge\cdots\wedge\partial x_n \Big\ra.
		\end{equation*}
		But,
		\begin{equation*}
			\bigwedge_{i=2}^n(\df x_i-\df x_1)\wedge\sum_{i=1}^na_i\df x_i=\sum_{i=1}^n \bigwedge_{k\neq i}\df x_k \wedge\sum_{j=1}^na_j\df x_j= \sum_{i=1}^n a_i\df x_1\wedge \cdots \wedge\df x_{n}.
		\end{equation*}
		Thus, $[f_v]\cdot[f_{\sigma^n}]= \sum_{i=1}^n \rho_{\sigma^n}(\df x_i)=\rho_{\sigma^n}(v)$, as $\rho_{\sigma^n}(\df x_i)=a_i\,(\mod 2)$ for all $1\leq i\leq n$.
\end{proof}

\subsection{Wave Cohomology of Mikhalkin--Zharkov}\label{sec:waves}

In \cite{mikhalkin_tropical_2014}, Mikhlakin and Zharkov introduced the sheaves of \emph{tropical waves} of a tropical manifold. Here, we recall their construction with coefficients in $\F_2$ for this theory proves itself useful to represent the boundary operators of the first page of the Renaudineau--Shaw spectral sequence.

Let $K$ be a primitive triangulation of a smooth lattice polytope $P$ in $\M{\R}$ and $X$ be the complex of Definition~\ref{dfn:X}.
\begin{dfn}\label{dfn:wave_sheaf}
	Let $p\geq 0$ be an integer. The sheaf of $p$-\emph{tropical waves} $\cW_p$ of $X$ is the subsheaf of $\bigwedge^p\N{2}$ whose groups are defined by $\cW_p(\sigma^q,Q)\coloneqq \bigwedge^p (T(\sigma^q)^\perp\cap N)\otimes \F_2$ for all coarse cells $[\sigma^q,Q]$ of $X$. 
\end{dfn}

The definition directly implies that $(\cW_p)_{p\geq0}$ forms a sheaf of graded algebras. This leads to cup-products on $(\Omega^q(X;\cW_p))_{p,q\geq 0}$ and in cohomology. Since we work in characteristic $2$, the cohomological cup-product is commutative. Following \cite[\S3.2]{mikhalkin_tropical_2014}, the wave cohomology ring acts on the tropical homology and cohomology. 

Indeed, let $e^q$ be a coarse cell of $X$, $w\in \cW_{p_1}(e^q)$, and $u\in \cF^X_{p_2}(e^q)$. If $\overline{w}$ denotes the image of $w$ by the projection $\bigwedge^{p_1}\N{2}\rightarrow \bigwedge^{p_1}\N{2}/\Sed(e^q)$, Definitions~\ref{dfn:trop_cosheaves}~and~\ref{dfn:wave_sheaf} imply that $\overline{w}\wedge u$ belongs to $\cF_{p_1+p_2}^X(e^q)$. We denote this product by $w\wedge u$. Similarly, we define a product $u\wedge w$ but, since we work over $\F_2$, $w\wedge u=u\wedge w$. Moreover, if $e^{q_1}\leq e^{q_2}$, $w\in \cW_{p_1}(e^{q_1})$, and $u\in \cF_{p_2}(e^{q_2})$, it respects the projection formula 
\begin{equation*}
	f(g(w)\wedge u)=w\wedge h(u) \;\mathrm{where}\; \left\{\begin{smallmatrix}f=\cF^X_{p_1+p_2}(e^{q_2}\leftarrow e^{q_1})\\ g=\cW_{p_1}(e^{q_1}\rightarrow e^{q_2}) \hspace{.4cm}\\ h=\cF^X_{p_2}(e^{q_2}\leftarrow e^{q_1}) \hspace{.5cm}\end{smallmatrix} \right. 
\end{equation*}
see \cite[Lemma 2]{mikhalkin_tropical_2014}. Overall, this induces morphisms of sheaves $\iota: \cW_{p_1}\otimes \cF_X^{p_2}\rightarrow \cF_X^{p_2-p_1}$ which are contractions of the wedges of a wave and a framing, i.e. for all cells $e^q$ of $X$, all $w\in \cW_{p_1}(e^q)$, all $\alpha\in\cF_X^{p_2}(e^q)$, and all $u\in \cF^X_{p_2-p_1}(e^q)$, we have
\begin{equation}\label{eq:iota}
	\la \iota(w\otimes \alpha),u\ra= \la \alpha,w\wedge u\ra.
\end{equation} 
Following Definition~\ref{dfn:cup}, for all waves $w\in \Omega^{q_1}(X;\cW_{p_1})$ and tropical cochains $\alpha\in \Omega^{q_2}(X;\cF_X^{p_2})$, we can form the cup-products $w\cupp\alpha$ and $\alpha\cupp w$ in $\Omega^{q_1+q_2}(X;\cF_X^{p_2-p_1})$ (where for the second expression we used $\iota$ precomposed by the swap). Since we work over $\F_2$, we have $[\alpha]\cupp[w]=[w]\cupp[\alpha]$ in cohomology. Moreover, we have “associativity” on the side of waves, i.e. $(w_1\cupp w_2)\cupp\alpha=w_1\cupp(w_2\cupp \alpha)$ and likewise on the other side. This follows the relation $\iota(w_1\otimes\iota(w_2\otimes\alpha))=\iota((w_1\wedge w_2)\otimes \alpha)$ easily derived from Equation~(\ref{eq:iota}). 

Following Definition~\ref{dfn:contra}, for all $w\in \Omega^{q_1}(X;\cW_{p_1})$ and all $c\in \Omega_{q_2}(X;\cF^X_{p_2})$, we have cap-products $w\cap c$ and $c\cap w$ in $\Omega_{q_2-q_1}(X;\cF^X_{p_1+p_2})$ defined by the relations
\begin{equation}\label{eq:cap_wave_rel}
	\la \alpha, w\capp c\ra = \la \alpha\cupp w,c\ra \quad \mathrm{and} \quad \la \alpha,c\capp w\ra=\la w\cupp\alpha, c\ra
\end{equation}
for all $\alpha\in \Omega^{q_2-q_1}(X;\cF_X^{p_1+p_2})$. They satisfy a Leibniz rule and induce homological cap-products. These homological cap product can be equivalently defined using the relations (\ref{eq:cap_wave_rel}) in homology. These are the pairings of \cite[Proposition 9]{mikhalkin_tropical_2014}.  Since $[\alpha]\cupp[w]=[w]\cupp[\alpha]$, we also have $[w]\capp[c]=[c]\capp[w]$ in homology. Equations~(\ref{eq:iota})~and~(\ref{eq:cap_wave_rel}) allow to compute the following formulæ for the cap-products of a wave $w\in \Omega^{q_1}(X;\cW_{p_1})$ and a tropical chain $c\in \Omega_{q_2}(X;\cF^X_{p_2})$:
\begin{equation}\label{eq:cap_wave}
	(w\capp c)_{\sigma^k,\sigma^l}=\sum_{\sigma^m\geq\sigma^l} c_{\sigma^k,\sigma^m}\wedge w(\sigma^l,\sigma^m) \quad \mathrm{and} \quad (c\capp w)_{\sigma^k,\sigma^l}=\sum_{\sigma^j\leq \sigma^k} w(\sigma^j,\sigma^k)\wedge c_{\sigma^j,\sigma^l},
\end{equation}
for all $\sigma^k\leq \sigma^l$ with $l-k=q_2-q_1$, $m-l=k-j=q_1$, and $m-k=l-j=q_2$.

\begin{lem}\label{lem:cohomo_wave_bounded_cell}
	If $e^p$ is a coarse bounded cell of $X$, $H^p(e^p;\cW_p)=\cW_p(e^p)$ has dimension $1$.
\end{lem}

\begin{proof}
	Following Definition~\ref{dfn:wave_sheaf}, if $e^q$ is a coarse bounded cell dual to an interior simplex $\sigma^{n+1-q}$, $\cW_p(e^q)$ has dimension $\scriptstyle \binom{q}{p}$. Hence, the complex $(C^q(e^p,\cW_p))_{q\geq 0}$ has trivial groups in degrees $q<p$ and $H^p(e^p,\cW_p)=\cW_p(e^p)$ which has dimension $1$. 
\end{proof}

\begin{dfn}\label{dfn:wave_numbers}
	Let $e^p$ be a coarse bounded cell of $X$ and $[w]\in H^p(X;\cW_p)$. We denote by $[w]|_{e^p}\in \cW_p(e^p)$ the vector representing the restriction of $[w]$ to $e^p$; cf. Lemma~\ref{lem:cohomo_wave_bounded_cell}. In addition, we denote by $[w]|_{e^p}^{\F_2}\in \F_2$ the number\footnote{An $\F_2$-vector space of dimension $1$ is canonically isomorphic to $\F_2$.} associated to $[w]|_{e^p}$ as a vector of the line $\cW_p(e^p)$.
\end{dfn}

As a direct consequence of Definition~\ref{dfn:wave_numbers}, we have the following dictionary between the cell restriction of a wave classe and the corresponding wave number. 

\begin{prop}\label{prop:dict_wave_res_wave_numb}
	Let $\dv$ and $\dv^*$ respectively denote the volume elements of $\N{2}$ and $\M{2}$. For all waves $[w]\in H^p(X,\cW_p)$ and all coarse bounded cells $e^p$ of $X$, we have
	\begin{equation*}
		[w]|_{e^p}^{\F_2}=\la \dv^*,[w]|_{e^p}\wedge \omega^*(\sigma^{n+1-p})\ra \quad\mathrm{and}\quad [w]|_{e^p}=[w]|^{\F_2}_{e^p} \,\omega(\sigma^{n+1-p})\cdot\dv,
	\end{equation*}
	 where $e^p=(\sigma^{n+1-p})^*$. 
\end{prop}

\subsection{Hessian Wave of a Cochain}

\begin{dfn}
	Let us denote by $\cZ^1$ the cellular sheaf on $Y$ that associates the group $Z^1(\sigma^q;\F_2)$ to the coarse cell $[\sigma^q,Q]$ and whose transition morphisms are given by the restriction of cocycles. 
\end{dfn}

\begin{dfn}\label{dfn:XS}
	Let $S\leq K$ be a subcomplex whose support is a closed $(n+1)$-ball, e.g. the closed star of a simplex. We denote by $Y\cap S$ the cubical subdivision of $S$ and by $X\cap S$ the intersection of the fine subdivision of $X$ with $Y\cap S$. In addition, we denote by $Y_S\leq Y\cap S$ and $X_S\leq X\cap S$ the subcomplexes made of the dual of the interior simplices of $S$ and interior simplices of dimension at least 1 of $S$ respectively. They can be seen as subcomplexes of the coarse structure of $Y$ and $X$ respectively; see Figure~\ref{fig:SX}.  
\end{dfn}

\begin{rem}
	$X_K$ is the subcomplex of bounded cells of $X$.
\end{rem}

\begin{figure}[H]
	\centering
	\begin{subfigure}[t]{.3\textwidth}
		\centering
		\input{SX1.tex}
		\caption{The subcomplex $S$.}
	\end{subfigure}
	\hfill
	\begin{subfigure}[t]{.3\textwidth}
		\centering
		\input{SX2.tex}
		\caption{The subcomplex $Y\cap S$.}
	\end{subfigure}
	\hfill
	\begin{subfigure}[t]{.3\textwidth}
		\centering
		\input{SX3.tex}
		\caption{The subcomplex $Y_S$.}
	\end{subfigure}
	\caption{The complexes $Y\cap S$ and $Y_S$ associated to a subcomplex $S\leq K$ satisfying Definition~\ref{dfn:XS}.}
	\label{fig:SX}
\end{figure}

\begin{prop}\label{prop:exact_seq_wave}
	We have the short exact sequence of sheaves on $X$
	\begin{equation*}
			0 \rightarrow \cW_1 \longrightarrow \N{2} \overset{z}{\longrightarrow} \cZ^1|_X \rightarrow 0
	\end{equation*}
	where $z(u)$ denotes the cochain $\sigma^1 \mapsto \la \omega(\sigma^1),u\ra$ for all edges $\sigma^1\leq \sigma^p$.
\end{prop}

\begin{proof}
	Lemma~\ref{lem:omega_closed} ensures the morphism $z$ is well defined. Let $p\geq 1$ and $e^q$ be the coarse cell of $X$ indexed by $\sigma^p\leq Q$ (so that $q=\dim Q-p$). A vector $u\in \N{2}$ induces a trivial cocyle $z(u)\in \cZ^1(e^q)=Z^1(\sigma^p;\F_2)$ if an only if $u$ is orthogonal to $(T\sigma^p\cap M)\otimes \F_2$, i.e. belongs to $\cW_1(e^q)$. Hence, $\cW_1$ is the kernel of $z$. A simple dimension argument shows that $z$ is also surjective. Indeed, with the same notations, $\dim\cZ^1(e^q)=(p+1)-1=p$, $\dim \N{2}=n+1$, and $\dim\cW_1(e^q)=n+1-p$.
\end{proof}

\begin{prop}\label{prop:H0_wave}
	For all $S\leq K$ as in Definition~\ref{dfn:XS}, $H^0(X\cap S;\cW_1)$ vanishes.
\end{prop}

\begin{proof}
	The morphism $H^0(X\cap S;\cW_1)\rightarrow H^0(X\cap S;\N{2})$ is injective. Since $X\cap S$ is connected, there must be, for every $w\in H^0(X\cap S;\cW_1)$, a vector $u\in\N{2}$ such that $w(v)=u$ for every vertex $v$ of $X\cap S$. This means $u$ belongs to $\bigcap_{v\in X\cap S}\cW_1(v)$. But, by definition, if $v=(\sigma^{n+1})^*$ for a simplex $\sigma^{n+1}\in S$, then $\cW_1(v)=((T\sigma^{n+1})^\perp\cap N)\otimes\F_2=0$. Therefore, $w=0$ and $H^0(X\cap S;\cW_1)$ vanishes.  
\end{proof}

\begin{prop}\label{prop:acyclic_Z1}
	For all $S\leq K$ as in Definition~\ref{dfn:XS}, the sheaf $\cZ^1|_{X\cap S}$ is acyclic, the restriction $H^0(X\cap S;\cZ^1)\rightarrow H^0(X_S;\cZ^1)$ is an isomorphism, and $H^0(X\cap S;\cZ^1)$ is isomorphic to $Z^1(S;\F_2)$.
\end{prop}

\begin{proof}
	Since $\cZ^1$ is supported on $X$, both restriction morphisms $H^k(Y\cap S;\cZ^1)\rightarrow H^k(X\cap S;\cZ^1)$ and $H^k(Y_S;\cZ^1)\rightarrow H^k(X_S;\cZ^1)$ are isomorphisms for all integers $k\geq 0$. The sheaf $\cZ^1$ fits the following exact sequence
	\begin{equation}\label{eq:seq_C0}
		0\rightarrow \F_2\longrightarrow \mathcal{C}^0 \overset{\df}{\longrightarrow} \cZ^1 \rightarrow 0, 
	\end{equation}
	where $\mathcal{C}^0$ is the sheaf sending the coarse cell of $Y$ indexed by $\sigma^p\leq Q$ to $C^0(\sigma^p;\F_2)$ and whose restriction morphisms are given by cochain restrictions. Therefore, by the cohomological long exact sequences, we have
	\begin{equation}\label{eq:cohomoCZ}
		\begin{tikzcd}
			0\ar[r]& \F_2\ar[r]\ar[d,"\cong"]& H^0(Y\cap S;\mathcal{C}^0) \ar[r,"\df"]\ar[d]& H^0(Y\cap S;\cZ^1) \ar[r]\ar[d]& 0 \\
			0\ar[r]& \F_2\ar[r]& H^0(Y_S;\mathcal{C}^0) \ar[r,"\df"]& H^0(Y_S;\cZ^1) \ar[r]& 0 
		\end{tikzcd}
	\end{equation} 
	and $H^k(Y\cap S;\mathcal{C}^0)\cong H^k(Y\cap S;\cZ^1)$ for all $k>0$ since $S$ is contractible and $Y\cap S$ retracts on $Y_S$. In the framework of Remark~\ref{rem:cellular_sheaves_alexandrov}, $\mathcal{C}^0|_{Y\cap S}$ is a flabby sheaf\footnote{I.e. a sheaf all of whose restriction morphisms are surjective, n.b. in the cellular setting it means that for all subsets of cells $E_1\leq E_2\leq \mathrm{Cell}\,Y\cap S$ the map $\mathrm{proj.lim}_{E_2}\mathcal{C}^0\rightarrow \mathrm{proj.lim}_{E_1}\mathcal{C}^0$ is onto. Flabby sheaves are acyclic.}. Hence, $\cZ^1|_{Y\cap S}$ and $\cZ^1|_{X\cap S}$ are both acyclic. 
	
	Moreover, an element $\alpha$ in $H^0(Y\cap S;\mathcal{C}^0)$ (resp. in $H^0(Y_S;\mathcal{C}^0)$) is equivalent to the data of a cochain $\alpha_{\sigma^p}\in C^0(\sigma^p;\F_2)$  for every simplex $\sigma^p$ of $S$ (resp. interior simplex of $S$), such that if $\sigma^{p+1}\geq \sigma^p$ in $S$, we have $\alpha_{\sigma^{p+1}}|_{\sigma^p}=\alpha_{\sigma^p}$. Thus, such an $\alpha$ is fully determined by its values on the maximal simplices of $S$, i.e. its $(n+1)$-simplices. Moreover, since $S$ deprived of its $(n-1)$-skeleton is connected, there must exist a unique cochain $\beta\in C^0(S;\F_2)$ such that $\alpha_{\sigma^p}=\beta|_{\sigma^p}$ for all $\sigma^p\in S$ (resp. all interior $\sigma^p\in S$). Therefore, $H^0(Y\cap S;\mathcal{C}^0)\rightarrow H^0(Y_S;\mathcal{C}^0)$ is an isomorphism. It implies that $H^0(Y\cap S;\cZ^1)\rightarrow H^0(Y_S;\cZ^1)$ and thus $H^0(X\cap S;\cZ^1)\rightarrow H^0(X_S;\cZ^1)$ are isomorphisms. 
	
	To prove the last statement, we note that the morphism
	\begin{equation*}
		C^0(S;\F_2)\rightarrow H^0(Y\cap S;\mathcal{C}^0): \alpha\mapsto \big(\,[\sigma^p,\sigma^p]\mapsto \alpha|_{\sigma^p}\big)  
	\end{equation*}
	is an isomorphism. Therefore, the diagram (\ref{eq:cohomoCZ}), together with the fact that the restriction morphism $H^0(Y\cap S;\cZ^1)\rightarrow H^0(X\cap S;\cZ^1)$ is an isomorphism, implies that
	\begin{equation*}
		Z^1(S;\F_2)\rightarrow H^0(X\cap S;\cZ^1): \alpha\mapsto \big(\,[\sigma^p,\sigma^p] \mapsto \alpha|_{\sigma^p}\big)  
	\end{equation*} 
	is also an isomorphism. 
\end{proof}
 
\begin{cor}\label{cor:exact_seq_hess}
	For all $S\leq K$ as in Definition~\ref{dfn:XS}, we have the exact sequence
	\begin{equation*}
		0 \rightarrow \N{2} \rightarrow Z^1(S;\F_2) \xrightarrow{\df} H^1(X\cap S;\cW_1) \rightarrow H^{1}(X\cap S;\N{2}) \rightarrow 0.
	\end{equation*}	
	Moreover, the restriction $H^1(X\cap S;\cW_1)\rightarrow H^1(X_S;\cW_1)$ is an isomorphism.
\end{cor}

\begin{proof}
	The exactness of the sequence is a direct consequence of Propositions~\ref{prop:exact_seq_wave},~\ref{prop:H0_wave},~and~\ref{prop:acyclic_Z1}. To prove the last statement, we note that $X\cap S$ retracts on $X_S$, thus $H^k(X\cap S;X_S;\N{2})$ vanishes for all $k\geq 0$. Therefore, by the long exact sequence associated with Proposition~\ref{prop:exact_seq_wave}, the connecting morphism $\df:H^k(X\cap S;X_S;\cZ^1)\rightarrow H^{k+1}(X\cap S;X_S;\cW_1)$ is an isomorphism. By Proposition~\ref{prop:acyclic_Z1}, $H^0(X\cap S;\cZ^1)\rightarrow H^0(X_S;\cZ^1)$ is an isomorphism. Since, in addition, $\cZ^1|_{X\cap S}$ is acyclic, this implies that both $H^0(X\cap S;X_S;\cZ^1)$ and $H^1(X\cap S;X_S;\cZ^1)$ vanish. Therefore, the restriction $H^1(X\cap S;\cW_1)\rightarrow H^1(X_S;\cW_1)$ is an isomorphism. 
\end{proof}

\begin{rem}\label{rem:integrable}
	If $n\geq 2$, $H^1(X;\N{2})=H^{0,1}(X;\F_2)\otimes\N{2}$ vanishes and $Z^1(K;\F_2)\rightarrow H^1(X;\cW_1)$ is onto; see \cite[Lemma 5.1]{renaudineau_bounding_2023}. If $X$ is a curve, let us express the associated morphism $i^*:H^1(X;\cW_1) \rightarrow H^{0,1}(X;\F_2)\otimes \N{2}$ of Corollary~\ref{cor:exact_seq_hess}. We can easily compute that for all waves $[w]\in H^1(X;\cW_1)$ and all $[f_v]\in H_{0,1}(X;\F_2)$, 
	\begin{equation*}
		\la i^*[w],[f_v]\ra=\sum_{\sigma^1\geq v} [w]|_{(\sigma^1)^*}\in \N{2}.
	\end{equation*}
	We say that a wave $[w]$ is \emph{integrable} if $i^*[w]=0$ (n.b. this concept is only relevant for curves).
\end{rem}

\begin{dfn}\label{dfn:hessian_wave}
	Let $S$ be as in Definition~\ref{dfn:XS}. The \emph{Hessian wave} $\D^2\varepsilon\in H^1(X\cap S;\cW_1)$ of $\varepsilon\in C^0(S;\F_2)$ is its image by the composition $\D^2:C^0(S;\F_2)\rightarrow Z^1(S;\F_2)\rightarrow H^1(X\cap S;\cW_1)$.
\end{dfn}

\begin{prop}\label{prop:kernel_hessian}
	For all $S\leq K$ as in Definition~\ref{dfn:XS}, we have the exact sequence
	\begin{equation*}
		0\rightarrow \Aff(\M{2}) \longrightarrow C^0(S;\F_2) \overset{\D^2}{\longrightarrow} H^1(X\cap S;\cW_1) \longrightarrow H^1(X\cap S;\N{2}) \rightarrow 0,
	\end{equation*}
	In particular, if $e^1$ is a bounded edge of $X$ dual to a simplex $\sigma^n\in K$, then $\D^2\varepsilon|_{e^1}=0$ if and only if there is an affine function $a:\M{2}\rightarrow \F_2$ such that $\varepsilon(v)=a(v\,\mod 2)$ for every vertex of the closed star of $\sigma^n$.  
\end{prop}

\begin{proof}
	The exact sequence is an elementary reformulation of the exact sequence of Corollary~\ref{cor:exact_seq_hess} using the Hessian operator. If $S$ denotes the closed star of $\sigma^n$, then it satisfies Definition~\ref{dfn:XS}. Moreover, $X_S=e^1$ and the restriction $H^1(X\cap S;\cW_1)\rightarrow H^1(e^1;\cW_1)$ is an isomorphism. Therefore, $\D^2\varepsilon|_{e^1}=0$ if and only if $\D^2\varepsilon|_{X\cap S}=0$. By the exact sequence, this is equivalent to $\varepsilon$ being affine over $S$. 
\end{proof}

\begin{prop}\label{prop:formula_hess}
	Let $e^1=(\sigma^n)^*$ be a coarse bounded edge of $X$, i.e. $\sigma^n$ is an interior simplex of $K$. We have 
	\begin{equation*}
		\D^2\varepsilon|^{\F_2}_{e^1}=\varepsilon(v_1)+\varepsilon(v_2)+ \sum_{v\in \sigma^n} \varepsilon(v)\rho_{\sigma^n}(v),
	\end{equation*}
	where $v_1,v_2$ are the two apices of the $(n+1)$-simplices of $K$ adjacent to $\sigma^n$. 
\end{prop}

\begin{proof}
	Let $\D\varepsilon\in \Omega^0(X;\N{2})$ denote the cochain 
	\begin{equation}\label{eq:Nabla_eps}
		\D\varepsilon : [\sigma^p,\sigma^p]\mapsto \sum_{v\in \sigma^p}\varepsilon(v)\omega(\sigma^p-v)\cdot\omega^*(\sigma^p).
	\end{equation}
	 This cochain satisfies 
	 \begin{equation*}
	 	\la \omega(\sigma^1),\D\varepsilon(\sigma^p,\sigma^p)\ra = \df\varepsilon(\sigma^1),
	 \end{equation*}
	 for all edges $\sigma^1\leq \sigma^p$. Thus, it lifts $\df\varepsilon$ seen in $\Omega^0(X;\cZ^1)$ into $\Omega^0(X;\N{2})$. Given Definition~\ref{dfn:hessian_wave}, $\D^2\varepsilon$ is the class of $\df\D\varepsilon$ seen as a cocycle in $\Omega^1(X;\cW_1)$. Hence, in the notations of the statement of the proposition,
	\begin{equation*}
		\D^2\varepsilon|_{e^1} = \sum_{i=1}^2\big(\D\varepsilon(v_i*\sigma^n,v_i*\sigma^n)+ \D\varepsilon(\sigma^n,\sigma^n)\big)= \sum_{i=1}^2\D\varepsilon(v_i*\sigma^n,v_i*\sigma^n). 
	\end{equation*}
	Expanding the right hand side using (\ref{eq:Nabla_eps}), we find
	\begin{align*}
		\sum_{i=1}^2\D\varepsilon(v_i*\sigma^n,v_i*\sigma^n) &= \sum_{i=1}^2 \left( \varepsilon(v_i)\omega(\sigma^n)\cdot \dv+  \sum_{v\in \sigma^n} \varepsilon(v)\omega(v_i*\sigma^n-v)\cdot \dv\right) \\
		&= \sum_{i=1}^2 \varepsilon(v_i)\omega(\sigma^n)\cdot \dv + \sum_{v\in \sigma^n} \varepsilon(v)\big((v_1-v_2)\wedge\omega(\sigma^n-v)\big)\cdot \dv 
	\end{align*}
	Following Lemma~\ref{lem:rho}, $v_1-v_2=\sum_{v'\in \sigma^n-v}\rho_{\sigma^n}(v')(v'-v)$ and we have 
	\begin{align*}
		\sum_{v\in \sigma^n} \varepsilon(v)\big((v_1-v_2)\wedge\omega(\sigma^n-v)\big)\cdot \dv &=  \sum_{v\in \sigma^n}\sum_{v'\in \sigma^n-v} \varepsilon(v)\big(\rho_{\sigma^n}(v')(v'-v)\wedge\omega(\sigma^n-v)\big)\cdot \dv.
	\end{align*}
	Moreover, since $(v'-v)\wedge\omega(\sigma^n-v)=\omega(\sigma^n)$ for all $v\neq v'\in \sigma^{n}$,
	\begin{align*}
		\sum_{v\in \sigma^n}\sum_{v'\in \sigma^n-v} \varepsilon(v)\big(\rho_{\sigma^n}(v')(v'-v)\wedge\omega(\sigma^n-v)\big)\cdot \dv &= \sum_{v\in \sigma^n}\sum_{v'\in \sigma^n-v} \varepsilon(v)\rho_{\sigma^n}(v')\omega(\sigma^n)\cdot \dv \\
		&= \sum_{v\in \sigma^n}\varepsilon(v)\rho_{\sigma^n}(v) \omega(\sigma^n)\cdot\dv,
	\end{align*}
	for $\sum_{v\in\sigma^n}\rho_{\sigma^n}(v)=0$; see Lemma~\ref{lem:rho} again. The line $\cW_1(e^1)$ is spanned by $\omega(\sigma^n)\cdot \dv$. Thus putting everything back together we find that 
	\begin{equation*}
		\D^2\varepsilon|^{\F_2}_{e^1}=\varepsilon(v_1)+\varepsilon(v_2) + \sum_{v\in \sigma^n}\varepsilon(v)\rho_{\sigma^n}(v), 
	\end{equation*}
	and the proposition follows.
\end{proof}

As a direct consequence of the formulæ of Lemma~\ref{lem:intersec_numb} and Propositions~\ref{prop:inter_rho} and \ref{prop:formula_hess}, we find the following identity between intersection numbers and wave numbers.

\begin{cor}\label{cor:wave_int}
	For all interior vertices $v$ and interior simplices $\sigma^n$ of $K$, we have
	\begin{equation*}
		[f_v]\cdot[f_{\sigma^n}]=\D^2\chi_v|^{\F_2}_{e^1},
	\end{equation*}
	where $e^1$ is the coarse bounded edge dual to $\sigma^n$. 
\end{cor}

\begin{rem}
	Let $K$ be convex and induced as the dual subdivision of a tropical hypersurface $\mathcal{X}(\T)$ defined by the tropical polynomial $f(x)=\max_{v\in P\cap M}(a(v)+\la v,x\ra)$. Mikhalkin and Zharkov introduced the \emph{eigenwave} of $\mathcal{X}(\T)$ as a distinguished element $\phi_{\mathcal{X}(\T)}$ in $H^1(\mathcal{X}(\T);\cW_1\otimes \R)=H^1(X;\cW_1\otimes \R)$ (here $\cW_1$ has to be interpreted with integral coefficients i.e. not tensorised by $\F_2$ as in Definition~\ref{dfn:wave_sheaf}); cf. \cite[\S5.1]{mikhalkin_tropical_2014}. The exact sequence of Corollary~\ref{cor:exact_seq_hess} is still valid over $\R$ and therefore one has a second natural wave, namely the Hessian $\D^2a$. Do we have $\phi_{\mathcal{X}(\T)}=\D^2a$?
\end{rem}

\subsection{Cap-products of Waves and Face Cycles}

In this paragraph, we establish useful formulæ to compute cap products of waves with the classes of face cycles. 

\begin{lem}\label{lem:cap_wave_dual_cell}
	Let $w\in \Omega^1(X;\cW_1)$, $\sigma^p$ be a simplex of $K$ with $p>0$, and $u\in \cF_q^X((\sigma^p)^*)$. We have
	\begin{equation*}
		\big(u\otimes (\sigma^p)^*\big)\capp w=\sum_{\sigma^{p+1}\geq \sigma^p} \big(w(\sigma^p,\sigma^{p+1})\wedge u \big)\otimes (\sigma^{p+1})^*.
	\end{equation*}
\end{lem}

\begin{proof}
	We apply formula (\ref{eq:cap_wave}),
	\begin{align*}
		\big(u\otimes (\sigma^p)^*\big)\capp w&=\sum_{\sigma^{n+1}\geq \sigma^p} \big(u\otimes [\sigma^p;\sigma^{n+1}]\big)\capp w \\
		& = \sum_{\sigma^{p+1}\geq \sigma^p} \sum_{\sigma^{n+1}\geq \sigma^{p+1}} \big(w(\sigma^p,\sigma^{p+1})\wedge u\big)\otimes [\sigma^{p+1};\sigma^{n+1}] \\
		&= \sum_{\sigma^{p+1}\geq \sigma^p} \big(w(\sigma^p,\sigma^{p+1})\wedge u\big)\otimes (\sigma^{p+1})^*,
	\end{align*}
	and the proposition follows.
\end{proof}

\begin{prop}\label{prop:w_cap_face} 
	Let $[w]\in H^1(X;\cW_1)$ and $\sigma^p$ be an interior simplex of $K$ with $p>0$, we have
	\begin{equation*}
		[w]\capp [f_{\sigma^p}]= \sum_{\sigma^{p+1}\geq \sigma^p} \la \omega(\sigma^{p+1}),w(\sigma^p;\sigma^{p+1})\wedge\omega^*(\sigma^p)\ra [f_{\sigma^{p+1}}],
\end{equation*}
for any cocycle $w$ representing $[w]$. 
\end{prop}

\begin{proof} Let $\sigma^p$ be an interior simplex with $p>0$ and $w\in \Omega^1(X,\cW_1)$ be a closed wave. We have
	\begin{align*}
		f_{\sigma^p}\capp w &= \sum_{\sigma^{p+1}\geq\sigma^p} (\omega^*(\sigma^p)\otimes (\sigma^{p+1})^*)\capp w\\
		&= \sum_{\sigma^{p+1}\geq\sigma^p} \sum_{\sigma^{p+2}\geq\sigma^{+1}} \big(w(\sigma^{p+1},\sigma^{p+2})\wedge \omega^*(\sigma^p)\big)\otimes(\sigma^{p+2})^* & \textnormal{by Lemma~\ref{lem:cap_wave_dual_cell}}\\
		&= \sum_{\sigma^{p+2}>\sigma^p} \left( \sum_{\substack{\sigma^{p+1}\leq \sigma^{p+2} \\ \sigma^{p+1}\geq \sigma^p}} w(\sigma^{p+1},\sigma^{p+2})\wedge \omega^*(\sigma^p)\right)\otimes(\sigma^{p+2})^*\\
		&= \sum_{\sigma^{p+2}>\sigma^p} \left( \sum_{\substack{\sigma^{p+1}\leq \sigma^{p+2} \\ \sigma^{p+1}\geq \sigma^p}} w(\sigma^{p},\sigma^{p+1})\wedge \omega^*(\sigma^p)\right)\otimes(\sigma^{p+2})^* & \textnormal{for $w$ is closed}\\
		&= \sum_{\sigma^{p+1}\geq \sigma^p} \big(w(\sigma^p,\sigma^{p+1})\wedge \omega^*(\sigma^p)\big)\otimes \partial(\sigma^{p+1})^*.
	\end{align*}
	We can write $w(\sigma^p,\sigma^{p+1})\wedge \omega^*(\sigma^p)=\la \omega(\sigma^{p+1}), w(\sigma^p,\sigma^{p+1})\wedge \omega^*(\sigma^p)\ra \omega^*(\sigma^{p+1}) + c_{\sigma^p,\sigma^{p+1}}$ with $c_{\sigma^p,\sigma^{p+1}}\in\cF^X_{p+1}((\sigma^{p+1})^*)$. Therefore,
	\begin{equation*}
		f_{\sigma^p}\capp w = \sum_{\sigma^{p+1}\geq\sigma^p} \la \omega(\sigma^{p+1}), w(\sigma^p,\sigma^{p+1})\wedge \omega^*(\sigma^p)\ra f_{\sigma^{p+1}} +\partial\left(\sum_{\sigma^{p+1}\geq\sigma^p} c_{\sigma^p,\sigma^{p+1}}\otimes(\sigma^{p+1})^*\right),
	\end{equation*}
	with $\sum_{\sigma^{p+1}\geq\sigma^p} c_{\sigma^p,\sigma^{p+1}}\otimes(\sigma^{p+1})^*\in C_{n-p}(X;\cF^X_{p+1})$. Since $[w]\capp[f_{\sigma^p}]=[f_{\sigma^{p}}]\capp[w]$, the proposition follows.
\end{proof}

The formula of Proposition~\ref{prop:w_cap_face} cannot make sense when $p=0$ as $w(v,\sigma^1)$ is not defined ($w$ is only defined on $X$). We can however provide a nice formula in the case of Hessian waves.  

\begin{prop}\label{prop:cap_0_face_cycle}
	Let $\varepsilon\in C^0(K;\F_2)$ and $v$ be an interior vertex of $K$. We have
	\begin{equation*}
		\D^2\varepsilon\capp[f_v]=\sum_{\sigma^1\geq v} \df\varepsilon(\sigma^1)[f_{\sigma^1}].
	\end{equation*}
\end{prop}

\begin{proof}
	Let us consider the cocycle $\df\D\varepsilon$ of (\ref{eq:Nabla_eps}) that represents $\D^2\varepsilon$. By Lemma~\ref{lem:cap_wave_dual_cell}, we have
	\begin{align*}
		f_{v}\capp \df\D\varepsilon 
		&= \sum_{\sigma^{1}\geq v}\sum_{\sigma^{2}\geq \sigma^{1}} \df\D\varepsilon(\sigma^{1},\sigma^2) \otimes (\sigma^{2})^* 
		 = \sum_{\sigma^{2}\geq v} \left(\sum_{\substack{\sigma^1\leq \sigma^{2} \\ \sigma^{1}\geq v}} \df\D\varepsilon(\sigma^{1},\sigma^2)\right) \otimes (\sigma^{2})^* \\
		 &=  \sum_{\sigma^{2}\geq v} \left(\sum_{\substack{\sigma^1\leq \sigma^{2} \\ \sigma^{1}\geq v}} \D\varepsilon(\sigma^{1},\sigma^1)+\D\varepsilon(\sigma^2,\sigma^2)\right) \otimes (\sigma^{2})^* 
		 =  \sum_{\sigma^{2}\geq v} \left(\sum_{\substack{\sigma^1\leq \sigma^{2} \\ \sigma^{1}\geq v}} \D\varepsilon(\sigma^{1},\sigma^1)\right) \otimes (\sigma^{2})^*, 
	\end{align*}
	for there are two edges between $v$ and $\sigma^2$. By definition, $\D\varepsilon(\sigma^1,\sigma^1)=\df\varepsilon(\sigma^1)\omega^*(\sigma^1)$. Hence,
	\begin{align*}
		\D^2\varepsilon \capp [f_{v}]=[f_{v}]\capp \D^2\varepsilon&= \sum_{\sigma^{1}\geq v} \df\varepsilon(\sigma^1)\left[ \sum_{\sigma^{2}\geq \sigma^{1}} \omega^*(\sigma^1) \otimes (\sigma^{2})^* \right] = \sum_{\sigma^1\geq v}\df\varepsilon(\sigma^1)[f_{\sigma^1}],
	\end{align*}
	and the proposition follows.
\end{proof}

Following Remark~\ref{rem:integrable}, the only restrictions on the cap-product we are able to compute with the formula of Proposition~\ref{prop:cap_0_face_cycle} are for curves. Indeed, if $X$ is of dimension at least 2, every wave is integrable. The next two formula are key computations in understanding the number of connected components of T-hypersurfaces.

\begin{prop}\label{prop:wave_int_curves}
	If $n=1$ and $v_1,v_2$ are two interior vertices of $K$ 
	\begin{equation*}
		[f_{v_1}]\cdot\D^2\varepsilon \capp[f_{v_2}]=\left\{\begin{array}{ll} \D^2\varepsilon|_{(\sigma^1)^*}^{\F_2} & \sigma^1=v_0*v_1 \\[5pt]
		\sum_{\sigma^1\geq v_1}\D^2\varepsilon|_{(\sigma^1)^*}^{\F_2} & v_2=v_1 \\[5pt]
		0 & \textnormal{otherwise.} \end{array}\right.
	\end{equation*}
\end{prop}

\begin{proof}
	Following Proposition~\ref{prop:cap_0_face_cycle}
	\begin{align*}
		[f_{v_1}]\cdot\D^2\varepsilon\capp[f_{v_2}] &= \sum_{\sigma^1\geq v_2}\df\varepsilon(\sigma^1)[f_{v_1}]\cdot [f_{\sigma^1}].
	\end{align*}
	Let us assume we are in the third case of the proposition i.e. that no triangle contains both $v_1$ and $v_2$. Hence, no triangle contains both $v_1$ and $\sigma^1$ for any $\sigma^1\geq v_2$. Therefore, by Lemma~\ref{lem:intersec_numb},
	\begin{equation*}
		[f_{v_1}]\cdot\D^2\varepsilon\capp[f_{v_2}]= \sum_{\sigma^1\geq v_2}\df\varepsilon(\sigma^1)\,0 =0.
	\end{equation*}
	If now $v_1=v_2$, Proposition~\ref{prop:inter_rho} ensures that 
	\begin{equation*}
		[f_{v_1}]\cdot\D^2\varepsilon\capp[f_{v_2}]= \sum_{\sigma^1\geq v_1}\df\varepsilon(\sigma^1)\rho_{\sigma^1}(v_1)=\sum_{\sigma^1\geq v_1}\big(\varepsilon(v_1)\rho_{\sigma^1}(v_1)+\varepsilon(\sigma^1-v_1)\rho_{\sigma^1}(\sigma^1-v_1)\big),
	\end{equation*}
	for $\sum_{v\in \sigma^1}\rho_{\sigma^1}(v)=0$; see Lemma~\ref{lem:rho}. Thus, following Proposition~\ref{prop:formula_hess},
	\begin{equation*}
		[f_{v_1}]\cdot\D^2\varepsilon\capp[f_{v_2}]=\sum_{\sigma^1\geq v_1}\left(\D^2\varepsilon|^{\F_2}_{(\sigma^1)^*}+\sum_{\sigma^2\geq\sigma^1}\varepsilon(\sigma^2-\sigma^1)\right)=\sum_{\sigma^1\geq v_1}\D^2\varepsilon|^{\F_2}_{(\sigma^1)^*},
	\end{equation*}
	for $\sum_{\sigma^1\geq v_1} \sum_{\sigma^2\geq\sigma^1}\varepsilon(\sigma^2-\sigma^1)$ is twice the sum of the values of $\varepsilon$ on the link of $v_1$, i.e. $0$. Let us now place ourselves in the last case, where there is an edge $\sigma^1$ whose vertices are $v_1$ and $v_2$. This edge is interior and is contained in two triangles $v_3*\sigma^1$ and $v_4*\sigma^1$. Therefore, following Lemma~\ref{lem:intersec_numb} and Proposition~\ref{prop:inter_rho}, we have
	\begin{equation*}
		[f_{v_1}]\cdot\D^2\varepsilon\capp[f_{v_2}]=\df\varepsilon(\sigma^1)[f_{v_1}]\cdot [f_{\sigma^1}] + \sum_{i=3}^4\df\varepsilon(v_2*v_i)[f_{v_1}]\cdot[f_{v_2*v_i}] = \df\varepsilon(\sigma^1)\rho_{\sigma^1}(v_1)+\sum_{i=3}^4\varepsilon(v_i). 
	\end{equation*}
	Since $\rho_{\sigma^1}(v_1)=\rho_{\sigma^1}(v_2)$; cf. Lemma~\ref{lem:rho}, we finally find $[f_{v_1}]\cdot\D^2\varepsilon\capp[f_{v_2}]=\D^2\varepsilon|^{\F_2}_{(\sigma^1)^*}$ following Proposition~\ref{prop:formula_hess}.
\end{proof}	

\begin{prop}\label{prop:int_numb_wave_numb}
	Let $n\geq 2$, $v$ be an interior vertex of $K$, $\sigma^{n-1}$ be an interior simplex dual to the coarse bounded cell $e^2\in X$, and $[w]\in H^1(X;\cW_1)$. We have
	\begin{equation*}
		[f_v]\cdot[w]\capp[f_{\sigma^{n-1}}]=\big(\D^2\chi_v\cupp[w]\big)\big|^{\F_2}_{e^2}.
	\end{equation*}
	In particular, $[f_v]\cdot[w]\capp[f_{\sigma^{n-1}}]$ vanishes if no maximal simplex contains both $v$ and ${\sigma^{n-1}}$, and $[f_v]\cdot[w]\capp[f_{\sigma^{n-1}}]=[w]|^{\F_2}_{(\sigma^n)^*}$ if $v*\sigma^{n-1}=\sigma^n$ belongs to $K$. 
\end{prop}

\begin{proof} As usual, $\dv$ and $\dv^*$ respectively denotes the volume element of $\N{2}$ and $\M{2}$. We have 
	\begin{align*}
		\big(\D^2\chi_v\cupp [w] \big)\big|_{e^2} &= \big([w]\cupp \D^2\chi_v \big)\big|_{e^2} \\[5pt]
		&=\sum_{\sigma^{n+1}\geq \sigma^{n-1}} \sum_{\substack{\sigma^n\leq \sigma^{n+1}\\ \sigma^n\geq \sigma^{n-1}}} w(\sigma^{n-1},\sigma^n)\wedge \df\D\chi_v(\sigma^n,\sigma^{n+1}) & \textnormal{by (\ref{eq:Nabla_eps})} \\
		&= \sum_{\sigma^n\geq \sigma^{n-1}} w(\sigma^{n-1},\sigma^n)\wedge \D^2\chi_v|_{(\sigma^n)^*} \\
		&=  \sum_{\sigma^n\geq \sigma^{n-1}} w(\sigma^{n-1},\sigma^n)\wedge (\omega(\sigma^n)\cdot\dv)\big([f_v]\cdot[f_{\sigma^n}]\big).
	\end{align*}
	With the last step following Proposition~\ref{prop:dict_wave_res_wave_numb} and Corollary~\ref{cor:wave_int}. Hence,
	\begin{align*}
		\big(\D^2\chi_v\cupp [w] \big)\big|^{\F_2}_{e^2} &= \sum_{\sigma^n\geq \sigma^{n-1}} \la\dv^*, w(\sigma^{n-1},\sigma^n)\wedge (\omega(\sigma^n)\cdot\dv)\wedge\omega^*(\sigma^{n-1})\ra \big([f_v]\cdot[f_{\sigma^n}]\big)\\
		&= \sum_{\sigma^n\geq \sigma^{n-1}} \la\omega(\sigma^n), w(\sigma^{n-1},\sigma^n)\wedge\omega^*(\sigma^{n-1})\ra \big([f_v]\cdot[f_{\sigma^n}]\big)\\
		&= [f_v]\cdot \sum_{\sigma^n\geq \sigma^{n-1}} \la\omega(\sigma^n), w(\sigma^{n-1},\sigma^n)\wedge\omega^*(\sigma^{n-1})\ra \cdot[f_{\sigma^n}]\\
		&= [f_v]\cdot[w]\capp[f_{\sigma^{n-1}}],
	\end{align*}
	by Proposition~\ref{prop:w_cap_face}. This computation implies the equality 
	\begin{equation}\label{eq:develop_inter_wave}
		[f_v]\cdot[w]\capp[f_{\sigma^{n-1}}]=\sum_{\sigma^n\geq \sigma^{n-1}} \la\omega(\sigma^n), w(\sigma^{n-1},\sigma^n)\wedge\omega^*(\sigma^{n-1})\ra \big([f_v]\cdot[f_{\sigma^n}]\big).
	\end{equation}
	Hence, if no $(n+1)$-simplex contains both $v$ and $\sigma^{n-1}$, then no $(n+1)$-simplex could contain $v$ and $\sigma^n$ for any $\sigma^n\geq \sigma^{n-1}$. Therefore, $[f_v]\cdot[w]\capp[f_{\sigma^{n-1}}]=0$ by Lemma~\ref{lem:intersec_numb} and Equation~(\ref{eq:develop_inter_wave}).
	
	Finally, if $v*\sigma^{n-1}$ is a simplex of $K$, then it is the only one containing both $v$ and $\sigma^{n-1}$. Moreover, by what we just proved, (\ref{eq:Nabla_eps}), and Proposition~\ref{prop:dict_wave_res_wave_numb}, we have
	\begin{align*}
		[f_v]\cdot[w]\capp[f_{\sigma^{n-1}}] &= (\D^2\chi_v\cupp[w])|^{\F_2}_{e^2} \\[5pt]
		&= \la \dv^*, (\df\D\chi_v\cupp w)(e^2)\wedge \omega^*(\sigma^{n-1})\ra \\
		&= \sum_{\sigma^{n+1}\geq \sigma^n\geq \sigma^{n-1}} \la \dv^*, \df\D\chi_v(\sigma^{n-1},\sigma^n)\wedge w(\sigma^n,\sigma^{n+1})\wedge \omega^*(\sigma^{n-1})\ra.
	\end{align*}
	But we have
	\begin{equation*}
		\df\D\chi_v(\sigma^{n-1},\sigma^n)=\left\{ \begin{array}{ll} 0 & \mathrm{if}\;\sigma^n\neq \sigma^{n-1}*v \\ \omega(\sigma^{n-1})\cdot\omega^*(\sigma^n) & \mathrm{if}\;\sigma^n=\sigma^{n-1}*v. \end{array} \right.
	\end{equation*}
	Thus, if $\sigma^n=\sigma^{n-1}*v$, we have $[f_v]\cdot[w]\capp[f_{\sigma^{n-1}}]=\la\dv^*, (\omega(\sigma^{n-1})\cdot\omega^*(\sigma^n))\wedge w((\sigma^n)^*)\wedge \omega^*(\sigma^{n-1}) \ra$, since $\sum_{\sigma^{n+1}\geq\sigma^n}w(\sigma^n,\sigma^{n+1})=w((\sigma^n)^*)$. Note that the expression is independent from the choice of $\omega^*$; cf. Lemma~\ref{lem:generating_family}. Moreover, there is a choice of $\omega^*(\sigma^n)$ and $\omega^*(\sigma^{n-1})$ such that $\omega^*(\sigma^n)=(\omega(\sigma^{n-1})\cdot\omega^*(\sigma^n))\wedge\omega^*(\sigma^{n-1}) $, so we may as well assume it is the case. Therefore,
	\begin{equation*}
		[f_v]\cdot[w]\capp[f_{\sigma^{n-1}}]=\la\dv^*, w((\sigma^n)^*)\wedge \omega^*(\sigma^{n}) \ra= \la\dv^*, [w]|_{(\sigma^n)*}\wedge \omega^*(\sigma^{n}) \ra=[w]|^{\F_2}_{(\sigma^n)^*},
	\end{equation*} 
	and the proposition follows.
\end{proof}

\subsection{Uniformity}

\begin{dfn}\label{dfn:X0} The subcomplex $Y_0\leq Y$ is defined as the union of the bounded maximal cells of $Y$, i.e. $Y_0\coloneqq \bigcup_{v\notin \partial K} v^*$. The subcomplex $X_0$ is the intersection $X\cap Y_0$. By definition, all cells of $X_0$ are bounded coarse cells of $X$. We also define $X_1\leq X_0$ as follows: a coarse cell $e^q\in X_0$ belongs to $X_1$ if and only if its dual simplex $\sigma^{n+1-q}\in K$ has an interior vertex $v$ whose opposite face $\sigma^{n+1-q}-v$ is also interior. Figure~\ref{fig:X0} depicts a bidimensional example.  
\end{dfn}

\begin{dfn}\label{dfn:rho_unif}
	The triangulation $K$ is said to be \emph{$\rho$-uniform} if there is a wave $\rho\in H^1(X_0;\cW_1)$ such that $\rho|_{\partial v^*}=\D^2\chi_v|_{\partial v^*}$ for every interior vertex $v\in K$. By definition, if it exists, $\rho$ is unique. 
\end{dfn}

\begin{prop}\label{prop:crit_rho_unif}
	$K$ is $\rho$-uniform if and only if for all interior simplices $\sigma^n$,
	\begin{equation*}
		\rho_{\sigma^n}:\{\textnormal{vertices of }\sigma^n\}\rightarrow \F_2
	\end{equation*}
	is constant on the subset of vertices of $\sigma^n$ contained in the interior of $K$ (n.b. the condition is vacuously true when every vertex of $\sigma^n$ is contained in the boundary of $K$). 
\end{prop}

\begin{proof}
	Let us first prove the necessity of the proposition. Let $\rho\in H^1(X_0;\cW_1)$ be a class satisfying the requirements of Definition~\ref{dfn:rho_unif} and $\sigma^n$ be an interior simplex of $K$ with two interior vertices $v_1,v_2$. Following Corollary~\ref{cor:wave_int} and Proposition~\ref{prop:inter_rho}, we have
	\begin{align*}
		\rho_{\sigma^n}(v_1)=[f_{v_1}]\cdot[f_{\sigma^n}]=\D^2\chi_{v_1}|^{\F_2}_{(\sigma^n)^*}=\D^2\chi_{v_1}|_{\partial v_1^*}|^{\F_2}_{(\sigma^n)^*}=\rho|_{\partial v_1^*}|^{\F_2}_{(\sigma^n)^*}=\rho|^{\F_2}_{(\sigma^n)^*}.
	\end{align*}
	Since the same is true for $v_2$, we have $\rho_{\sigma^n}(v_1)=\rho_{\sigma^n}(v_2)$ and the condition is necessary. Conversely let us assume the condition satisfied. We consider $\rho\in C^1(X_0;\cW_1)$ defined by 
	\begin{equation*}
		\rho((\sigma^n)^*)\coloneqq \rho_{\sigma^n}(v)\,\omega(\sigma^n)\cdot \dv,
	\end{equation*}
	where $\dv$ is the volume element of $\N{2}$ and $v$ is any vertex of $\sigma^n\setminus\partial K$. This is well defined by hypothesis and the definition of $X_0$. Let us show this is a cocyle. If $n=1$, there is nothing to prove, so we assume $n\geq 2$. We consider $e^2$ a coarse cell of $X_0$. It is dual to a simplex $\sigma^{n-1}\in K$ possessing at least an interior vertex $v$. We have
	\begin{align*}
		\df\rho(e^2)&= \sum_{\sigma^n\geq \sigma^{n-1}}\rho((\sigma^n)^*)  = \sum_{\sigma^n\geq \sigma^{n-1}}\rho_{\sigma^n}(v)\,\omega(\sigma^n)\cdot \dv.
	\end{align*}
	Following Corollary~\ref{cor:wave_int} and Proposition~\ref{prop:inter_rho} again, $\rho((\sigma^n)^*)=\rho_{\sigma^n}(v)\,\omega(\sigma^n)\cdot \dv= \D^2\chi_v|_{(\sigma^n)^*}$. By (\ref{eq:Nabla_eps}), $\df\D\chi_v$ is a cocycle representing the class $\D^2\chi_v$, thus $\rho_{\sigma^n}(v)\,\omega(\sigma^n)\cdot \dv=\df\D\chi_v((\sigma^n)^*)$. Therefore, we have
	\begin{align*}
		\df\rho(e^2)&= \sum_{\sigma^n\geq \sigma^{n-1}}  \df\D\chi_v((\sigma^{n})^*)=\df^2\D\chi_v (e^2)=0.
	\end{align*}
	Furthermore, the property $\rho((\sigma^n)^*)= \D^2\chi_v|_{(\sigma^n)^*}$ (also valid when $n=1$) establishes that if $e^1$ is a coarse cell of $X$ in $\partial v^*$, for $v$ an interior vertex of $K$, then $\rho(e^1)=\D^2\chi_v|_{e^1}$. Hence, the class of $\rho$ restricted to $\partial v^*$ equals $\D^2\chi_v|_{\partial v^*}$ and the condition is sufficient.
\end{proof}

\begin{figure}[H]
	\centering
	\begin{subfigure}[t]{.45\textwidth}
		\centering
		\input{Y01.tex}
		\caption{The primitive triangulation K.}
	\end{subfigure}
	\hfill
	\begin{subfigure}[t]{.45\textwidth}
		\centering
		\input{Y02.tex}
		\caption{The subcomplex $Y_0$ is made of the grey cells while $X_1$ is the bold egde.} 
	\end{subfigure}
	\caption{A triangulation $K$ and its dual subdivision $Y$ with the subcomplexes of Definition~\ref{dfn:X0}.}
	\label{fig:X0}
\end{figure}

\begin{prop}\label{prop:rho_int_curves}
	If $n=1$, there is a canonical integrable wave $\tilde{\rho}\in H^1(X;\cW_1)$ that satisfies $$\tilde{\rho}|_{\partial v^*}=\D^2\chi_v|_{\partial v^*}$$ for all interior vertices $v$. Therefore, every $2$-dimensional primitive triangulation is $\rho$-uniform.
\end{prop}

\begin{proof}
	Following Lemma~\ref{lem:rho}, $\rho_{\sigma^1}(v_1)=\rho_{\sigma^1}(v_2)$ for every interior edge $\sigma^1=v_1*v_2$. This allows us to define the wave $\tilde{\rho}\in H^1(X,\cW_1)$ as the class of the cocycle
	\begin{equation*}
		e^1\mapsto\left\{\begin{array}{ll} \rho_{\sigma^1}(v)\,\omega(\sigma^1)\cdot\dv & \textnormal{if $e^1=(\sigma^1)^*$ where $\sigma^1$ is interior and $v\in\sigma^1$} \\
		0 & \textnormal{otherwise}. \end{array} \right.
	\end{equation*}
	N.b. $X$ is a curve so every $1$-cochain is a cocycle. Hence, following the proof of Proposition~\ref{prop:crit_rho_unif}, we know that $\tilde{\rho}|_{\partial v^*}=\D^2\chi_v|_{\partial v^*}$ for all interior vertices $v$ of $K$. Thus, it only remains to show that it is integrable. Following Remark~\ref{rem:integrable}, we have established that $\tilde{\rho}$ is integrable if and only if $i^*(\tilde{\rho})$ vanishes where $i^*:H^1(X,\cW_1)\rightarrow H^1(X,\N{2})$ comes from the inclusion $i:\cW_1\subset \N{2}$. If $v$ is an interior vertex of $K$, we have
	\begin{equation*}
		\la i^*(\tilde{\rho}),[f_v]\ra = \sum_{\sigma^1\geq v} \tilde{\rho}|_{(\sigma^1)^*} = \sum_{\sigma^1\geq v} \D^2\chi_v|_{(\sigma^1)^*} = \la i^*(\D^2\chi_v),[f_v]\ra =0
	\end{equation*}
	for $\D^2\chi_v$ is integrable. Since this is true for all $v$, $i^*(\tilde{\rho})=0$ and $\tilde{\rho}$ is integrable. 
\end{proof}

\section{Boundaries of the First Page}

The aim of this section is to give an effective computation of the first page of the Renaudineau--Shaw spectral sequence.

\subsection{Lifting Face Cycles} 

Let $\sigma^{n+1}$ be a simplex of $K$. Recall we have an isomorphism $\N{2}\rightarrow Z^1(\sigma^{n+1};\F_2)$ sending $u\in\N{2}$ to the cocycle $\sigma^1\mapsto \la\omega(\sigma^1),u\ra$; see Proposition~\ref{prop:exact_seq_wave}. 

\begin{dfn}
	For all simplices $\sigma^p<\sigma^{n+1}$, we denote by $\oo(\sigma^p,\sigma^{n+1})$ the unique vector of $\N{2}$ satisfying $\la\omega(\sigma^1), \oo(\sigma^p,\sigma^{n+1})\ra = \df(\varepsilon + \chi_{\sigma^p})(\sigma^1)$ for all edges $\sigma^1$ of $\sigma^{n+1}$.
\end{dfn}

\begin{lem}\label{lem:origin}
	Let $\sigma^p\leq\sigma^{p+1}\leq\sigma^{n+1}$ be three interior simplices of $K$. If $U\leq \N{2}$ is a supplementary subspace of $((T\sigma^p)^\perp\cap N)\otimes\F_2$, then $\oo(\sigma^p,\sigma^{n+1})+U\subset \Arg(\sigma^{p+1},\sigma^{n+1})$. 
\end{lem}

\begin{proof}
	Following Definition~\ref{dfn:t-hyp}, $\Arg(\sigma^{p+1},\sigma^{n+1})$ is the set of vectors $u\in \N{2}$ for which there is an edge $\sigma^1$ of $\sigma^{p+1}$ such that $\la\omega(\sigma^1),u\ra\neq\df\varepsilon(\sigma^1)$. If $\sigma^1$ is an edge of $\sigma^{p+1}$ not contained in $\sigma^p$, we have
	$$\la \omega(\sigma^1), \oo(\sigma^p,\sigma^{n+1})\ra= \df\varepsilon(\sigma^1) + \df\chi_{\sigma^p}(\sigma^1) = \df\varepsilon(\sigma^1) + 1, $$
	and $\oo(\sigma^p,\sigma^{n+1})\in \Arg(\sigma^{p+1},\sigma^{n+1})$. Let us consider $U$ a supplementary subspace of $((T\sigma^p)^\perp\cap N)\otimes\F_2$. If $\sigma^p$ is a vertex, $U=0$ and there is nothing more to prove. Let us assume $p\geq 1$. If $u\in U\setminus 0$, $u$ is not orthogonal to the tangent space of $\sigma^p$. Therefore, there is an edge $\sigma^1\leq \sigma^p$ for which $\la \omega(\sigma^1),u\ra=1$. Thus,
	$$\la \omega(\sigma^1), \oo(\sigma^p,\sigma^{n+1})+u\ra= \df\varepsilon(\sigma^1) + \df\chi_{\sigma^p}(\sigma^1) + \la\omega(\sigma^1),u\ra = \df\varepsilon(\sigma^1) + 0 + 1, $$
	and $\oo(\sigma^p,\sigma^{n+1})+u\in \Arg(\sigma^{p+1},\sigma^{n+1})$ for $\sigma^1$ is also an edge of $\sigma^{p+1}$.
\end{proof}

\begin{lem}
	Let $\sigma^p<\sigma^n$ be two interior simplices of $K$, we have
	$$\sum_{\sigma^{n+1}\geq \sigma^n} \oo(\sigma^p;\sigma^{n+1})=\D^2(\varepsilon+\chi_{\sigma^p})|_{(\sigma^n)^*}.$$ 
\end{lem}

\begin{proof}
	The definition of $\oo(\sigma^p,\sigma^{n+1})$ directly implies the left-hand-side sum $s$ is orthogonal to $\sigma^n$. Let $v_1,v_2$ be the apices of the two $(n+1)$-simplices adjacent to $\sigma^n$. For every vertex $v_0\in \sigma^n$, we have
	\begin{align*}
		\la v_1-v_0,s\ra &= \la v_1-v_0,\oo(\sigma^p,v_1*\sigma^n)+\oo(\sigma^p,v_2*\sigma^n)\ra \\
		&= \df\varepsilon(v_1,v_0)+1+ \la v_1-v_0,\oo(\sigma^p,v_2*\sigma^n)\ra \\
		&= \df\varepsilon(v_1,v_0)+1+ \df\varepsilon(v_2,v_0)+1 + \la v_2-v_1,\oo(\sigma^p,v_2*\sigma^n)\ra \\
		&=\varepsilon(v_1)+\varepsilon(v_2)+\sum_{v\in\sigma^n-v_0}\rho_{\sigma^n}(v)\df(\varepsilon+\chi_{\sigma^p})(v-v_0) & \textnormal{by Lemma~\ref{lem:rho}}\\
		&=\varepsilon(v_1)+\varepsilon(v_2)+\sum_{v\in\sigma^n}\rho_{\sigma^n}(v)(\varepsilon+\chi_{\sigma^p})(v) & \textnormal{by Lemma~\ref{lem:rho}} \\ 
		&= \la v_1-v_0,\D^2(\varepsilon+\chi_{\sigma^p})|_{(\sigma^n)^*}\ra & \textnormal{by Proposition~\ref{prop:formula_hess}.}
	\end{align*}
	Since $\{v-v_0:v\in (v_1*\sigma^n)-v_0\}$ is a basis of $\N{2}$ and $\D^2(\varepsilon+\chi_{\sigma^p})|_{(\sigma^n)^*}$ is orthogonal to $\sigma^n$, the lemma follows. 
\end{proof}

\begin{dfn}
	Let $\sigma^p$ be an interior simplex of $K$, and $\tilde{f}_{\sigma^p}\in \Omega_{n-p}(X;\K^\X_{p})$ be the chain defined by
	\begin{equation*}
		\tilde{f}_{\sigma^p} \coloneqq \sum_{\sigma^{p+1}\geq\sigma^p} \sum_{\sigma^{n+1}\geq\sigma^{p+1}}\x^{\oo(\sigma^p,\sigma^{n+1})} \prod_{i=1}^p(1+\x^{u_i})\otimes[\sigma^{p+1},\sigma^{n+1}],
	\end{equation*}
	where $\omega^*(\sigma^p)=u_1\wedge\cdots\wedge u_p$. Since $\la u_1,\dots,u_p\ra$ is a supplementary subspace of $((T\sigma^p)^\perp\cap N)\otimes\F_2$, Lemma~\ref{lem:origin} garanties it is well defined. 
\end{dfn}

\begin{lem}\label{lem:lift}
	The chain $\tilde{f}_{\sigma^p}$ is a Borel--Viro lift of $f_{\sigma^p}$. 
\end{lem}

\begin{proof}
	By Proposition~\ref{prop:borel_viro}, we have
	$$\bv_p(\tilde{f}_{\sigma^p})=\sum_{\sigma^{p+1}\geq\sigma^p}\sum_{\sigma^{n+1}\geq\sigma^{p+1}} \bigwedge_{i=1}^pu_i\otimes [\sigma^{p+1},\sigma^{n+1}]=\omega^*(\sigma^p)\otimes\partial(\sigma^p)^*=f_{\sigma^p},$$
	and the proposition follows.
\end{proof}

\begin{prop}\label{prop:partial_eps_as_wave}
	For all interior simplices $\sigma^p$, we have $\partial_\varepsilon [f_{\sigma^p}]=\D^2(\varepsilon+\chi_{\sigma^p})\capp[f_{\sigma^p}]$.
\end{prop}

\begin{proof}
	By the defining property of spectral sequences, $\partial_\varepsilon[f_{\sigma^p}]$ is the homology class represented by $\bv_{p+1}(\partial \tilde{f}_{\sigma^p})$. For all $\sigma^a\leq\sigma^b$ with $b-a=n-(p+1)$,
	\begin{equation*}
		(\partial \tilde{f}_{\sigma^p})_{\sigma^a,\sigma^b} =\sum_{\sigma^{a-1}\leq \sigma^a} (\tilde{f}_{\sigma^p})_{\sigma^{a-1},\sigma^b} +\sum_{\sigma^{b+1}\geq \sigma^b} (\tilde{f}_{\sigma^p})_{\sigma^a,\sigma^{b+1}}
	\end{equation*}
	This vanishes if $\sigma^a$ does not contain $\sigma^p$ or $b<n$. If $\sigma^a\geq\sigma^p$ and $b=n+1$, i.e. $a=p+2$, the second term vanishes and
	\begin{align*}
		(\partial \tilde{f}_{\sigma^p})_{\sigma^{p+2},\sigma^{n+1}} &=\sum_{\sigma^{p+2}\geq \sigma^{p+1}\geq \sigma^p} \x^{\oo(\sigma^p,\sigma^{n+1})} \prod_{i=1}^p(1+\x^{v_i}) = 2\x^{\oo(\sigma^p,\sigma^{n+1})} \prod_{i=1}^p(1+\x^{v_i}).
	\end{align*}
	And $(\partial \tilde{f}_{\sigma^p})_{\sigma^{p+2},\sigma^{n+1}}=0$ for there are exactly two $(p+1)$-dimensional simplices between $\sigma^p$ and $\sigma^{p+2}$. If now $b=n$, the first term vanishes and
	\begin{align*}
		(\partial \tilde{f}_{\sigma^p})_{\sigma^{p+1},\sigma^n} &=\sum_{\sigma^{n+1}\geq \sigma^n} \x^{\oo(\sigma^p,\sigma^{n+1})} \prod_{i=1}^p(1+\x^{v_i}).
	\end{align*}
	Thus, applying the Borel--Viro morphism, the class $\partial_\varepsilon[f_\sigma]$ is represented by the cycle
	\begin{align*}
		\sum_{\sigma^n\geq \sigma^{p+1}\geq \sigma^p} \omega^*(\sigma^p)\wedge\sum_{\sigma^{n+1}\geq\sigma^n} \oo(\sigma^p;\sigma^{n+1})&\otimes [\sigma^{p+1};\sigma^n]\\&=\sum_{\sigma^n\geq \sigma^{p+1}\geq \sigma^p} \omega^*(\sigma^p)\wedge\D^2(\varepsilon+\chi_{\sigma^p})|_{(\sigma^n)^*} \otimes [\sigma^{p+1};\sigma^n].
	\end{align*}
	Let now $w$ be a cocycle representing the wave $\D^2(\varepsilon+\chi_{\sigma^p})$. By (\ref{eq:cap_wave}), we have
	\begin{align*}
		w\capp f_{\sigma^p} &= \sum_{\sigma^n\geq \sigma^{p+1}\geq \sigma^p}\omega^*(\sigma^p)\wedge\sum_{\sigma^{n+1}\geq\sigma^n} w(\sigma^n;\sigma^{n+1})\otimes [\sigma^{p+1};\sigma^n] \\ &= \sum_{\sigma^n\geq \sigma^{p+1}\geq \sigma^p} \omega^*(\sigma^p)\wedge\D^2(\varepsilon+\chi_{\sigma^p})|_{(\sigma^n)^*}\otimes [\sigma^{p+1};\sigma^n],
	\end{align*}
	and the proposition follows. 
\end{proof} 

\begin{lem}\label{lem:auto_adj_wave}
	Let $\varepsilon$ be a sign distribution and $\sigma^p,\sigma^q$ be two interior simplices with $p+q=n-1$, we have $(\D^2\varepsilon\capp[f_{\sigma^p}])\cdot [f_{\sigma^q}]=[f_{\sigma^p}]\cdot\D^2\varepsilon\capp[f_{\sigma^q}]$.
\end{lem}

\begin{proof}
	By Proposition~\ref{prop:partial_eps_as_wave}, $\D^2\varepsilon\capp[f_{\sigma^p}]=\partial_0[f_{\sigma^p}] + \partial_\varepsilon[f_{\sigma^p}]$ and $\D^2\varepsilon\capp[f_{\sigma^q}]=\partial_0[f_{\sigma^q}] + \partial_\varepsilon[f_{\sigma^q}]$. Therefore, by Equation~(\ref{eq:auto_adj})
	\begin{equation*}
		(\D^2\varepsilon\capp[f_{\sigma^p}])\cdot [f_{\sigma^q}]=((\partial_0+\partial_\varepsilon)[f_{\sigma^p}])\cdot [f_{\sigma^q}]=[f_{\sigma^p}]\cdot ((\partial_0+\partial_\varepsilon)[f_{\sigma^q}])=[f_{\sigma^p}]\cdot\D^2\varepsilon\capp[f_{\sigma^q}],
	\end{equation*}
	and the lemma follows.
\end{proof}

\subsection{Intersection Matrices} We use the tropical intersection product to  represent the boundary operators of the first page as matrices and provide formul\ae\,for the dimensions of the groups of the second page of the Renaudineau--Shaw spectral sequence. 

\begin{dfn}\label{dfn:inter_mat}
	Let $\varepsilon$ be a sign distribution on $K$ and $p,q\geq 0$ be two integers with $p+q=n-1$, $M_{p,q}(\varepsilon)$ denotes the matrix whose rows are indexed by the interior $p$-simplices of $K$, whose columns are indexed by interior $q$-simplices of $K$, and whose $(\sigma^p,\sigma^q)$ coefficient is $[f_{\sigma^q}]\cdot\partial_\varepsilon [f_{\sigma^p}]$. Its rank is denoted by $m_{p,q}(\varepsilon)$. 
\end{dfn}

Poincaré duality, cf. \cite[Theorem 5.29]{chenal_poincare_2026}, implies that $M_{p,q}(\varepsilon)$ and $M_{q,p}(\varepsilon)$ are transpose of each other. In particular, $m_{p,q}(\varepsilon)=m_{q,p}(\varepsilon)$ for all $p+q=n-1$.

\begin{dfn}\label{dfn:Radical}
	We assume that $n=2k$. The \emph{radical} of $X$, denoted by $\Rad(X)$, is the radical of the restriction of the intersection form to the kernel of ${i_{k,k}:H_{k,k}(X;\F_2)\rightarrow H_{k,k}(Y;\F_2)}$, i.e.
	\begin{equation*}
		\Rad(X)\coloneqq\big\{[a]\in H_{k,k}(X;\F_2)\;|\;i_{k,k}[a]=0 \textnormal{ and } [a]\cdot[b]=0,\,\forall [b]\in \ker(i_{k,k})\big\}.
	\end{equation*}
	We denote its dimension by $\rad(X)$.
\end{dfn}

\begin{prop}
	We assume that $n=2k$. Let $\omega_X\in H^{1,1}(Y;\F_2)$ be Poincaré dual to the image of the fundamental class of $X$. The radical $\Rad(X)$ is generated by the classes $\alpha|_X\capp [X]$ for all $\alpha\in H^{k,k}(Y;\F_2)$ such that $\alpha\cupp\omega_X=0$.
\end{prop}

\begin{proof}
	The Poincaré duality in $X$ provides an isomorphism between the group $H_{k,k}(X;\F_2)$ and its dual $H^{k,k}(X;\F_2)$. Under this isomorphism, the orthogonal of $\ker i_{k,k}$ is the image of the restriction $i^{k,k}$ i.e. $H^{k,k}(Y;\F_2)|_X$. Moreover $i_{k,k} (\alpha|_X\capp[X])=(\alpha\cupp\omega_X)\capp[Y]$. Therefore,
	\begin{equation*}
		\Rad(X)=\ker(i_{k,k})^\perp\cap\ker(i_{k,k})=\Big\langle \alpha|_X\capp[X] \colon \alpha\in H^{k,k}(Y;\F_2) \textnormal{ s.t. } \alpha\cup\omega_X=0\Big\rangle,
	\end{equation*}
	and the proposition follows.
\end{proof}

\begin{prop}\label{prop:rank_M}
	Let $\varepsilon\in C^0(K;\F_2)$ and $p,q\geq 0$ be integers with $p+q=n-1$. If $n$ is odd, $m_{p,q}(\varepsilon)$ equals the rank of $\partial_\varepsilon:H_{q,n-q}(X;\F_2)\rightarrow H_{n-p,p}(X;\F_2)$. If $n=2k$, it is also the case except when $p$ or $q$ equals $k$. In one case, $m_{k,k-1}(\varepsilon)$ is the rank of 
	\begin{equation*}
		H_{k-1,k+1}(X;\F_2) \overset{\partial_\varepsilon}{\longrightarrow} H_{k,k}(X;\F_2) \longrightarrow H_{k,k}(X;\F_2)/\Rad(X).
	\end{equation*}
	In the other, $m_{k-1,k}(\varepsilon)$ is the rank of the restriction of $\partial_\varepsilon$ to the kernel of $i_{k,k}$.
\end{prop}

\begin{proof}
	Following Proposition~\ref{prop:face_cycles_kernel}, $m_{p,q}(\varepsilon)$ is the rank of the morphism
	\begin{equation}\label{eq:M_replace}
		\ker(i_{n,n-p})\overset{\subset}{\longrightarrow} H_{p,n-p}(X;\F_2)\overset{\partial_\varepsilon}{\longrightarrow} H_{n-q,q}(X;\F_2) \longrightarrow H_{n-q,q}(X;\F_2)/\ker(i_{q,n-q})^\perp,
	\end{equation}
	where the orthogonal is taken with respect to the intersection pairing. Furthermore, by the Tropical Lefschetz Hyperplane Section Theorem \cite[Theorem 1.2]{arnal_lefschetz_2021}, \cite[Proposition 3.2]{brugalle_combinatorial_2024} or \cite[Theorem 2]{chenal_poincare-lefschetz_2025}, $\ker(i_{p,n-p})=H_{p,n-p}(X;\F_2)$ unless $n=2p$. Therefore, every statement, except the one about $m_{k-1,k}(\varepsilon)$, follows from (\ref{eq:M_replace}). For the remaining case, note that $\partial_\varepsilon H_{k-1,k+1}(X;\F_2)$ is contained in $\ker(i_{k,k})$. Therefore, 
	\begin{equation*}
		H_{k-1,k}(X;\F_2)\overset{\partial_\varepsilon}{\longrightarrow} H_{k,k}(X;\F_2) \longrightarrow H_{k,k}(X;\F_2)/\ker(i_{k,k})^\perp,
	\end{equation*}
	has the same rank as $H_{k-1,k}(X;\F_2)\overset{\partial_\varepsilon}{\longrightarrow} H_{k,k}(X;\F_2) \longrightarrow H_{k,k}(X;\F_2)/\Rad(X)$ and the proposition follows. 
\end{proof}

\begin{prop}\label{prop:dimE2}
Let $\varepsilon\in C^0(K;\F_2)$ be a sign distribution. If $n$ is odd, then for all non-negative integers $p+q=n$, we have
			\begin{equation*}
				\dim E^2_{p,q}(\X) = h^{p,q}(X) - \big(m_{q-1,p}(\varepsilon)+ m_{q,p-1}(\varepsilon)\big). 
			\end{equation*}
	If $n=2k$, then for all non-negative integers $p+q=n$, we have
	\begin{equation*}
		\dim E^2_{p,q}(\X)= h^{p,q}(X)-(m_{p-1,q}(\varepsilon)+m_{p,q-1}(\varepsilon))-\left\{\begin{array}{ll}
			\delta & p=k\pm1 \\
			2\delta & p=k\\
			0 & \mathrm{Otherwise,}
		\end{array}\right.
	\end{equation*}
	where $\delta$ is the dimension of the cokernel of ${i_k:H_k(\X;\F_2)\rightarrow H_k(\R P;\F_2)}$. Moreover, $\delta$ cannot exceed $\rad(X)$.
\end{prop}

\begin{proof}
	Let $p+q=n$, we have $\dim E^2_{p,q}(\X)=h^{p,q}(X)-\rk \partial^1_{p,q}-\rk\partial^1_{p+1,q-1}$. Thus, the statement follows Proposition~\ref{prop:rank_M} for $n$ odd or $n=2k$ and $p$ different from $k-1$, $k$, or $k+1$. Indeed, remember that under the Borel--Viro isomorphisms, cf. Definition~\ref{dfn:partial_varepsilon}, the boundary operator $\partial^1_{p,q}:E^1_{p,q}(\X)\rightarrow E^1_{p+1,q-1}(\X)$ is equivalent to $\partial_{\varepsilon}:H_{p,q}(X;\F_2)\rightarrow H_{p+1,q-1}(X;\F_2)$.
	
	Let us now assume that $n=2k$ and write $\rk \partial^1_{k,k}=m_{k,k-1}(\varepsilon)+\delta$. Following Proposition~\ref{prop:rank_M}, we have $0\leq \delta\leq \rad(X)$. Moreover, by Poincaré duality of the spectral sequence, cf. \cite[Theorem 5.29]{chenal_poincare_2026}, we have $\rk \partial^1_{k-1,k+1}=\rk \partial^1_{k,k}=m_{k,k-1}(\varepsilon)+\delta=m_{k-1,k}(\varepsilon)+\delta$ and the second part of the statement follows.
	
	Only remains to show that the number $\delta$ equals the rank $r$ of the cokernel of the morphism $i_k:H_k(\X;\F_2)\rightarrow H_k(\R Y;\F_2)$. Note that we have following commutative diagram
	\begin{equation*}
		\begin{tikzcd}
			H_k(\X;\F_2)\ar[d] \ar[r] & H_k(\R Y;\F_2)\ar[d,"\cong" ] \\ 
			E^2_{k,k}(\X) \ar[r] & E^2_{k,k}(\R Y)
		\end{tikzcd}
	\end{equation*}
	So that $r$ equals the rank of cokernel of the induced map $\ker(\partial^1_{k,k})\rightarrow E^2_{k,k}(\R Y) =E^1_{k,k}(\R Y)$. Moreover, we have the commutative diagram with exact rows
	\begin{equation*}
		\begin{tikzcd}
			0 \ar[r] & \ker(i^1_{k,k}) \ar[d] \ar[r] & E^1_{k,k}(\X) \ar[d,"\partial^1_{k,k}"] \ar[r] & E^1_{k,k}(\R Y) \ar[d] \ar[r] & 0 \\
			0 \ar[r] & \ker(i^1_{k+1,k-1}) \ar[r] & E^1_{k+1,k-1}(\X) \ar[r] & 0 \ar[r] & 0 
		\end{tikzcd}
	\end{equation*}
	Hence, by the Snake Lemma we have the exact sequence
	\begin{equation*}
		\ker(\partial^1_{k,k})\rightarrow E^2_{k,k}(\R Y) \rightarrow \mathrm{coker}(\partial^1_{k,k}|_{\ker(i^1_{k,k})})  \rightarrow \mathrm{coker}(\partial^1_{k,k}) \rightarrow 0, 
	\end{equation*}
	and, by Proposition~\ref{prop:rank_M}, $r=h_{k+1,k-1}(X)-m_{k,k-1}(\varepsilon)-(h_{k+1,k-1}(X)-\rk\partial^1_{k,k})=\delta$ and the proposition follows.
\end{proof}

\section{Vanishings and Number of Connected Components}

\subsection{Twisted Edges}\label{sec:twists} Recall that, by Proposition~\ref{prop:exact_seq_wave}, for all $\sigma^{n+1}\in K$, we have an isomorphism $z:\N{2}\rightarrow Z^1(\sigma^{n+1};\F_2)$ sending $u$ to the cochain $z(u):\sigma^1\mapsto \la \omega(\sigma^1),u\ra$, and we see $P$ as a cell of $\R P$ via $\mu$.

\begin{prop}\label{prop:dual}
	Let $\sigma^{n+1}\in K$ and $u\in \N{2}$. The hypersurface $(u\cdot\X)\cap\sigma^{n+1}$ is dual to the cocycle $\df\varepsilon|_{\sigma^{n+1}}+z(u)$. 
\end{prop}

\begin{proof}
	By Definition~\ref{dfn:t-hyp}, $(u\cdot\X)\cap\sigma^{n+1}$ is the dual hypersurface of a cocycle $\zeta$. Let $\sigma^1\leq\sigma^{n+1}$, the dual cell of $\sigma^1$ belongs to $u\cdot\X$ if and only if the dual cell of $u\cdot\sigma^1$ belongs to $\X$. That is the case if and only if $\df\varepsilon(\sigma^1)+\la\omega(\sigma^1),u\ra=1$; cf. Definition~\ref{dfn:t-hyp}. Therefore, $\zeta=\df\varepsilon|_{\sigma^{n+1}}+z(u)$ and the proposition follows. 
\end{proof}

Let $e^1$ be a coarse bounded edge of $X$. It is dual to an interior simplex $\sigma^n\in K$. Hence, $e^1$ is contained in the maximal coarse cells of $Y$ dual to the vertices of $\sigma^n$. 

\begin{dfn}[Twists]
	Let $\varepsilon$ be a sign distribution, $e^1$ be a bounded edge of $X$ and $v$ be a vertex of the dual simplex of $e^1$. We say that $\X$ is \emph{twisted along $e^1$ in the direction of $v$} if $\D^2(\varepsilon+\chi_v)|_{e^1}^{\F_2}=1$. 
\end{dfn}

Let $e^1$ be a bounded edge of $X$ dual to the interior simplex $\sigma^n\in K$. The set of directions in which $\X$ can be twisted along $e^1$ is strongly influenced by the geometry of $K$. Indeed, using Proposition~\ref{prop:inter_rho} and Corollary~\ref{cor:wave_int}, we find that $\D^2(\varepsilon+\chi_v)|_{e^1}^{\F_2}=\D^2\varepsilon|_{e^1}^{\F_2}+\rho_{\sigma^n}(v)$ for all vertices $v\in \sigma^n$. Thus, $\X$ is twisted along $e^1$ in the direction of $v$ if and only if $\rho_{\sigma^n}(v)=1+\D^2\varepsilon|_{e^1}^{\F_2}$. The function $\rho_{\sigma^n}:\{\textnormal{Vertices of }\sigma^n\}\rightarrow \F_2$ partition the vertices of $\sigma^n$ into two sets $\rho_{\sigma^n}^{-1}(\{0\})$ and $\rho_{\sigma^n}^{-1}(\{1\})$. Thus, we see that the set of directions in which $\X$ is twisted along $e^1$ is necessarily one of these two sets. We will now give a geometrical condition for twistedness.

\begin{prop}\label{prop:separation}
	Let $S$ be the closed star of $\sigma^n$ in $K$. The T-hypersurface $\X$ is not twisted along $e^1$ in the direction of $v\in \sigma^n$ if and only if there exists $u\in \N{2}$ such that $v$ is isolated from the other vertices of $S$ in $S\setminus(u\cdot\X)$. 
\end{prop}

\begin{proof}
	Following Proposition~\ref{prop:kernel_hessian}, $\D^2(\varepsilon+\chi_v)|^{\F_2}_{e^1}=0$ if and only if there exists an affine sign distribution $a$ such that $(\varepsilon+a)|_S=\chi_v|_S$. Note that if $u$ is the slope of $a$ then $\R X_{\varepsilon+a}=u\cdot\X$. Hence, we have $S\setminus (u\cdot\X)=S\setminus \R X_{\chi_v}$ in which $v$ is necessarily isolated for $\chi_v(v)=1$ while $\chi_v(v')=0$ for any other vertex $v'\in S$.
	
	Conversely, since $v$ is connected to every other vertex of $S$, $v$ is isolated in $S\setminus \R X_\eta$ only if $\eta|_S=\chi_v|_S$ or $\eta|_S=(\chi_v+1)|_S$. Therefore, if $v$ is isolated in $S\setminus (u\cdot \X)$, then for any affine sign distribution $a$ with slope $u$, $(\varepsilon+a)|_S$ equals $\chi_v|_S$ up to a constant, i.e. $\D^2(\varepsilon+\chi_v)|_{e^1}^{\F_2}=0$ and $\X$ is not twisted along $e^1$ in the direction of $v$.
\end{proof}

Let us assume $\X$ is a curve. We consider an interior edge $\sigma^1$ of $K$ dual to a bounded edge $e^1$ of $X$. Following Lemma~\ref{lem:rho}, the function $\rho_{\sigma^1}$ is constant. Hence $\X$ can either be twisted along $e^1$ in all directions or in none. This corresponds to the usual notion of twist; cf. \cite[\S7]{renaudineau_bounding_2023} for instance.

\begin{ex}
	We represented two examples of curves in Figure~\ref{fig:twists_curves}. We depicted a bounded edge $e^1$ of $X$ as well as the closed star $S\leq K$ of its dual edge $\sigma^1$ in Figure~\ref{subfig:edge_and_closed_star}. In both Figures~\ref{subfig:untwisted_inter} and \ref{subfig:twisted_inter}, we gave the intersection $(u\cdot \X)\cap S$ for $u\in \N{2}=\F_2^2$ as well as the function $v\in S\mapsto \varepsilon(v) +\la v,u\ra$ to the differential of which $(u\cdot \X)\cap S$ is dual; cf. Proposition~\ref{prop:dual}. In Figure~\ref{subfig:untwisted_inter}, $\X$ is not twisted along $e^1$, while in Figure~\ref{subfig:twisted_inter} $\X$ twisted along $e^1$. See Proposition~\ref{prop:separation}.
\end{ex}

\begin{figure}[H]
	\centering
	\begin{subfigure}[t]{0.3\textwidth}
		\centering
		\input{untwisted_edge.tex}
		\caption{The four pieces $(u\cdot\X)\cap S$ with $e^1$ untwisted.}
		\label{subfig:untwisted_inter}
	\end{subfigure}
	\hfill
	\begin{subfigure}[t]{0.3\textwidth}
		\centering
		\input{the_edge.tex}
		\caption{The closed star $S$ of $\sigma$ with $X\cap S$ and the dual edge $e^1$.}
		\label{subfig:edge_and_closed_star}
	\end{subfigure}
	\hfill
	\begin{subfigure}[t]{0.3\textwidth}
		\centering
		\input{twisted_edge.tex}
		\caption{The four pieces $(u\cdot\X)\cap S$ with $e^1$ twisted.}
		\label{subfig:twisted_inter}
	\end{subfigure}
	\caption{The behaviour of a curve $\X$ in the four lifts of the closed star $S$ of an edge $\sigma^1\leq K$.}
	\label{fig:twists_curves}
\end{figure}  

	Let now $\X$ be a surface and $\sigma^2$ be an interior triangle of $K$ dual to a coarse bounded edge $e^1$ of $X$. As we discussed earlier, the set of directions in which $\X$ can be twisted along $e^1$ is determined by the geometry of $K$. If we assume that the function $\rho_{\sigma^2}$ vanishes, then $\X$ is twisted along $e^1$ in every direction or in none.

\begin{ex}
	Let $n=2$, $(\df x,\df y,\df z)$ be a basis of $\M{2}$ and $(\partial x, \partial y,\partial z)$ be the dual basis of $\N{2}$. Assume the bounded edge $e^1\in X$ is dual to the triangle $\sigma^2\in K$ given as the convex hull of $\{0,\df x,\df y\}$, and that the closed star $S$ of $\sigma^2$ is the convex hull of $\{0,\df x,\df y,\pm\df z\}$; see Figure~\ref{fig:S_in_surf}. We are in the configuration where $\rho_{\sigma^2}=0$. Therefore, a T-surface $\X$ is twisted along $e^1$ in all directions or in none.
	
	\begin{figure}[H]
		\centering
		\begin{subfigure}[t]{0.45\textwidth}
			\centering
			\input{star_dim_3.tex}
			\caption{The closed star $S$ of $\sigma^2$.}
		\end{subfigure}
		\hfill
		\centering
			\begin{subfigure}[t]{0.45\textwidth}
			\centering
			\input{XS_star_dim3.tex}
			\caption{The portion $X\cap S$ containing $e^1$.}
		\end{subfigure}
		\caption{The closed star of a triangle and the neighbourhood of its dual edge.}
		\label{fig:S_in_surf}
	\end{figure}
	
	In Figures~\ref{fig:surfaces_untwisted_nappes} and \ref{fig:surfaces_twisted_nappes}, we depicted the eight pieces of T-hypersurfaces $(u\cdot \X)\cap S$ as well as the function $v\in S\mapsto \varepsilon( v) +\la v,u\ra$ for all $u\in \N{2}$ and two different sign distributions. In Figure~\ref{fig:surfaces_untwisted_nappes}, $\varepsilon=0$ and $\X$ is twisted in no direction along $e^1$. Whereas in Figure~\ref{fig:surfaces_twisted_nappes}, $\varepsilon=\chi_{\df z}$ and $\X$ is twisted in every direction along $e^1$.
	
	\begin{figure}[H]
	\centering
	\begin{subfigure}[t]{0.22\textwidth}
	\centering
	\input{untwisted_edge_surf_000.tex}
	\caption{$u=(0,0,0)$}
\end{subfigure}
\begin{subfigure}[t]{0.22\textwidth}
	\centering
	\input{untwisted_edge_surf_100.tex}
	\caption{$u=(1,0,0)$}
\end{subfigure}
\begin{subfigure}[t]{0.22\textwidth}
	\centering
	\input{untwisted_edge_surf_110.tex}
	\caption{$u=(1,1,0)$}
\end{subfigure}
\begin{subfigure}[t]{0.22\textwidth}
	\centering
	\input{untwisted_edge_surf_010.tex}
	\caption{$u=(0,1,0)$}
\end{subfigure}
\begin{subfigure}[t]{0.22\textwidth}
	\centering
	\input{untwisted_edge_surf_001.tex}
	\caption{$u=(0,0,1)$}
\end{subfigure}
\begin{subfigure}[t]{0.22\textwidth}
	\centering
	\input{untwisted_edge_surf_101.tex}
	\caption{$u=(1,0,1)$}
\end{subfigure}
\begin{subfigure}[t]{0.22\textwidth}
	\centering
	\input{untwisted_edge_surf_111.tex}
	\caption{$u=(1,1,1)$}
\end{subfigure}
\begin{subfigure}[t]{0.22\textwidth}
	\centering
	\input{untwisted_edge_surf_011.tex}
	\caption{$u=(0,1,1)$}
\end{subfigure}
	\caption{The eight intersections $(u\cdot \X)\cap S$ for all $u\in \N{2}$ and $\varepsilon=0$.}
	\label{fig:surfaces_untwisted_nappes}
\end{figure}
	
	\begin{figure}[H]
	\centering
	\begin{subfigure}[t]{0.22\textwidth}
	\centering
	\input{twisted_edge_surf_000.tex}
	\caption{$u=(0,0,0)$}
\end{subfigure}
\begin{subfigure}[t]{0.22\textwidth}
	\centering
	\input{twisted_edge_surf_100.tex}
	\caption{$u=(1,0,0)$}
\end{subfigure}
\begin{subfigure}[t]{0.22\textwidth}
	\centering
	\input{twisted_edge_surf_110.tex}
	\caption{$u=(1,1,0)$}
\end{subfigure}
\begin{subfigure}[t]{0.22\textwidth}
	\centering
	\input{twisted_edge_surf_010.tex}
	\caption{$u=(0,1,0)$}
\end{subfigure}
\begin{subfigure}[t]{0.22\textwidth}
	\centering
	\input{twisted_edge_surf_001.tex}
	\caption{$u=(0,0,1)$}
\end{subfigure}
\begin{subfigure}[t]{0.22\textwidth}
	\centering
	\input{twisted_edge_surf_101.tex}
	\caption{$u=(1,0,1)$}
\end{subfigure}
\begin{subfigure}[t]{0.22\textwidth}
	\centering
	\input{twisted_edge_surf_111.tex}
	\caption{$u=(1,1,1)$}
\end{subfigure}
\begin{subfigure}[t]{0.22\textwidth}
	\centering
	\input{twisted_edge_surf_011.tex}
	\caption{$u=(0,1,1)$}
\end{subfigure}
	\caption{The eight intersections $(u\cdot \X)\cap S$ for all $u\in \N{2}$ and $\varepsilon=\chi_{\df z}$.}
	\label{fig:surfaces_twisted_nappes}
\end{figure}
\end{ex}

\subsection{Number of Connected Components}

\paragraph{Maximal Number of Connected Components} We recall from \cite[Theorem~1.4]{renaudineau_bounding_2023} that
\begin{equation}\label{eq:ineq_b0_RS}
	b_0(\X;\F_2)\leq 1+h^{n,0}(X).
\end{equation} When a T-hypersurface achieves this upper bound, we say that it has a maximal number of connected components. In this paragraph, we characterise the T-hypersurfaces achieving this upper bound. For instance, our computations allow us to reprove Haas' Theorem for curves.

\begin{thm}\label{thm:haas_classic} \textnormal{\cite[Theorem 7.3.0.10]{haas_real_1997} or \cite[Theorem 4.6]{bertrand_haas_2017}.}
	Let $n=1$. The T-curve $\X$ is maximal (i.e. has a maximal number of connected components) if and only if $\D^2\varepsilon|_{X_1}=\rho|_{X_1}$ and $\sum_{e^1\in\partial v^*} (\D^2\varepsilon+\rho)|^{\F_2}_{e^1}=0$ for all interior vertices $v\in K$.
\end{thm}

\begin{proof}
	From \cite[Theorem 6.1]{renaudineau_bounding_2023}, $\X$ is maximal if and only if $\partial_\varepsilon:H_{0,1}(X;\F_2)\rightarrow H_{1,0}(X;\F_2)$ vanishes, i.e. $m_{0,0}(\varepsilon)=0$. This is equivalent to the vanishing of $[f_{v_1}]\cdot \partial_\varepsilon [f_{v_2}]=0$ for all couples of interior vertices $v_1,v_2\in K$. By Proposition~\ref{prop:partial_eps_as_wave}, $[f_{v_1}]\cdot \partial_\varepsilon [f_{v_2}]=[f_{v_1}]\cdot \D^2(\varepsilon+\chi_{v_2})\capp [f_{v_2}]$. Thus, by Proposition~\ref{prop:wave_int_curves}, we have 
	\begin{equation*}
		[f_{v_1}]\cdot\partial_\varepsilon[f_{v_2}]=\left\{\begin{array}{ll} \D^2(\varepsilon+\chi_{v_2})|_{(\sigma^1)^*}^{\F_2} & \sigma^1=v_1*v_2 \\[5pt]
		\sum_{\sigma^1\geq v_1}\D^2(\varepsilon+\chi_{v_1})|_{(\sigma^1)^*}^{\F_2} & v_2=v_1 \\[5pt]
		0 & \textnormal{otherwise} \end{array}\right.
	\end{equation*}
	Since $n=1$, $K$ is $\rho$-uniform and we can write
	\begin{equation*}
		[f_{v_1}]\cdot\partial_\varepsilon[f_{v_2}]=\left\{\begin{array}{ll} (\D^2\varepsilon+\rho)|_{(\sigma^1)^*}^{\F_2} & \sigma^1=v_1*v_2 \\[5pt]
		\sum_{\sigma^1\geq v_1}(\D^2\varepsilon+\rho)|_{(\sigma^1)^*}^{\F_2} & v_2=v_1 \\[5pt]
		0 & \textnormal{otherwise} \end{array}\right.
	\end{equation*}
	Hence, $\D^2\varepsilon|_{X_1}=\rho|_{X_1}$ if and only if $[f_{v_1}]\cdot\partial_\varepsilon[f_{v_2}]=0$ for all interior vertices $v_1,v_2$ whose join belong to $K$. Likewise, for all interior vertex $v$, $\sum_{e^1\in\partial v^*} (\D^2\varepsilon+\rho)|^{\F_2}_{e^1}=0$ if and only if $[f_v]\cdot\partial_\varepsilon[f_v]=0$. Since $[f_{v_1}]\cdot\partial_\varepsilon[f_{v_2}]$ always vanishes in the other cases, the theorem follows.  
\end{proof}

\begin{rem}\label{rem:Dirichlet}
	In the special case of curves, one can compute $b_0(\X;\F_2)$ as the dimension of the space of solutions of a discrete Dirichlet problem. Let $\varepsilon\in C^0(K;\F_2)$ be a sign distribution. We consider the graph $\Gamma_\varepsilon$ obtained from the $1$-skeleton of $K$ by removing the edges $\sigma^1\in K$ with more than one boundary vertex or such that $(\D^2\varepsilon+\rho)|^{\F_2}_{(\sigma^1)^*}=0$. The boundary $\partial \Gamma_\varepsilon$ is the subgraph spanned by the vertices of $\Gamma_\varepsilon$ contained in $\partial K$. See Figure~\ref{fig:Gamma} for an example. Let us consider $\Delta:C^0(\Gamma_\varepsilon;\F_2)\rightarrow C^0(\Gamma_\varepsilon;\F_2)$ the Laplace operator 
	\begin{equation*}
		\Delta f:v\mapsto \sum_{\sigma^1\geq v} \df f(\sigma^1).
	\end{equation*}  
	We denote by $\mathcal{H}$ the space of solutions of the discrete Dirichlet problem
	\begin{equation*}
		\Delta f|_{\Gamma_\varepsilon\setminus\partial \Gamma_\varepsilon}=0 \quad\mathrm{and}\quad
			f|_{\partial \Gamma_\varepsilon}=0.
	\end{equation*}
	Then $\mathcal{H}$ has the same dimension as the kernel of $M_{0,0}(\varepsilon)$. Moreover, by Proposition~\ref{prop:dimE2}, we have $b_0(\X;\F_2)=1+\dim \mathcal{H}$. 
\end{rem}

\begin{figure}[h]
	\centering
	\begin{subfigure}[t]{.3\textwidth}
		\centering
		\input{Gamma1.tex}
		\caption{The triangulation $K$ and the sign distribution $\varepsilon$.}
	\end{subfigure}
	\hfill
	\begin{subfigure}[t]{.3\textwidth}
		\centering
		\input{Gamma2.tex}
		\caption{The twisted edges of $\X$.} 
	\end{subfigure}
	\hfill
	\begin{subfigure}[t]{.3\textwidth}
		\centering
		\input{Gamma3.tex}
		\caption{The graph $\Gamma_\varepsilon$ and its boundary $\partial \Gamma_\varepsilon$ (in white).} 
	\end{subfigure}
	\caption{The graph $\Gamma_\varepsilon$ of Remark~\ref{rem:Dirichlet}. Here, the curve $\X$ is connected.}
	\label{fig:Gamma}
\end{figure}

The novelty is for T-hypersurfaces of higher dimensions.

\begin{thm}\label{thm:gen_Haas}
	Let $n\geq 2$. The hypersurface $\X$ has a maximal number of connected components with respect to the Renaudineau--Shaw inequality if and only if \textnormal{(1)} $K$ is $\rho$-uniform; \textnormal{(2)} the morphism $i_1:H_1(\X;\F_2)\rightarrow H_1(\R Y;\F_2)$ is onto; and \textnormal{(3)} the sign distribution $\varepsilon$ satisfies the equations
	\begin{equation}\label{eq:equation_max_numb_comp}
		\D^2\varepsilon|_{X_1}=\rho|_{X_1}\quad \mathrm{and}\quad \D^2\varepsilon|_{X_0}\cupp\rho=\rho^2.
	\end{equation}
\end{thm}

\begin{proof}
	By the Renaudineau--Shaw spectral sequence, we have 
	\begin{equation*}
		b_0(\X;\F_2)=\dim E^\infty_{0,0}(\X) + \dim E^\infty_{n,0}(\X).
	\end{equation*}
	Moreover, $E^\infty_{0,0}(\X)=E^1_{0,0}(\X)\cong H_{0,0}(X;\F_2)$, thus, $b_0(\X;\F_2)=1 + \dim E^\infty_{n,0}(\X)$. Therefore, $\X$ has a maximal number of connected components if and only if $E^\infty_{n,0}(\X)=E^1_{n,0}(\X)$. Following \cite[Theorem 6.1]{renaudineau_bounding_2023} and \cite[Theorem 5.29]{chenal_poincare_2026}, this is equivalent to the vanishing of $\partial^1_{n-1,1}$ and $\partial^{n-1}_{1,1}$.
	
	If $n=2$, there is only one operator and the condition is equivalent to $E^2_{1,1}(\X)=E^1_{1,1}(\X)$. With Proposition~\ref{prop:dimE2}, this amounts to $m_{1,0}(\varepsilon) = \delta=0$ with $\delta$ the dimension of the cokernel of $i_1:H_1(\X;\F_2)\rightarrow H_1(\R Y;\F_2)$.
	
	If $n\geq 3$, Proposition~\ref{prop:rank_M} ensures that $\partial^1_{n-1,1}$ vanishes if and only if $m_{n-1,0}(\varepsilon)=0$. In addition, by \cite[Theorem 6.2]{chenal_poincare_2026}, $\partial^{n-1}_{1,1}=0$ if and only if $i_1$ is onto. 
	
	Therefore, for every $n\geq 2$, $\X$ has a maximal number of connected components if and only if $i_1$ is onto and $m_{n-1,0}(\varepsilon)=0$.
	
	Now we can conclude the proof by showing that $m_{n-1,0}(\varepsilon)=0$ if and only if $K$ is $\rho$-uniform and $\varepsilon$ is a solution of (\ref{eq:equation_max_numb_comp}). By Definition~\ref{dfn:inter_mat}, $m_{n-1,0}(\varepsilon)$ vanishes if and only if $[f_v]\cdot\partial_\varepsilon[f_{\sigma^{n-1}}]=0$ for all interior simplices $v,\sigma^{n-1}\in K$. Using (\ref{eq:auto_adj}), Propositions~\ref{prop:int_numb_wave_numb} and \ref{prop:partial_eps_as_wave}, and Lemma~\ref{lem:auto_adj_wave}, we have
	\begin{align*}
		[f_v]\cdot\partial_\varepsilon [f_{\sigma^{n-1}}]&=\partial_\varepsilon[f_v]\cdot [f_{\sigma^{n-1}}]=(\D^2(\varepsilon+\chi_v)\capp[f_v])\cdot [f_{\sigma^{n-1}}]= [f_v]\cdot \D^2(\varepsilon+\chi_v)\capp[f_{\sigma^{n-1}}]\\&=(\D^2\chi_v\cupp\D^2(\varepsilon+\chi_v))|_{e^2}^{\F_2},
	\end{align*}
	where $e^2=(\sigma^{n-1})^*$. So $m_{n-1,0}(\varepsilon)$ vanishes if and only if $(\D^2\chi_v\cupp\D^2(\varepsilon+\chi_v))|_{e^2}^{\F_2}=0$ for all interior vertices $v\in K$ and all coarse bounded cells $e^2\in X$. 
	
	Let us assume $m_{n-1,0}(\varepsilon)=0$, and let $e^1$ be an edge of $\partial v^*$. It is dual to an interior simplex $\sigma^n\in K$. We let $\sigma^{n-1}=\sigma^n-v$ and $e^2$ be the dual cell of $\sigma^{n-1}$. If $e^2$ is bounded, by the second part of Proposition~\ref{prop:int_numb_wave_numb}, we have 
	\begin{equation*}
		(\D^2\chi_v\cupp\D^2(\varepsilon+\chi_v))|_{e^2}^{\F_2}=\D^2(\varepsilon+\chi_v)|^{\F_2}_{e^1}=\D^2\varepsilon|^{\F_2}_{e^1}+\rho_{\sigma^n}(v).
	\end{equation*}
	Let us assume that $\sigma^n$ has more than one interior vertices which we denote by $v_1,...,v_k$, with $k\geq 2$. Thus, $\sigma^{n-1}_i=\sigma^n-v_i$ is an interior simplex of $K$ and its dual $2$-cell is bounded. Thus, if $m_{n-1,0}(\varepsilon)=0$, we have
	\begin{equation*}
		\rho_{\sigma^n}(v_i)=\D^2\varepsilon|^{\F_2}_{e^1},\quad \textnormal{for all }1\leq i\leq k.
	\end{equation*}
	Therefore, by Proposition~\ref{prop:crit_rho_unif}, $K$ is $\rho$-uniform. Moreover, if $e^1\in X_1$, then, by Definition~\ref{dfn:X0}, its dual simplex $\sigma^n$ contains an interior vertex $v$ whose opposite face $\sigma^{n-1}=\sigma^n-v$ is also interior. In this case,
	\begin{equation*}
		(\D^2\varepsilon|_{X_0}+\rho)|^{\F_2}_{e^1}=(\D^2\varepsilon|_{X_0}+\rho)|_{\partial v^*}|^{\F_2}_{e^1}=\D^2(\varepsilon+\chi_v)|^{\F_2}_{e^1}=(\D^2\chi_v\cupp\D^2(\varepsilon+\chi_v))|_{(\sigma^{n-1})^*}^{\F_2}=0.
	\end{equation*}
	Thus, in addition of ensuring the $\rho$-uniformity of $K$, $m_{n-1,0}(\varepsilon)=0$ also implies $\D^2\varepsilon|_{X_1}=\rho|_{X_1}$. Let now $e^2\in X_0$, then there is an interior vertex $v\in K$ for which $e^2\in \partial v^*$. Therefore,
	\begin{equation*}
		(\rho\cupp(\D^2\varepsilon+\rho))|^{\F_2}_{e^2}=(\rho\cupp(\D^2\varepsilon+\rho))|_{\partial v^*}|^{\F_2}_{e^2}= (\D^2\chi_v\cupp\D^2(\varepsilon+\chi_v))|_{e^2}^{\F_2}=0,
	\end{equation*}
	and the last equation is proved.
	
	Conversely, let us assume that $K$ is $\rho$-uniform and that $\varepsilon$ is a solution of (\ref{eq:equation_max_numb_comp}). To establish $m_{n-1,0}(\varepsilon)=0$, we need to prove that $(\D^2\chi_v\cupp\D^2(\varepsilon+\chi_v))|_{e^2}^{\F_2}=0$ for every interior vertex $v\in K$ and every coarse bounded cell $e^2\in X$. With Proposition~\ref{prop:int_numb_wave_numb}, we only need to prove this when there is a simplex $\sigma^{n+1}\in K$ containing both $v$ and $\sigma^{n-1}$, the dual simplex of $e^2$. There are two cases, either $\sigma^n=v*\sigma^{n-1}$ is a simplex of $K$ or $v\in \sigma^{n-1}$. In the first case, $e^1=(\sigma^n)^*$ belongs to $X_1$ and we have
	\begin{equation*}
		(\D^2\chi_v\cupp\D^2(\varepsilon+\chi_v))|_{e^2}^{\F_2}=\D^2(\varepsilon+\chi_v)|^{\F_2}_{e^1}=\D^2(\varepsilon+\rho)|^{\F_2}_{e^1}=0,
	\end{equation*}
	by the first equation of (\ref{eq:equation_max_numb_comp}). In the second case, $e^2$ belongs to $\partial v^*$ and then
	\begin{equation*}
		(\D^2\chi_v\cupp\D^2(\varepsilon+\chi_v))|_{e^2}^{\F_2}=(\rho\cupp(\D^2\varepsilon+\rho))|_{\partial v^*}|^{\F_2}_{e^2}= (\rho\cupp(\D^2\varepsilon+\rho))|^{\F_2}_{e^2} =0,
	\end{equation*}
	by the second equation of (\ref{eq:equation_max_numb_comp}). 
\end{proof}

\begin{rem}\label{rem:haas_twists}
	As announced in the introduction, both Theorem~\ref{thm:haas_classic} and Theorem~\ref{thm:gen_Haas} can be rephrased in terms of twists. Since for curves a bounded edge is either twisted in all directions or in none, we usually say that it is twisted or untwisted without mention of the direction. The relevant generalisation of this terminology for T-hypersurfaces over $\rho$-uniform triangulations would be to say that $\X$ is \emph{regurlarly untwisted} along an edge $e^1\in X_0$ if $(\D^2\varepsilon+\rho)|_{e^1}^{\F_2}=0$. Likewise, we say that $e^1$ is a \emph{regularly twisted} edge if $(\D^2\varepsilon+\rho)|_{e^1}^{\F_2}=1$. Therefore, $e^1$ is regularly untwisted if and only if $\X$ is untwisted along $e^1$ in the direction of any interior vertex contained in the dual simplex.
	
	 Then, if $\X$ has a maximal number of connected components, the equation $\D^2\varepsilon|_{X_1}=\rho|_{X_1}$ implies that the edges of $X_1$ are all regularly untwisted. 
	 
	 Moreover, in the case of curves, the equations $\sum_{e^1\in \partial v^*}(\D^2\varepsilon+\rho)|_{e^1}^{\F_2}$ for all $v\in K\setminus \partial K$ impose that the number of regular twists contained in a cycle of $X$ is always even. In higher dimensions, the equation $\D^2\varepsilon|_{X_0}\cupp\rho=\rho^2$ can be interpreted similarly. Indeed, if $e^2$ is a cell of $X_0$ dual to a simplex $\sigma^{n-1}$ containing an interior vertex $v$, we have
	 \begin{equation*}
	 	(\rho\cupp(\D^2\varepsilon+\rho))|^{\F_2}_{e^2}=\sum_{\sigma^n\geq  \sigma^{n-1}} \la \omega(\sigma^n-v),\omega^*(\sigma^{n-1})\ra (\D^2\varepsilon+\rho))|_{(\sigma^n)^*}^{\F_2}.
	 \end{equation*}
	 The number $(\rho\cupp(\D^2\varepsilon+\rho))|^{\F_2}_{e^2}$ can then be read as a weighted number of regular twists along the boundary of $e^2$. The equation $\D^2\varepsilon|_{X_0}\cupp\rho=\rho^2$ then asks for these weighted numbers of regular twists to be even. This expansion can however be somewhat unsatisfactory for the weights depend on both the chosen vertex $v$ and the choice of $\omega^*$. This is why we prefer the canonical form $\D^2\varepsilon|_{X_0}\cupp\rho=\rho^2$ of these equations. 
\end{rem}

\begin{prop}\label{prop:simple_harnack}
	Let $K$ be a $\rho$-uniform triangulation. If $\varepsilon$ is a solution of
	\begin{equation}\label{eq:simple_harnack}
		\D^2\varepsilon|_{X_0}=\rho,
	\end{equation}
	then $\X$ is made of \textnormal{(1)} $h^{n,0}(X)$ spheres contained in the torus of $\mathcal{Y}(\R)$ where they bound disjoint balls; and \textnormal{(2)} an additional connected component that intersects every toric divisor of $\mathcal{Y}(\R)$. When (\ref{eq:simple_harnack}) has a solution, we say that $K$ is \emph{simply integrable}.
\end{prop}

\begin{proof}
	Let $\varepsilon$ be a solution of (\ref{eq:simple_harnack}), $v$ be an interior vertex of $K$ and $S$ be the closed star of $v$. We have $X_S=\partial v^*$. Following Proposition~\ref{prop:kernel_hessian}, $\D^2\varepsilon|_{X_0}=\rho$ implies that $\D^2(\varepsilon|_{S}+\chi_v|_{S})=0$ and there is an affine sign distribution $a$ such that $\varepsilon|_S=\chi_v|_{S}+a|_{S}$. Thus, if $u\in \N{2}$ is the slope of $a$, we have $$(u\cdot S)\cap \X=u\cdot(S\cap (u\cdot\X))=u\cdot(S\cap \R X_{\chi_v}).$$
	The space $S\cap \R X_{\chi_v}$ consists of an $n$-sphere contained in the interior of $S$ and isolating $v$ from the other vertices of $S$. As a consequence, $\X$ is made of at least $h^{n,0}(X)$ spheres contained in the torus of $ \mathcal{Y}(\R)$ and all bounding disjoint balls. Moreover, note that if $Q$ is a face of $P$, then $\X\cap \R Q$ is the T-hypersurface builded on the triangulation $K\cap Q$ with sign distribution $\varepsilon|_{Q}$. Hence, $\X$ has at least $h^{n,0}(X)+1$ connected components one of which intersecting every toric divisor and the proposition follows. 
\end{proof}

\begin{prop}\label{prop:solvability_simple_harnack}
	Let $n\geq 2$. A $\rho$-uniform triangulation $K$ of a non-singular polytope $P$ is simply integrable if and only if $\df \rho=0$ where $\df$ denotes the connecting morphism
	\begin{equation*}
		\df: H^1(X_0;\cW_1) \rightarrow H^2(X;X_0;\cW_1),
	\end{equation*}
	of the long exact sequence associated to the couple $X_0\subset X$. This is especially the case when the subgraph of the $1$-skeleton of $K$ spanned by the interior vertices is connected. 
\end{prop}

\begin{proof}
	By the property of the cohomology long exact sequence, $\df \rho=0$ if and only if there is a wave $\tilde{\rho}\in H^1(X;\cW_1)$ such that $\rho=\tilde{\rho}|_{X_0}$. Since $\D^2:C^0(K;\F_2)\rightarrow H^1(X;\cW_1)$ is onto, cf. Corollary~\ref{cor:exact_seq_hess}, $\tilde{\rho}$ can always be taken as the Hessian of a sign distribution $\varepsilon$. Hence, $\df\rho=0$ if and only if $K$ is simply integrable. 
	
	For the second part, with the additional assumption on $K$, we will show that $H^2(X;X_0;\cW_1)$ vanishes. In this case, $X_0$ is connected and the restrictions $H^k(X;\F_2)\rightarrow H^k(X_0;\F_2)$ are isomorphisms for every integer $k\ge 0$ (when $X_0$ is disconnected this is only true for $k>0$). Since we have the exact sequence
	\begin{equation*}
		H^1(X;X_0;\N{2}) \rightarrow H^1(X;X_0;\cZ^1) \rightarrow H^2(X;X_0;\cW_1) \rightarrow H^2(X;X_0;\N{2}),
	\end{equation*}
	see Proposition~\ref{prop:exact_seq_wave}, the group $H^2(X;X_0;\cW_1)$ is isomorphic to $H^1(X;X_0;\cZ^1)$. Furthermore, we also have the following exact sequence
	\begin{equation*}
		H^1(X;X_0;\F_2) \rightarrow H^1(X;X_0;\mathcal{C}^0) \rightarrow H^1(X;X_0;\cZ^1) \rightarrow H^2(X;X_0;\F_2),
	\end{equation*} 
	see Equation~(\ref{eq:seq_C0}). Therefore, $H^2(X;X_0;\cW_1)$ is even isomorphic to $H^1(X;X_0;\mathcal{C}^0)$. In addition, since $\mathcal{C}^0$ is flabby, $H^1(X;X_0;\mathcal{C}^0)$ is the cokernel of the restriction $H^0(X;\mathcal{C}^0)\rightarrow H^0(X_0;\mathcal{C}^0)$. We can directly compute $H^0(X_0;\mathcal{C}^0)$. By definition, we have
	\begin{equation*}
		0\rightarrow H^0(X_0;\mathcal{C}^0) \rightarrow \prod_{\substack{\sigma^{n+1}\in K \\ \textnormal{with an int. vert.}}} C^0(\sigma^{n+1};\F_2) \rightarrow \prod_{\substack{\sigma^{n}\in K \\ \textnormal{with an int. vert.}}} C^0(\sigma^{n};\F_2).
	\end{equation*}
	Therefore, if $S_v$ denotes the closed star of $v$ in $K$, and $G$ is the subgraph of the $1$-skeleton of $K$ spanned by the interior vertices, we have the following commutative diagram 
	\begin{equation*}
		\begin{tikzcd}
			H^0(X;\mathcal{C}^0) \ar[d,"\textnormal{restriction}", swap] \ar[r,"\cong"]& C^0(K;\F_2)\ar[d,"\textnormal{prod.~of restrictions}"] \\
			H^0(X_0;\mathcal{C}^0) \ar[r,"\cong", swap] & \displaystyle \prod_{C\in \pi_0(G)} C^0({\textstyle\bigcup_{v\in C}} S_v;\F_2),
		\end{tikzcd}
	\end{equation*}
	and the restriction $H^0(X;\mathcal{C}^0)\rightarrow H^0(X_0;\mathcal{C}^0)$ is onto when $G$ is connected. This implies that $H^2(X;X_0;\cW_1)$ vanishes and the proposition follows. 
\end{proof}

Proposition~\ref{prop:rho_int_curves} implies that (\ref{eq:simple_harnack}) always has a solution when $n=1$, more precisely we have the following theorem:

\begin{thm}\label{thm:quadratic_distrib}
	Let $n=1$ and $\tilde{\rho}\in H^1(X;\cW_1)$ be the integrable wave of Proposition~\ref{prop:rho_int_curves}. The solutions of the equation
	\begin{equation}\label{eq:simple_harnack_2}
		\D^2\varepsilon=\tilde{\rho},
	\end{equation}
	are the restrictions of quadratic polynomial functions $\M{2}\rightarrow \F_2$ to the vertices of $K$. As such, every bidimensional triangulation is simply integrable. 
\end{thm}

\begin{proof}
	Let $h$ be a quadratic sign distribution. Now let $v_0$ and $v_1$ be the two vertices of an interior edge of $K$, and $u_0,u_1$ be the apices of the two triangles of $K$ adjacent to $v_0*v_1$. Following Proposition~\ref{prop:kernel_hessian}, affine maps are in the kernel of $\D^2$. Since $h$ is quadratic, we can assume, to compute $\D^2h|^{\F_2}_{v_0*v_1}$, that $h$ is given by $(x;y)\mapsto xy$ where $(x,y)$ are the affine coordinates in the basis $(v_1+v_0,u_0+v_0)$ with origin at $v_0$. We have $y(u_1)=1$ since $(v_1+v_0,u_1+v_0)$ is also a basis of $\N{2}$. Hence, we have:
	\begin{equation*}
			u_0+u_1=x(u_1)(v_1+v_0),
	\end{equation*}
	and, by Lemma~\ref{lem:rho}, $\rho_{v_0*v_1}(v_i)$ equals $x(u_1)$ for any $i\in\{0;1\}$, i.e. $\tilde{\rho}|^{\F_2}_{v_0*v_1}=x(u_1)$. Then, following Proposition~\ref{prop:formula_hess}, we have
	\begin{equation*}
		\begin{split}
			\D^2h|^{\F_2}_{v_0*v_1} &=x(u_0)y(u_0)+x(u_1)y(u_1)+x(u_1) \big(x(v_0)y(v_0)+x(v_1)y(v_1)\big) \\
			&= 0\times 1 +x(u_1)\times 1+ x(u_1)(0\times 0 + 1\times 0)\\
			&= x(u_1)=\tilde{\rho}|^{\F_2}_{v_0*v_1}.
		\end{split}
	\end{equation*}
	And $\D^2h$ equals $\tilde{\rho}$. Now by Proposition~\ref{prop:kernel_hessian}, the kernel of $\D^2$ is exactly the space of affine sign distributions, hence, the theorem follows.
\end{proof}

\begin{rem} Since $\M{2}$ has dimension $2$, there is essentially one quadratic polynomial function $\M{2}\rightarrow \F_2$. It is a consequence of Fermat's little theorem. Hence, all the T-curves $\X$ obtained from a solution $\varepsilon$ of (\ref{eq:simple_harnack_2}) are reflections of one another.

Itenberg and Viro proved in \cite[\S Patchworking Harnack Curves]{itenberg_patchworking_1996} that a particular sign distribution always produces a maximal curve of a given topological type for all convex triangulations of the standard triangle of size $d\geq 1$: ${\{(x,y)\in\R^2_+\;|\; x+y\leq d\}}$. They named their sign distribution the \emph{Harnack distribution}. It is given by the formula $\varepsilon(x,y)\coloneqq 1+(1+x)(1+y)\,(\mod 2)$, hence, by a quadratic polynomial.
\end{rem}

We conclude by providing an upper bound on the number of connected components of a T-hypersurface constructed on a non-$\rho$-uniform triangulation.

\begin{dfn}\label{dfn:kappa}
	Let $K$ be a primitive triangulation of an $(n+1)$-dimensional polytope. The distance between two simplices is the minimum of the graph length of their vertices in the 1-skeleton of $K$. We denote by $\kappa(K)$ the maximal number of $n$-simplices of $K$ at which it fails to be $\rho$-uniform and that are all at distance at least $2$ of each others.
\end{dfn}

\begin{prop}\label{prop:sharper_upper_bound}
	For all $\varepsilon$, we have  $b_0(\X;\F_2)\leq h^{n,0}(X)+1 - \kappa(K)$.
\end{prop}

\begin{proof}
	Let $\sigma^{n}_1$,...,$\sigma^{n}_k$ be simplices of $K$ at which it fails to be $\rho$-uniform and that are all at distance at least $2$ of each others. Let $\varepsilon$ be a sign distribution on $K$. Since $K$ fails to be $\rho$-uniform at the simplices $\sigma^{n}_1$,...,$\sigma^{n}_k$, we can find, for all $1\leq i\leq k$, two different interior vertices $v_i,v'_i$ of $\sigma^{n}_i$ for which $\rho_{\sigma^{n}_i}(v_i)$ and $\rho_{\sigma^{n}_i}(v_i)$ are different. Following Proposition~\ref{prop:int_numb_wave_numb}, either $[f_{v_i}]\cdot \partial_\varepsilon [f_{\sigma^n_i-v_i}]$ or $[f_{v'_i}]\cdot \partial_\varepsilon [f_{\sigma^n_i-v'_i}]$ equals 1. We assume that it happens on $v_i$, for all $1\leq i\leq k$. The submatrix
	\begin{equation*}
		\Big([f_{v_i}]\cdot \partial_\varepsilon [f_{\sigma^n_j-v_j}]\Big)_{1\leq i,j\leq k}
	\end{equation*}
	of $M_{n-1,0}(\varepsilon)$ is the identity matrix of size $k$. Indeed, for all $i$ different from $j$, $v_i$ and $\sigma^{n}_j-v_j$ are at distance at least 2. Therefore, they cannot be both included in the same simplex, so $[f_{v_i}]\cdot \partial_\varepsilon [f_{\sigma^n_j-v_j}]=0$ by Proposition~\ref{prop:int_numb_wave_numb}. Hence, the rank of $M_{n-1,0}(\varepsilon)$ is at least equal to $k$. Consequently, the rank of $M_{n-1,0}(\varepsilon)$ is at least equal to $\kappa(K)$, and the upper bound follows from Proposition~\ref{prop:dimE2}.
\end{proof}

\subsection{Growth Rate of the Expected Number of Connected Components} \label{sec:growth}

Let $P$ be a smooth $(n+1)$-dimensional lattice polytope of $\M{\R}$. We consider the sequence $(\X(d))_{d\geq 1}$ of random T-hypersurfaces of “degrees” $(dP)_{d\geq 1}$ constructed as follows: for all $d\geq 1$,\begin{enumerate}
	\item We choose a random primitive triangulation $K(d)$ of $dP$ following any distribution.
	\item We choose independently and uniformly at random a sign distribution $\varepsilon:dP\cap M\rightarrow \F_2$.
	\item We form the T-hypersurface $\X(d)$ associated to $(K(d),\varepsilon)$.
\end{enumerate}
We denote by $X(d)$ the random cell complex associated to $K(d)$. The aim of this section is to show that the growth rate of the sequence of expectations $(\mathbb{E}[b_0(\X(d)])_{d\geq 1}$ is of the order of $d^{n+1}$. Let $L(P;x)\in\Z[x]$ denote the Ehrhart polynomial of $P$. It is the unique polynomial such that $L(P;d)=\card (dP\cap M)$ for all non-negative integers $d$; see \cite[Théorème 1]{ehrhart_sur_1962}.

\begin{prop}\label{prop:growth_h_n0}
	For all integers $d\geq 1$, we have
	\begin{equation*}
		h^{n,0}\big(X(d)\big)=\card(\mathrm{rel.int.}(dP)\cap M) =\vol_\Z(P)d^{n+1}+o(d^{n+1}),
	\end{equation*}
	where $\vol_\Z(P)$ is the integral volume of $P$. That is to say the Lebesgue measure of $P$, normalised so that a parallelogram on a basis of $M$ has volume $1$.
\end{prop}

\begin{proof}
	The first equality is a widely known fact for Hodge numbers of smooth hypersurfaces of toric varieties; cf. \cite[\S 5.5]{danilov_newton_1987}. The equality for tropical Hodge numbers follows \cite[Corollary 1.9]{arnal_lefschetz_2021} and \cite[Theorem 1.4]{brugalle_combinatorial_2024}. By Ehrhart-Macdonald reciprocity, cf. \cite[Theorem 4.6]{macdonald_polynomials_1971}, it equals $(-1)^{n+1}L(P;-d)$. Moreover, the polynomial $L(P;x)$ has degree $n+1$ and its leading coefficient is $\vol_\Z(P)$ so the proposition follows.
\end{proof}

\begin{prop}\label{prop:asymp_f1}
	For all integers $0\leq q \leq n+1$ and $d\geq 1$, the random variable $f_q(K(d))\coloneqq \dim C^q(K(d);\F_2)$ is deterministic. Furthermore, the function ${d\mapsto f_1(K(d))}$ is a polynomial of degree $n+1$. 
\end{prop}

\begin{proof}
	Let $d\geq 1$ be an integer and $K$ be a primitive triangulation of $dP$. The polytope $dP$ is partinonned into the relative interiors of its simplices. Therefore, by Ehrhart-Macdonald reciprocity, we have
	\begin{equation*}
		L(dP,k)= \sum_{q\geq 0} \sum_{\sigma^q\in K} (-1)^q L(\sigma^q,-k)= \sum_{q\geq 0} \sum_{\sigma^q\in K} \frac{1}{q!}\prod_{i=1}^q (k-i) = \sum_{g\geq 0} \frac{f_q(K)}{q!} \prod_{i=1}^q (k-i),
	\end{equation*}
	with $L(\sigma^q,k)={\scriptstyle \binom{q+k}{q}}$. Since the polynomials $(\tfrac{1}{q!}\prod_{i=1}^q(k-i))_{q\geq 0}$ form a basis of $\Q [k]$, the numbers $(f_q(K))_{q\geq 0}$ do not depend on the particular choice of the primitive triangulation $K$ of $dP$ and the $f_q(K(d))$'s are deterministic random variables. Furthermore, by \cite[Théorème 1]{ehrhart_sur_1962}, we have
	 \begin{equation*}
		\sum_{k\geq 0}L(dP,k)t^k= \frac{\sum_{q=0}^{n+1} h^*_q(dP) t^q}{(1-t)^{n+2}}.
	\end{equation*}
	Thus, $f_1(K(d))$ is a linear combination of the $h^*_q(dP)$'s. In particular, there is an $M$-translation invariant valuation $\varphi$ on the integral polytopes of $\M{\R}$ such that $f_1(K(d))$ equals $\varphi(dP)$ for all $d\geq 1$. By \cite[Theorem~5]{mcmullen_valuations_1977}, $d\mapsto f_1(K(d))$ is given by a polynomial of degree at most $n+1$. However, the number  $f_1(K(d))$ is not smaller than $\frac{1}{2}f_0(K(d))$ since no vertex is isolated. Hence, $d\mapsto f_1(K(d))$ is of degree $n+1$ for ${d\mapsto f_0(K(d))=\card( dP\cap M )=L(P,d)}$ is polynomial of degree $n+1$. 
\end{proof}

\begin{thm}\label{thm:growth_rate}
	In Landau's notation, $\mathbb{E}[b_0(\X(d))]=\Theta(d^{n+1})$, i.e. it is eventually bounded above and below by positive multiples of $d^{n+1}$.
\end{thm}

\begin{proof}
	By the Renaudineau--Shaw inequality~(\ref{eq:ineq_b0_RS}) and Proposition~\ref{prop:growth_h_n0}, we know $\mathbb{E}[b_0(\X(d))]$ is bounded from above by a multiple of $d^{n+1}$. Only remains to prove the lower bound. To do so, we will show that
	\begin{equation*}
        \liminf_{d\rightarrow+\infty}\mathbb{E}\left[\frac{b_0(\X(d))}{1+h^{n,0}(X(d))}\right] > 0.
    \end{equation*}
    This will be enough thanks to Proposition~\ref{prop:growth_h_n0}. Let $d\geq 1$ be an integer and $v\in \relint(dP)\cap M$ be an interior lattice point. We denote by $S_v(K(d))$ the closed star of $v$ in $K(d)$, and $\xi_v$ the Bernoulli random variable whose value is $1$ if and only if $\X(d)$ contains a lift of the sphere $\partial v^*$. By Propositions~\ref{prop:separation}~and~\ref{prop:simple_harnack}, $\xi_v(K(d),\varepsilon)$ equals 1 if and only if $\D^2\varepsilon |_{\partial v^*}=\D^2\chi_v|_{\partial v^*}$ in $H^1(\partial v^*;\cW_1)$. We have the following commutative diagram with exact rows and columns:
	\begin{equation*}
		\begin{tikzcd}[row sep=small]
			& & 0 \ar[d] & \\
			& 0 \ar[d] & C^0(K(d);S_v(K(d));\F_2) \ar[d] & \\
			0 \ar[r] & \Aff(\N{2}) \ar[r] \ar[d] & C^0(K(d);\F_2) \ar[r,"\D^2"] \ar[d] & H^1(X(d);\cW_1) \ar[d] \\
			0 \ar[r] & \Aff(\N{2}) \ar[r] \ar[d] & C^0(S_v(K(d));\F_2) \ar[r,"\D^2"] \ar[d] & H^1(\partial v^*;\cW_1) \\
			 & 0 & 0 \\
		\end{tikzcd}
	\end{equation*}
	Hence, we have the following exact sequence
	\begin{equation*}
			0 \rightarrow C^0(K(d);S_v(K(d));\F_2)\oplus \Aff(\N{2}) \rightarrow C^0(K(d);\F_2) \xrightarrow{|_{\partial v^*}\circ \D^2} H^1(\partial v^*;\cW_1).
	\end{equation*}
	Thus, the equation $\D^2\varepsilon|_{\partial v^*}=\D^2\chi_v|_{\partial v^*}$ has exactly $2^{n+2+f_0(K(d))-f_0(S_v(K(d)))}$ solutions. Therefore, the conditional expectation of $\xi_v$ knowing $K(d)$ satisfies
	\begin{equation*}
		\mathbb{E}[\xi_v|K(d)]=2^{n+1-\deg_{K(d)}(v)},
	\end{equation*}
	where $\deg_{K(d)}(v)$ is the degree of the vertex $v$ in $K(d)$, i.e. $\deg_{K(d)}(v)=f_0(S_v(K(d)))-1$. Moreover, we have the following inequality
	\begin{equation*}
		b_0(\X(d))\geq \sum_{v\in \relint(dP)\cap M} \xi_v,
	\end{equation*}
	which implies
	\begin{equation*}
		\mathbb{E}\left[\left.\frac{b_0(\X(d))}{1+h^{n,0}(X(d))}\right|K(d)\right] \geq \sum_{v\in \relint (dP)\cap M}\frac{ 2^{n+1-\deg_{K(d)}(v)}}{1+h^{n,0}\big(X(d)\big)}.
	\end{equation*}
	By convexity of the function $x\mapsto 2^{-x}$, we have
	\begin{equation*}
		\sum_{v\in \relint (dP)\cap M} 2^{-\deg_{K(d)}(v)}\geq  h^{n,0}(X(d)) 2^{-\frac{\sum_{v}\deg_{K(d)}(v)}{h^{n,0}(X(d))}} \geq h^{n,0}(X(d)) 2^{-2f_1(K(d))/h^{n,0}(X(d))},
	\end{equation*}
	for $\card( \relint (dP)\cap M)=h^{n,0}(X(d))$ and
	\begin{equation*}
		\sum_{v\in \relint (dP)\cap M} \deg_{K(d)}(v)\leq \sum_{v\in dP\cap M}\deg_{K(d)}(v)=2f_1(K(d)). 
	\end{equation*}
	Therefore, we have the deterministic lower bound
	\begin{equation*}
		\mathbb{E}\left[\left.\frac{b_0(\X(d))}{1+h^{n,0}(X(d))}\right|K(d)\right] \geq \frac{h^{n,0}(X(d))}{1+h^{n,0}(X(d))} 2^{n+1-2f_1(K(d))/h^{n,0}(X(d))}.
	\end{equation*}
	By Propositions~\ref{prop:growth_h_n0}~and~\ref{prop:asymp_f1}, the sequence $(f_1(K(d))/h^{n,0}(X(d)))_{d\geq 1}$ converges toward a number $c>0$ as $d$ grows to $+\infty$. Hence,
	\begin{align*}
		\liminf_{d\rightarrow+\infty}\mathbb{E}\left[\frac{b_0(\X(d))}{1+h^{n,0}(X(d))}\right] &= \liminf_{d\rightarrow+\infty}\mathbb{E}\Bigg[\mathbb{E}\left[\left.\frac{b_0(\X(d))}{1+h^{n,0}(X(d))}\right|K(d)\right]\Bigg] \geq 2^{n+1-2c} >0,
	\end{align*}
	and the theorem follows.
\end{proof}

\section{Examples and Non-Examples}\label{sec:examples}

    Let $n$ and $d$ be two positive integers. We denote the standard $(n+1)$-simplex of size $d$ by the symbol $\P^{n+1}_d$. In this context, $\M{}=\Z^{n+1}$, we denote by $(e_1,...,e_{n+1})$ is the canonical basis of $\M{}$, and $\P^{n+1}_d$ is the convex hull of $\{0,de_1,...,de_{n+1}\}$.
    
\paragraph{Regular Triangulations of Knudsen} In the third chapter of \cite{kempf_toroidal_1973}, F. Knudsen introduces a family of primitive convex triangulations of $(\P^n_d)_{n,d\geq 1}$ called \emph{regular triangulations}. For all $n,d\geq 1$, the regular triangulation $R^n_d$ of $\P^n_d$ is obtained by cutting $\P^n_d$ by the following  hyperplanes:
\begin{equation*}
	H^{i,j}_k\coloneqq \left\{ \sum_{l=i}^j x_l = k\right\},
\end{equation*}
for all $1\leq i\leq j\leq n$, and all $1\leq k \leq d-1$. \cite[Lemma~2.12]{kempf_toroidal_1973} implies that $R^n_d$ is never $\rho$-uniform when $n$ is bigger than $2$ and $d$ is big enough.

\paragraph{Viro Triangulations} In \cite[Definition~7.3]{chenal_poincare_2026}, we described a method of construction of triangulations of $\P^{n+1}_{d+1}$ from a triangulation of $\P^{n+1}_{d}$ and of $\P^{n}_{d+1}$: If $K$ is a triangulation of $\P^{n+1}_d$ and $L$ is a triangulation of $\P^{n}_{d+1}$, $K+L$ is the unique triangulation of $\P^{n+1}_{d+1}$ for which:
\begin{enumerate}
	\item The restriction of $K+L$ to $\P^{n+1}_{d+1}\cap\{x_{n+1}\geq 1\}$ is the translation of $K$;
	\item For all integers $0\leq i\leq n$, the restriction of the triangulation $K+L$ to the join of simplices ${[e_{n+1},e_{n+1}+de_1,...,e_{n+1}+de_i]*[(d+1)e_i,...,(d+1)e_{n}]}$, is the join of the restriction of the translation of $K$ to ${[e_{n+1},e_{n+1}+de_1,...,e_{n+1}+de_i]}$ with the restriction of $L$ to ${[(d+1)e_i,...,(d+1)e_{n}]}$.
\end{enumerate}
If both $K$ and $L$ are primitive then so is $K+L$. By \cite[Proposition~7.5]{chenal_poincare_2026} the same holds for the convexity. We recursively defined the Viro triangulations $(V^n_d)_{n,d\geq 1}$ by the formulæ
\begin{equation*}
	\left\{ \begin{array}{ll}
		V^1_d \textnormal{ is the only primitive triangulation of }\P^1_d & \forall d\geq 1\\
		V^n_1=\P^n_1 & \forall n\geq 1 \\
		V^{n+1}_{d+1}=V^{n+1}_d+V^n_{d+1} & \forall n,d\geq 1.
	\end{array} \right.
\end{equation*}
This process can be generalised. Let $L=(L_d)_{d\geq 1}$ be a sequence of primitive triangulations of $(\P^n_d)_{d\geq 1}$, we denote by $(V_d(L))_{d\geq 1}$ the sequence of triangulations of $(\P^{n+1}_d)_{d\geq 1}$ defined by the formulæ
\begin{equation*}
	\left\{ \begin{array}{ll}
		V_1(L)=\P^n_1 &  \\
		V_{d+1}(L)=V_d(L)+L_{d+1} & \forall d\geq 1.
	\end{array} \right.
\end{equation*}

\begin{prop}\label{prp:Viro_gen_non_rho_uniform}
	Let $L=(L_d)_{d\geq 1}$ be a sequence of primitive triangulations of $(\P^{2n}_d)_{d\geq 1}$ and $d$ be an even integer. If $L_d$ has a $2n$-simplex that contains only interior vertices, then $V_D(L)$ is not $\rho$-uniform for all $D\geq d+1$.
\end{prop}
    
\begin{proof}
	Let $d$ be such an integer and $\sigma^{2n}$ be a $2n$-simplex of $L_d$ whose vertices are all interior. Using Proposition~\ref{prop:crit_rho_unif}, we see that $V_{d+1}(L)$ does not fail to be $\rho$-uniform at $\sigma^{2n}$ if and only if the two apices of the two $(2n+1)$-simplices adjacent to $\sigma^{2n}$ are congruent modulo 2. Indeed, since $\sigma^{2n}$ has an odd number of vertices, and $\sum_{v\in \sigma^{2n}}\rho_{\sigma^{2n}}(v)$ vanishes, $\rho_{\sigma^{2n}}$ can only be constant to $0$. Now the two apices of the two $(2n+1)$-simplices adjacent to $\sigma^{2n}$ are $2e_{2n+1}$ and $(d+1)e_{2n}$. Since $d$ is even, they are not congruent and $V_{d+1}(L)$ fails to be $\rho$-uniform at $\sigma^{2n}$. Since all the triangulations $(V_D(L))_{D\geq d+1}$ contain $V_{d+1}(L)$, they all fail to be $\rho$-uniform at $\sigma^{2n}$.
\end{proof}      

\paragraph{Itenberg-Viro Triangulations} In the last example, we always performed the subdivision of the prism $\P^{n+1}_{d+1}\cap\{0\leq x_{n+1} \leq 1\}$ in the same order. However, we can also reverse the order of the joins.

\begin{dfn}\label{dfn:K-L}
	If $K$ is a triangulation of $\P^{n+1}_d$ and $L$ is a triangulation of $\P^{n}_{d+1}$ we define $K+(-1)L$ as the unique triangulation of $\P^{n+1}_{d+1}$ for which:
	\begin{enumerate}
	    \item The restriction of $K+(-1)L$ to $\P^{n+1}_{d+1}\cap\{x_{n+1}\geq 1\}$ is the translation of $K$;
	    \item For all integers $0\leq i\leq n$, the restriction of the triangulation $K+(-1)L$ to the join of simplices ${[e_{n+1}+de_i,...,e_{n+1}+de_n]*[0,(d+1)e_1,...,(d+1)e_{i}]}$ is the join of the restriction of the translation of $K$ to ${[e_{n+1}+de_i,...,e_{n+1}+de_n]}$ with the restriction of $L$ to the face ${[0,(d+1)e_1,...,(d+1)e_{i}]}$.
	\end{enumerate}
\end{dfn}

As a direct consequence of the definition if both $K$ and $L$ are primitive then so is $K+(-1)L$. The same holds for the convexity by a small adaptation of \cite[Proposition~7.5]{chenal_poincare_2026}. Following \cite[\S 5]{itenberg_topology_1997}, we define the \emph{Itenberg-Viro triangulations} $(IV^n_d)_{n,d\geq 1}$ by the formulæ:
\begin{equation*}
	\left\{ \begin{array}{ll}
		IV^1_d \textnormal{ is the only primitive triangulation of }\P^1_d&\forall d\geq 1\\
		IV^n_1=\P^n_1 & \forall n\geq 1 \\
		IV^{n+1}_{d+1}=IV^{n+1}_d+(-1)^dIV^n_{d+1} & \forall n,d\geq 1.
	\end{array} \right.
\end{equation*}
An example of such a triangulation is illustrated in Figure~\ref{fig:trig_IV}. As in the case of Viro triangulations this process can be generalised from a sequence of primitive triangulations $(L_d)_{d\geq 1}$ of $(\P^n_d)_{d\geq 1}$ to a sequence $(IV_d(L))_{d\geq 1}$ of triangulations of $(\P^{n+1}_d)_{d\geq 1}$ by the formulæ:
\begin{equation*}
	\left\{ \begin{array}{ll}
		IV_1(L)=\P^{n+1}_1 &  \\
		IV_{d+1}(L)=IV_d(L)+(-1)^dL_{d+1} & \forall d\geq 1.
	\end{array} \right.
\end{equation*}
\begin{figure}[H]
	\centering
	\begin{subfigure}[t]{0.45\textwidth}
		\centering
		\begin{tikzpicture}
			\fill (0,0) circle (0.05);
			\fill (1,0) circle (0.05);
			\fill (2,0) circle (0.05);
			\fill (3,0) circle (0.05);
			\fill (4,0) circle (0.05);
			\fill (0,1) circle (0.05);
			\fill (1,1) circle (0.05);
			\fill (2,1) circle (0.05);
			\fill (3,1) circle (0.05);
			\fill (0,2) circle (0.05);
			\fill (1,2) circle (0.05);
			\fill (2,2) circle (0.05);
			\fill (0,3) circle (0.05);
			\fill (1,3) circle (0.05);
			\fill (0,4) circle (0.05);		
			\draw (0,0) -- (0,4) -- (4,0) -- cycle;
			\draw (0,1) -- (3,1);
			\draw (0,2) -- (2,2);
			\draw (0,3) -- (1,3);
			\draw[very thick] (0,2) -- (3,1);
			\draw[very thick] (3,1) -- (0,0);
			\draw[very thick] (1,3) -- (0,2);
			\draw (1,3) -- (1,2);
			\draw (1,2) -- (3,1);
			\draw (0,2) -- (2,1) -- (0,0);
			\draw (0,2) -- (1,1) -- (0,0);
			\draw (3,1) -- (1,0);
			\draw (3,1) -- (2,0);
			\draw (3,1) -- (3,0);
		\end{tikzpicture}
	\end{subfigure}
	\caption{The triangulation $IV^2_4$.}
	\label{fig:trig_IV}
\end{figure}
	
\begin{prop}\label{prop:rho_unif_IV}
	The triangulations $L_d$ are $\rho$-uniform for all $d\geq 1$ if and only if the triangulations $IV_d(L)$ are $\rho$-uniform for all $d\geq 1$.
\end{prop}

\begin{proof}
	Let $K$ be a triangulation of $\P^{n+2}_d$ and $0\leq q\leq n+2$ be an integer. We denote the set of $q$-simplices of $K$ having at least one vertex in the interior $K\setminus\partial K$ by $A_q(K)$. If $L$ is a triangulation of $\P^{n+1}_{d+1}$, we denote by $A^\pm_q(K;L)$ the set of $q$-simplices of $K+e_{n+1}$ that belong to $A_q(K\pm L)$ but not to $A_q(K)+e_{n+1}$. Let us denote by $\underline{K}$ the restriction of $K$ to the bottom face $[0;de_1;...;de_{n+1}]$ of $\P^{n+2}_d$. We have
	\begin{equation}\label{eq:A+}
		A_{n+1}(K+L)=(A_{n+1}(K)+e_{n+1}) \,\amalg\, A_{n+1}^+(K;L) \,\amalg\, \big([A_{n}(\underline{K})+e_{n+2}]*(d+1)e_{n+1}\big).
	\end{equation}
	Indeed, the simplices of $A_{n+1}(K+L)$ that are neither translates of simplices of $A_{n+1}(K)$ nor in $A_{n+1}^+(K;L)$ come from joins of simplices of $\underline{K}$ with simplices of $L$. By construction, the only such simplices not having all their vertices in $\partial(K+L)$ belong to the subdivision of the join ${\big(e_{n+2}+[0;de_1;...;de_{n+1}]\big) * (d+1)e_{n+1}}$. Since $(d+1)e_{n+1}$ belongs to $\partial(K+L)$, the remaining simplices of $A_{n+1}(K+L)$ are in the set ${\big([A_{n}(\underline{K})+e_{n+2}]*(d+1)e_{n+1}\big)}$. Likewise, we have
	\begin{equation}\label{eq:A-}
		A_{n+1}(K-L)=(A_{n+1}(K)+e_{n+2}) \,\amalg\, A_{n+1}^-(K;L) \,\amalg\, \big([A_{n}(\underline{K})+e_{n+2}]*0\big).
	\end{equation}
	Now let $L'$ be a triangulation of $\P^{n+1}_{d+2}$. We have
	\begin{equation}\label{eq:A-+}
		A^-_{n+1}(K+L;L')=(A_{n+1}(L)+e_{n+2})\,\amalg\,\big([A_{n}(L)+e_{n+2}]*2e_{n+2}\big),
	\end{equation}
	and similarly
	\begin{equation}\label{eq:A+-}
		A^+_{n+1}(K-L;L')=(A_{n+1}(L)+e_{n+2})\,\amalg\,\big([A_{n}(L)+e_{n+2}]*(2e_{n+2}+de_{n+1})\big).
	\end{equation}
	Therefore, combining (\ref{eq:A-}) and (\ref{eq:A-+}), we have
	\begin{align*}
		A_{n+1}(K+L-L') =& (A_{n+1}(K+L)+e_{n+2}) \, \amalg \, (A_{n+1}(L)+e_{n+2}) \\
		& \amalg\,\big([A_{n}(L)+e_{n+2}]*2e_{n+2}\big) \, \amalg\, \big([A_{n}(L)+e_{n+2}]*0\big),
	\end{align*}
	since $\underline{K+ L}$ equals $L$. And likewise with (\ref{eq:A+}) and (\ref{eq:A+-})
	\begin{align*}
			A_{n+1}(K-L+L') =& (A_{n+1}(K-L)+e_{n+2}) \, \amalg \, (A_{n+1}(L)+e_{n+2}) \\
			&  \amalg\,\big([A_{n}(L)+e_{n+2}]*(2e_{n+2}+de_{n+1})\big) \, \amalg\, \big([A_{n}(L)+e_{n+2}]*(d+2)e_{n+1}\big).
	\end{align*}
	If we assume $K+L$ to be $\rho$-uniform, $K+L-L'$ is $\rho$-uniform as well if and only if the simplices of $(A_{n+1}(L)+e_{n+2})$, $[A_{n}(L)+e_{n+2}]*2e_{n+2}$, and $[A_{n}(L)+e_{n+2}]*0$ satisfy the requirement of Proposition~\ref{prop:crit_rho_unif}.
	
	Let $\sigma^{n+1}$ be a simplex of the first set. It bounds the two simplices $\sigma^{n+1}*2e_{n+2}$ and $\sigma^{n+1}*0$. By Lemma~\ref{lem:rho}, we have
	\begin{equation*}
		\sum_{v\in \sigma^n}\rho_{\sigma^n}(v)v = 2e_{n+2} \;(\mod 2),
	\end{equation*}
	and therefore, $\rho_{\sigma^{n+1}}=0$.
	
	Now let $\sigma^{n}$ be a simplex of $A_{n}(L)$. It bounds two $(n+1)$-simplices of $L$ whose apices we denote by $v_1$ and $v_2$. Therefore, the simplex $(\sigma^{n}+e_{n+2})*2e_{n+2}$ bounds the two $(n+2)$-simplices $(\sigma^{n}*v_1+e_{n+2})*2e_{n+2}$ and $(\sigma^{n}*v_2+e_{n+2})*2e_{n+2}$ of $K+L-L'$. From Lemma~\ref{lem:rho}, we have:
	\begin{equation*}
		\rho_{(\sigma^{n}+e_{n+2})*2e_{n+2}}(2e_{n+2})+\sum_{v\in \sigma^n}\rho_{(\sigma^{n}+e_{n+2})*2e_{n+2}}(v+e_{n+2}) = v_2-v_1\;(\mod 2).
	\end{equation*}
	Since $v_2-v_1$ belongs the the hyperplane $\{x_{n+2}=0\}$, $\rho_{(\sigma^{n}+e_{n+2})*2e_{n+2}}(2e_{n+2})=0$ and 
	\begin{equation*}
		\rho_{(\sigma^{n}+e_{n+2})*2e_{n+2}}(v+e_{n+2})=\rho^L_{\sigma^n}(v),
	\end{equation*}
	for all $v\in \sigma^n$ ($\rho^L$ denotes the function $\rho$ of the triangulation $L$). Since $2e_{n+2}$ belongs to the boundary $\partial(K+L-L')$, $\rho_{(\sigma^{n}+e_{n+2})*2e_{n+2}}$ is constant on the interior vertices $(\sigma^{n}+e_{n+2})*2e_{n+2}$ if and only if $\rho^L_{\sigma^{n}}$ is constant on the interior vertices of $\sigma^n$ in $L$.
	
	The exact same reasoning applies to the simplices of the third set $[A_{n}(L)+e_{n+2}]*0$ and we find that $K+L-L'$ is $\rho$-uniform if and only if $K$ and $L$ are both $\rho$-uniform. By the same argument, we can also prove that $K-L+L'$ is $\rho$-uniform if and only if $K$ and $L$ are both $\rho$-uniform. The only difference in this case is the two apices with which we suspend the horizontal $(n+1)$-simplices of $L$. They are $2e_{n+2}+de_{n+1}$ and $(d+2)e_{n+1}$. Their difference is also $0$ modulo $2$ so $\rho_{\sigma^{n+1}}$ vanishes as in the case of $K+L-L'$. The proposition follows from the construction of the triangulations $(IV_d(L))_{d\geq 1}$.
\end{proof}

\begin{rem}\label{rem:rho_cones}
	In the last proof we showed that if $\sigma^{n+1}$ is a simplex of the triangulation $IV_d(L)$ of $\P^{n+2}_d$ written as a join $v*\sigma^{n}$ with $\sigma^{n}$ in a floor $L_k$, then $\rho_{v*\sigma^{n}}(v)=0$ and $\rho_{v*\sigma^{n}}(v')=\rho^{L_k}_{\sigma^{n}}(v')$ for all $v'\in \sigma^n$.
\end{rem}

\begin{dfn} \cite[Definition 6.5]{chenal_poincare_2026}. 
	Let $\ell(\X)$ be the supremum of the integers $\ell\geq 0$ for which $H_k(\X;\F_2)\rightarrow H_k(\R Y;\F_2)$ is surjective for all $k\leq \ell$.
\end{dfn}

\begin{prop}
	For all $n,d\geq 1$, for all $\varepsilon\in C^0(IV_d^n;\F_2)$, the number $\ell(\X)$ is at least equal to $\lfloor\frac{n-1}{2}\rfloor$.
\end{prop}

\begin{proof}
	This is \cite[Theorem~7.9]{chenal_poincare_2026} for Itenberg--Viro triangulations. The crucial points in the proof of this theorem were \cite[Proposition~7.8]{chenal_poincare_2026} i.e. the fact that the restriction of the Viro triangulation to the face opposite to $0$ is still a Viro triangulation and that $K+L$ contains a cone over $L$. However, the same is true for the Itenberg--Viro triangulation, hence the proposition holds. 
\end{proof}
    
\begin{thm}\label{thm:ItViro}
	For all $n\geq 2$, and all $d\geq 1$, the Itenberg--Viro triangulation $IV^n_d$ is simply integrable. In particular, the Harnack distribution
	\begin{equation*}
		h^n(x_1,...,x_n)\coloneqq \sum_{i<j}x_ix_j,
	\end{equation*}
	is a solution of (\ref{eq:simple_harnack}).
\end{thm}
    
\begin{proof}
    We will prove, by recursion on $n$, that $h^n$ is a solution of (\ref{eq:simple_harnack}) for all $d\geq 1$. Theorem~\ref{thm:quadratic_distrib} ensures that the property is true for $n=2$.
    
    Let us now assume that the property is true for some $n\geq2$. Let $1\leq d\leq D$ be two integers and $\sigma^n$ be a simplex of $IV^{n+1}_D$ containing at least one interior vertex and with non-empty intersection with the horizontal hyperplane ${\{x_{n+1}=D-d\}}$. Then, either $\sigma^n$ belongs to $IV^n_d$ or is a cone over an $(n-1)$-simplex of $IV^{n}_d$ with at least one interior vertex in $\P^n_d$. 
    
    If $d$ is even, the apex of such a cone is $(D-d-1)e_{n+1}$ or $(D-d+1)e_{n+1}$. And if $d$ is odd, it has to be $(D-d-1)e_{n+1} + (d+1)e_n$ or $(D-d+1)e_{n+1} + (d-1)e_n$. Either way we showed, cf. Remark~\ref{rem:rho_cones}, that, considering such a cone $v*\sigma^{n-1}$ with $\sigma^{n-1}\in IV^n_d$, we have
    \begin{equation}\label{eq:rho_recursion}
    	\rho_{v*\sigma^{n-1}}(v)=0 \quad\mathrm{and}\quad\rho_{v*\sigma^{n-1}}|_{\sigma^{n-1}}=\rho^{IV^n_d}_{\sigma^{n-1}}
    \end{equation}
    Moreover, (\ref{eq:rho_recursion}) with Proposition~\ref{prop:formula_hess} imply that if $e^1_1$ is the dual edge of $v*\sigma^{n-1}$ in $IV^{n+1}_D$ and $e^1_2$ is the dual edge of $\sigma^{n-1}$ in $IV_d^n$, we have
    \begin{equation*}
    	\D^2h^{n+1}|^{\F_2}_{e^1_1}=\D^2(h^{n+1}|_{IV^n_d})|^{\F_2}_{e^1_2}.
	\end{equation*}
	However, $h^{n+1}|_{IV^n_d}$ only differs from $h^n$ by an affine function. Hence, by recursion hypothesis and (\ref{eq:rho_recursion}), we have
	\begin{equation*}
		\D^2h^{n+1}|^{\F_2}_{e^1_1}=\D^2(h^{n+1}|_{IV^n_d})|^{\F_2}_{e^1_2}=\rho^{IV^n_d}|^{\F_2}_{e^1_2}=\rho|^{\F_2}_{e^1_1}.
	\end{equation*}
	Now, if $\sigma^n$ belongs to $IV^n_d$, the two $(n+1)$-simplices of $IV^{n+1}_D$ that are adjacent to $\sigma^n$ have their apices congruent modulo 2. It implies that $\rho_{\sigma^n}$ vanishes. Thus, if $e^1$ is the dual edge of $\sigma^n$ in $IV^{n+1}_D$, $\rho|_{e^1}^{\F_2}=0$. These two $(n+1)$-simplices are identical modulo 2, and the value of $h^{n+1}$ is the same on the two apices since they have congruent coordinates. Hence, $\D^2h^{n+1}|^{\F_2}_{e^1}$ vanishes as well. Therefore, we proved that $h^{n+1}$ is a solution of (\ref{eq:simple_harnack}).
\end{proof}
    
\begin{rem}
	The hypersurface $\R X_{h^n}$ of $\R\P^n_d$ has a maximal number of connected components but also intersects all the toric divisors of $\R\P^n_d$ along T-hypersurfaces that have a maximal number of connected components.  
\end{rem}

\bibliographystyle{alpha}
\bibliography{On_the_First_Page_of_the_Renaudineau_Shaw_Spectral_Sequence}
\end{document}